\input amstex
\documentstyle {amsppt}
\magnification=1200
\vsize=9.5truein
\hsize=6.5truein
\nopagenumbers
\nologo

\def\cc{\Cal C_{\ge Z}}

\topmatter

\title
The Chabauty space of closed subgroups of
\\ 
the three-dimensional Heisenberg group 
\endtitle

\author
Martin R. Bridson, Pierre de la Harpe, and Victor Kleptsyn
\endauthor

\leftheadtext
{Martin R Bridson, Pierre de la Harpe and Victor Kleptsyn}
\rightheadtext
{The space of subgroups of the Heisenberg group}

\abstract When equipped with
the natural topology first defined by Chabauty, 
the closed subgroups of a locally compact group $G$ 
form a compact space $\Cal C(G)$. 
We analyse the structure of $\Cal C(G)$ 
for some low-dimensional Lie groups, 
concentrating mostly on 
the 3-dimensional Heisenberg group $H$. 
We prove that $\Cal C(H)$ is a 6-dimensional space
that is path--connected but not locally connected.
The lattices in $H$ form a dense open subset
$\Cal L(H) \subset \Cal C(H)$ that is the disjoint union
of an infinite sequence of pairwise--homeomorphic
aspherical manifolds of dimension six, 
each a torus bundle over 
$(\bold S^3 \smallsetminus  T) \times \bold R$,
where $T$ denotes a trefoil knot. 
The complement of $\Cal L(H)$ in $\Cal C(H)$
is also described explicitly. 
The subspace of $\Cal C(H)$
consisting of subgroups that contain the centre $Z(H)$ 
is homeomorphic to the 4--sphere, 
and we prove  that this is a weak retract of $\Cal C(H)$.
\endabstract



\thanks
This work was funded by the Swiss National Science Foundation.
Bridson's work was also supported
in part by a Senior Fellowship from the EPSRC (UK) 
and a Royal Society Wolfson Research Merit Award,
while Kleptsyn's work was  supported
in part by the Russian scholarship RFBR 07-01-00017-a
and in part by the French--Russian scholarship 
RFBR CNRS-L-a 05-01-02801.
We thank each of these organisations. 
Bridson particularly thanks 
his colleagues at the University of Geneva and the EPFL (Lausanne)
for their hospitality throughout this project.
\endthanks

\address
Martin R Bridson, Section de Math\'ematiques,
Universit\'e de Gen\`eve, 
\newline
{\it on leave from}
Mathematics Department,
Imperial College
London, UK. \newline 
{\it Current address:}
Mathematical Institute,
24-29 St Giles', Oxford, OX1 3LB, UK.
\newline
E--mail: bridson\@maths.ox.ac.uk
\endaddress

\address
Pierre de la Harpe, Section de Math\'ematiques,
Universit\'e de Gen\`eve, C.P. 64,
\newline
CH--1211 Gen\`eve 4.
\newline
E--mail: Pierre.delaHarpe\@unige.ch
\endaddress

\address
Victor Kleptsyn, Section de Math\'ematiques,
Universit\'e de Gen\`eve, C.P. 64,
\newline
CH--1211 Gen\`eve 4.
\newline 
{\it Current address:}
UMR 6625 du CNRS, IRMAR, Campus de Beaulieu,
F--35042 Rennes cedex.
\newline
E--mail: Victor.Kleptsyn\@univ-rennes1.fr
\endaddress

\endtopmatter

\document

\head{\bf
1.~Introduction
}\endhead
\medskip

Let $G$ be a locally compact 
topological group.
We denote by $\Cal C (G)$ the {\it space of closed subgroups} of $G$
equipped with the {\it Chabauty topology};
this is a compact space. A basis of neighbourhoods
for a closed subgroup $C\in \Cal C (G)$ is formed by the subsets
$$
\Cal V_{K,U}(C) \, = \,
\{ D \in \Cal C (G) \mid
D \cap K \subset CU
\hskip.3cm \text{and} \hskip.3cm
C \cap K \subset DU \},
\tag{1.1}
$$
where $K\subseteq G$ is compact 
and $U$ is an open neighbourhood of the identity $e \in G$. 
We will consider the induced  topology on various subspaces
of $\Cal C(G)$
including the {\it space of lattices} $\Cal L (G)$,
the larger {\it space of discrete subgroups} $\Cal D (G)$, 
the {\it space of abelian closed subgroups} $\Cal A (G)$,
and the {\it space of normal closed subgroups} $\Cal N (G)$.

$\Cal C (G)$ is contained in the set of all
closed subspaces of the topological space $G$.
There is a natural topology on this larger set,  the study of which goes back in some form or another to
Hausdorff and Vietoris; we recall some of its basic properties
in Section~3. The Chabauty topology coincides with the restriction of this
topology to  $\Cal C (G)$.
\par

Chabauty's purpose in introducing this topology 
was to generalize a criterion of Mahler that allows one to prove 
that appropriate sets of lattices are relatively compact. 
Mahler's original work concerned lattices in $\bold R^n$,
while Chabauty's generalization deals with lattices
in a large class of locally compact groups  \cite{Chaba--50};
this result is known as the
{\it Mahler-Chabauty compactness  criterion.}
Two early references concerning the Chabauty topology are
\cite{MacSw--60} and \cite{Bourb--60, chapitre~8, \S~5}.
The Chabauty topology is also induced by a topology defined 
by Fell on spaces of closed subsets \cite{Fell--62}.
\par

The topology has been useful for the study 
of Fuchsian groups \cite{Harve--77, in particular Section 2}
and more generally 
of discrete subgroups of other semisimple groups 
\cite{GuiR\'e--06}, 
as well as hyperbolic manifolds, 
see \cite{Thurs--80, Chapter~9} and
\cite{CaEpG--87, in particular Section~I.3.1}; 
for Fuchsian and Kleinian groups,
convergence with respect to the Chabauty topology
is known as {\it geometric convergence.}
Observe that, whenever a space $X$ can be identified
with a family of closed subgroups of a locally compact group $G$,
as is  the case for Riemannian symmetric spaces
or Bruhat--Tits buildings, 
the closure of $X$ in $\Cal C (G)$
provides a natural compactification of $X$,
which is in many cases an efficient tool to study $X$
(this is  explicit in, for example, 
\cite{GuJiT--96, chap. IX}, \cite{GuiR\'e--06}, 
and \cite{BorJi--06, \S\S~I.17 and III.17--19}).
We also draw the reader's attention to  \cite{Ghys--07} where, 
by studying a certain flow on the space 
$\Cal L^{\operatorname{umod}} (\bold C)$ 
of coarea-1 lattices in $\bold C$, 
Ghys discovered fascinating phenomena linking 
dynamics, knot theory and arithmetic.
\par

When the group $G$ is discrete, the Chabauty topology
on $\Cal C (G)$ coincides with the restriction to $\Cal C (G)$
of the product topology on the compact space $2^G$ of all subsets of~$G$.
In particular, 
if $F_k$ is the  free group on $k$ generators,
$\Cal N (F_k)$ is the
{\it space of marked groups on~$k$ generators},
which has been the subject of much work since its appareance in
\cite{Gromo--81, final remarks}; see among others
\cite{Grigo--84},
\cite{Champ--00},
\cite{ChaGu--05},
\cite{CoGuP--07},
and \cite{CoGuP}.

\bigskip
\centerline{* * * * * * * * * * * * *}
\medskip

Our purpose in this work is to describe some  examples of spaces $\Cal C (G)$. 
A theme that the reader will observe quickly is that
for the most easily understood of low--dimensional Lie groups
the space $\Cal C(G)$ can be complicated
and yet beguilingly tractable.
\par

As a preliminary exercise, we consider in Section~2 
the group of orientation preserving affine transformations
of the real line
$$
\operatorname{Aff} \, = \,
\left\{ \left( \matrix e^s & t \\ 0 & 1 \endmatrix \right)
\, \Big\vert \,
s,t \in \bold R \right\} .
\tag{1.2}
$$
Up to isomorphism, this is the unique non--abelian
connected real Lie group of dimension~2.

\bigskip

\proclaim{1.1.~Proposition} The space $\Cal C (\operatorname{Aff})$
is a compact, contractible  space of dimension $2$.  
More precisely, it is the union of a closed 2-disc $\bold D$
and a compact interval $I$ attached by one of its end points
to a point on the boundary circle $\partial\bold D$. 

The disc $\bold D$ arises as a cone over the real projective line formed by
the subgroups isomorphic to $\bold R$; the cone point is the
trivial subgroup $\{e\}$.

The subspace $\Cal N (\operatorname{Aff})$ consisting
of the closed normal subgroups of $\operatorname{Aff}$
is the disjoint union of the point $\{e\}$ and the closed
interval $I$; in particular, $\Cal N (\operatorname{Aff})$ is disconnected.
\endproclaim

In Section~3, we recall some general facts on the Chabauty topology.
In Section~4 we revisit the following result of \cite{HubPo--79}, proving
a more detailed version than stated here (Theorem~4.6 versus Theorem~1.2).
This result is needed in our analysis of the 3--dimensional Heisenberg
group, which is our main interest.
\bigskip

{\bf Some notation.}
Given  locally compact groups $G, S,  \hdots, T$,
we denote by $\Cal C_{S,\hdots, T} (G)$ the subspace of $\Cal C (G)$
of closed subgroups isomorphic to one of $S,  \hdots, T$.
We write $\Sigma X$ to denote the suspension of a space $X$.

We write $\approx$ to indicate homeomorphism of topological spaces
and $\cong$ to indicate isomorphism of groups.
\bigskip

\proclaim{1.2.~Theorem (Hubbard and Pourezza)}
The space $\Cal C (\bold C)$ is homeomorphic to
the four--dimensional sphere $\bold  S^4$.

More precisely, there is a homeomorphism  
$\Cal C(\bold C) \longrightarrow \bold S^4 = \Sigma \bold S^3$ 
with the following properties: it identifies
$\Cal C_{\bold R} (\bold C)$ with a trefoil knot $T$
in the equator $\bold S^3\subset\bold S^4$; it
maps $\{e\}$ and $\bold C$ to the south  and north poles 
$\bold s, \bold n \in \bold S^4$, respectively; it
sends the space of lattices $\Cal L(\bold C)$ 
onto the complement of the (not locally--flat) $2$--sphere 
$\Sigma T\subset\Sigma\bold S^3$;
and it identifies $\Cal C_{\bold Z} (\bold C)$ and 
$\Cal C_{\bold R \oplus \bold Z} (\bold C)$
with the connected components of 
$\Sigma T \smallsetminus (T \cup \{\bold s, \, \bold n\} )$.
\endproclaim

We want to emphasize  how remarkable it is 
that $\Cal C (\bold R^2)$ is a manifold: 
even when one knows the homeomorphism type of the subspaces
corresponding to the different types of subgroups, 
it is by no means obvious that they assemble to form a sphere. 
For $n\ge 3$, the space $\Cal C (\bold R^n)$ is much wilder:
the space of lattices
$\Cal L (\bold R^n) \approx GL_n(\bold R)/GL_n(\bold Z)$
continues to provide a dense connected open set of dimension $n^2$
and the lower--dimensional subspaces
$\Cal C_{\bold R^a \oplus \bold Z^b}(\bold R^n)$
are easy to describe, but the manner in which they assemble is not.
\par

Another space which can be described 
in some detail is $\Cal C (SO(3))$.
This space is neither connected nor locally connected
(like $\Cal C (\bold Z)$).
In more detail,
the space of torus subgroups is the base space of the Hopf--like
fibration $SO(3) \longrightarrow \bold P^2$,
and if $T$ is such a torus subgroup then any neighbourhood of $T$
in $\Cal C (SO(3))$ contains finite cyclic groups
of arbitrarily large orders.
\par

Beginning in Section~5, we concentrate on the $3$--dimensional Heisenberg
group $H$. For the most part we shall work with the following model of this
group:
$$
\aligned
&\text{as a set} \hskip.2cm
H \, = \, \bold C \times \bold R,
\\
&\text{the product is given by} \hskip.2cm
(z , t)(z' , t') \, = \,
\Big(z+z' , t+t' + \frac{1}{2}\operatorname{Im}(\overline{z}z')\Big) ,
\\
&\text{and there is a homomorphism} \hskip.2cm
p \, : \, (z,t) \, \longmapsto \, z 
\hskip.2cm \text{from} \hskip.2cm H 
\hskip.2cm \text{onto} \hskip.2cm \bold C .
\endaligned
\tag{1.3}
$$
This  model (1.3) for the Heisenberg group
is isomorphic to the standard matrix model
via the isomorphism
$$
\bold C \times \bold R \, \ni \,
(x+iy,t)
\hskip.2cm \longmapsto \hskip.2cm
\exp
\left( \matrix 0 & x & t \\ 0 & 0 & y \\ 0 & 0 & 0 \endmatrix \right)
\, = \,
\left( \matrix
1 & x & t + \frac{1}{2}xy \\ 0 & 1 & y \\ 0 & 0 & 1
\endmatrix \right)
\, \in \,
\left( \matrix
1 & \bold R & \bold R \\ 0 & 1 & \bold R \\ 0 & 0 & 1
\endmatrix \right).
$$
$H$ is also isomorphic to the semidirect product
$\bold R \ltimes_{u}  \bold R^2$ 
associated to the unipotent action
$u : \bold R \times \bold R^2 \longrightarrow \bold R^2$
defined by
$$
\left( x \hskip.1cm , \hskip.1cm 
\left( \matrix t \\ y \endmatrix \right) 
\right)
\hskip.2cm \longmapsto \hskip.2cm
\left( \matrix 1 & x \\ 0 & 1 \endmatrix \right)
\left( \matrix t \\ y \endmatrix \right)
=
\left( \matrix t+xy \\ y \endmatrix \right) .
\tag{1.4}
$$

\bigskip

Our examination of the structure of $\Cal C (H)$ 
spans Sections 5 to 8.
We summarize the main results of these Sections in the following
two theorems, the first of which provides an overview of the
global structure of $\Cal C (H)$ and the second of which gives more
information about the strata corresponding to the various
types of subgroups as well as the manner in which 
$\operatorname{Aut}(H)$ acts on them. 
\par

We define subspaces $\Cal L_n(H) \subset \Cal L(H)$
by declaring that a lattice $\Lambda$ is in $\Cal L_n(H)$ 
if $[\Lambda,\Lambda]$ has index $n$ in $\Lambda \cap Z(H)$.
The subset  $\Cal L_{\infty}(H)\subset 
\Cal C (H)$ is defined to consist of the subgroups  
$p^{-1}(L)$ with $L \in \Cal L (\bold C)$.
And $\Cal L_{!!} (H) = \Cal L (H) \cup \Cal L_{\infty}(H)$.
We define  $\cc (H)$ as the space of closed subgroups of $H$
which contain the centre $Z(H)$;
observe that $\Cal L_{\infty}(H) \subset \cc(H)$.
\par

We write  $\bold P^1$ to denote
the real projective line and
$\bold P^2$ for the real projective plane.
\par
The real projective line $\bold P^1$ is of course
homeomorphic to the circle $\bold S^1$.
However, in this work, we try to mark the difference
between homeomorphism with a circle
related to the choice of a {\it unit vector} in a plane
and homeomorphism with a circle
related to the choice of a {\it line} in a plane;
see for example the covering
$\bold S^1 \approx \widehat{\Cal C}_{\bold R^2}(H)
\, \longrightarrow \,
\Cal C_{\bold R^2}(H) \approx \bold P^1$
of (6.5) below.

\bigskip

\proclaim{1.3.~Theorem}  The compact space $\Cal C (H)$ 
is arc--connected but not locally connected.
It can be expressed as the union of 
the following three subspaces.
\roster
\item"(i)"
$\Cal L(H)$, which is open and dense in $\Cal C (H)$;
this has countably many connected components $\Cal L_n(H)$, 
each of which is homeomorphic to a fixed aspherical $6$--manifold that is a 
$2$--torus bundle over $\Cal L (\bold C) \approx GL_2(\bold R) / GL_2(\bold Z)$.
\item"(ii)"
$\Cal A(H)$, which is homeomorphic to  the space obtained from
$\bold S^4\times\bold P^1$ by fixing a tame arc $I\subset\bold S^4$
and collapsing each of the circles 
$\{\{i\} \times \bold P^1 : i\in I\}$
to a point.  
\item"(iii)"
$\cc(H)$, from which there
 is a natural homeomorphism to  $\bold S^4$;
the complement of $\Cal L_{\infty}(H)$ in $\cc (H)$ is a 
$2$--sphere $\Sigma^2\subset \bold S^4$
(which fails to be locally flat at  two points).
\endroster
The union $\Cal A(H) \cup \cc(H)$ is the complement of $\Cal L (H)$ in $\Cal C (H)$.
The intersection
$$
\Cal A(H) \cap \cc(H) =
\{ C \in \cc (H) \mid p(C) \subset
\Cal C_{\{0\},\bold Z,\bold R}(\bold C) \}
$$ 
is  a closed $2$--disc in $\Sigma^2$. 
The space $\Cal L_{!!} (H) = \Cal L (H) \cup \Cal L_{\infty}(H)$
is precisely $\{ C \in \Cal C(H) \mid p(C) \in \Cal L (\bold C) \}$.

\smallskip

$\cc (H)$ is a weak retract of $\Cal C (H)$:
there exists a continuous map 
$f : \Cal C (H) \longrightarrow \bold S^4$, 
constant on $\Cal A(H)$, such that 
$f \circ j \simeq {\operatorname{id}}_{\bold S^4}$, where
$j : \bold S^4 \longrightarrow \cc(H)$ is the homeomorphism of (iii),
and where $\simeq$ denotes  homotopy equivalence.
In particular, $\pi_4(\Cal C (H))$ surjects onto~$\bold Z$.

\smallskip

The subspace $\Cal N (H)$ of normal closed subgroups of $H$
is the union of $\cc (H)$ (which 
is homeomorphic to $\Cal C (\bold C)
\approx \bold S^4$)
and  the interval $\{ C \in \Cal C (H) \mid C \subset Z(H) \}$,
attached to the sphere $\cc (H)$ by one of its endpoints.
\endproclaim 
 
For our second compendium of results concerning $\Cal C (H)$ we need
the following notation.  
We denote by $\bold K$ the Klein bottle
(i.e. the total space of the  non--trivial 
$\bold P^1$ bundle over $\bold P^1$). 
We write  
$p_*^{-1}\left( \Cal C_{\bold R \oplus \bold Z}(\bold C) \right)$
for the  subspace of $\Cal C (H)$ consisting of 
closed subgroups $C$ of $H$ 
with $p(C)$ a subgroup of $\bold C$
isomorphic to $\bold R \oplus \bold Z$;
observe that
$p_*^{-1}\left( \Cal C_{\bold R \oplus \bold Z}(\bold C) \right)
\subset \cc(H)$.

\bigskip

\proclaim{1.4.~Theorem} 
The spaces $\Cal L_n(H)$ are homeomorphic  
to a common aspherical homogeneous space, namely the quotient
of the  $6$--dimensional automorphism group
$\operatorname{Aut}(H)  \cong
\bold R^2  \rtimes  GL_2(\bold R)$
by the discrete subgroup $\bold Z^2\rtimes GL_2(\bold Z)$.
\par

The frontier of $\Cal L_n(H)$, which is independent of $n$,
consists of the following subspaces:
\roster
\item"(i)"
the trivial group $\{e\}$;
\item"(ii)"
$\Cal C_{\bold R}(H) \approx \bold P^2$;  
\item"(iii)"
$\Cal C_{\bold Z}(H) \approx \bold P^2 \times ]0,\infty[$;
\item"(iv)"
$\Cal C_{\bold R^2}(H) \approx \bold P^1$;
\item"(v)"
$\Cal C_{\bold R \oplus \bold Z}(H) \approx \bold K \times ]0,\infty[$,
which is a $(\bold P^1 \times ]0,\infty[)$--bundle over $\bold P^1$;
\item"(vi)"
$\Cal C_{\bold Z^2}(H)$,
which is a $(\bold S^4 \smallsetminus \Sigma^2)$--bundle over $\bold P^1$;
\item"(vii)"
$p_*^{-1}\left( \Cal C_{\bold R \oplus \bold Z}(\bold C) \right)$;
\item"(viii)"
the full group $H$.
\endroster
In particular, the frontier of $\Cal L_n(H)$, which is independent of $n$, 
is the union of  $\Cal A(H)$
and the complement $\Sigma^2$ of $\Cal L_{\infty}(H)$ in $\cc (H)$;
the part $\Cal A (H)$ is itself  the union of the subspaces {\rm (i)} to {\rm (vi)},
and $\Sigma^2 \smallsetminus (\Sigma^2 \cap \Cal A (H))$ 
is itself the union of the subspaces {\rm (vii)} and {\rm (viii)}.
The frontier of $\bigcup_{n=1}^{\infty} \Cal L_n(H)$ further contains
\roster
\item"(ix)"  
$\Cal L_{\infty}(H)$.
\endroster
Each of these spaces, except~{\rm (vi)}, consists of finitely many
$\operatorname{Aut}(H)$--orbits.
\endproclaim

The action of $\operatorname{Aut}(H)$ on $\Cal C_{\bold Z^2}(H)$,
which has uncountably many orbits and which is minimal,
is described in the proof of Proposition~6.1,
in terms of the standard action of $SL_2(\bold Z)$ on $\bold P^1$.

Observe that, as $\Cal L (H)$ is open dense, 
the spaces (i) to (ix) of Theorem~1.4
together with the spaces $\Cal L_n (H)$ for $n \ge 1$
constitute a partition of $\Cal C (H)$.

\par

The subspaces of $\Cal C (H)$ are described more precisely in 
Section~6 (for $\Cal A (H)$ and $\cc (H)$), 
Section~7 (for $\Cal L_n (H)$),
and Section~8 (for $\Cal L_{!!} (H)$ and $\Cal C (H) \smallsetminus \Cal A (H)$).

\bigskip
\head{\bf
2.~First examples, including the affine group $\operatorname{Aff}$
}\endhead
\medskip

\par

The space $\Cal C (G)$ has a straightforward description
when $G$ is a 1-dimensional Lie group:
$\Cal C (\bold R)$  is homeomorphic to the closed
interval $[0,\infty]$,
with $\lambda \in ]0,\infty[$ corresponding
to the subgroup $\bold Z \lambda^{-1}$; and $\Cal C(\bold R / \bold Z)$
is homeomorphic to $\Cal C(\bold Z)$, which  is homeomorphic to
$\{1, \frac{1}{2}, \frac{1}{3}, \cdots, 0 \}\subset [0,1]$,
with $\frac{1}{n}$ corresponding
to the subgroup of index $n$ and $0$ corresponding to 
 $\{0\} = \bigcap_{n=1}^{\infty} \bold Z n$. 
 (In passing, we note that 
 $\Cal C(\bold Z)$ is also homeomorphic to $\Cal C(\bold Z_p)$
 for an arbitrary prime $p$, where $\bold Z_p$ denotes the
 additive group of $p$--adic integers.)
\par

The direct observation $\Cal C(\bold R / \bold Z)\approx
\Cal C(\bold Z)$ illustrates the more general fact that
 if $G$ is a locally compact abelian group
with Pontryagin dual $\widehat G$, then there is a
homeomorphism
$\Cal C (G) \longrightarrow \Cal C (\widehat G)$
associating to each closed subgroup $C \subset G$
its orthogonal
$C^{\perp} = \{\chi \in \widehat G
\hskip.1cm \vert \hskip.1cm \chi(c) = 1
\hskip.1cm \text{for all} \hskip.1cm
c \in C \}$.
In the case $G=\bold R$, if we choose an isomorphism
of $\widehat{\bold R}$ with $\bold R$,
the above homeomorphism defines 
an involution of the interval $[0,\infty]$
that exchanges the two endpoints.

\medskip
\head{\bf
Closed subgroups of the affine group $\operatorname{Aff}$
}\endhead

We now turn to  an example that is only marginally more involved,
namely the affine group  of the real line,
$\operatorname{Aff}$ as described in (1.2).
The commutator subgroup $[\operatorname{Aff},\operatorname{Aff}]$
is the translation subgroup, described by the equation $s=0$.
The following properties are straightforward to check.
\roster
\item"---"
The Lie algebra $\eufm{aff}$ of $\operatorname{Aff}$ 
consists of the set of matrices
$\left( \matrix x & y \\ 0 & 0 \endmatrix \right)$
with $x,y \in \bold R$.
\item"---"
The exponential mapping $\eufm{aff} \longrightarrow \operatorname{Aff}$
is given by
$\exp \left( \matrix x & y \\ 0 & 0 \endmatrix \right)
\, = \,
\left( \matrix e^x & y \, \frac{e^x-1}{x} \\ 0 & 1 \endmatrix \right)$;
it is a diffeomorphism.
\item"---"
Any one--dimensional subspace of $\eufm{aff} $ distinct from
$[\eufm{aff} ,\eufm{aff} ]$ is conjugate under the adjoint representation
to the subspace with equation $y=0$; more precisely, when $x_0 \ne 0$,
we have
$$
\left( \matrix 1 & y_0/x_0 \\ 0 & 1 \endmatrix \right)
\left( \matrix sx_0 & sy_0 \\ 0 & 0 \endmatrix \right)
\left( \matrix 1 & -y_0/x_0 \\ 0 & 1 \endmatrix \right)
\, = \,
\left( \matrix sx_0 & 0 \\ 0 & 0 \endmatrix \right)
$$
for any $s \in \bold R$.
\item"---"
Any $g \in \operatorname{Aff}$, $g \ne e$,
lies in a unique one--parameter subgroup of $\operatorname{Aff}$
which is its centralizer $Z_{\operatorname{Aff}}(g)$.
\item"---"
Every non--abelian closed subgroup of $\operatorname{Aff}$
contains $[\operatorname{Aff},\operatorname{Aff}]$.
\endroster
Consequently, closed subgroups of $\operatorname{Aff}$ can be listed as
follows:
\roster
\item"(i)"
the trivial group $\{e\}$;
\item"(ii)"
the infinite cyclic subgroups;
\item"(iii)"
the one--parameter subgroups;
\item"(iv)"
the groups generated by $[\operatorname{Aff},\operatorname{Aff}]$
and one element $g=
\left( \matrix \lambda_g & t \\ 0 & 1 \endmatrix \right)$ with $\lambda_g>0$;
\item"(v)"
the group $\operatorname{Aff}$ itself.
\endroster
For the Chabauty topology, abelian closed subgroups,
which are the subgroups of types (i), (ii), and (iii),
and which are also the unimodular
\footnote{
Compare with Theorem 1.i in \cite{Bourb--60, chapitre~8, \S~5}:
if $G$ is a locally compact group, the subset
of unimodular closed subgroups is closed in  $\Cal C (G)$.
}
closed subgroups,
constitute a closed subspace $\bold D$ of $\Cal C (\operatorname{Aff})$
homeomorphic to the 2-disc; it is natural to regard it as a cone
$$
\big( \bold P^1 \times [0,\infty] \big)
\hskip.1cm / \hskip.1cm
\big( (x,0) \sim (y,0) \big)
$$
over the real projective line.
The vertex corresponds to the trivial group $\{e\}$,
the points in $\bold P^1 \times \{\infty\}$
to the one parameter subgroups isomorphic to $\bold R$,
and the other points to the infinite cyclic subgroups
(each infinite cyclic subgroup being contained
in a unique one parameter subgroup).
The parameter $\lambda_g$ defines a homeomorphism
from the subspace formed by subgroups of type (iv) 
to the open interval $]0,\infty[$. The closure of this
subspace is a compact interval
with endpoints $\operatorname{Aff}$ and 
$[\operatorname{Aff},\operatorname{Aff}]$. It intersects
$\bold D$ in a single point, namely 
the point in  $\bold P^1 \times \{\infty\}$
corresponding to $[\operatorname{Aff},\operatorname{Aff}]$.

The proof of Proposition~1.1 is complete.

\head{\bf
3.~A reminder concerning the Chabauty topology
}\endhead
\medskip

Aspects of
Proposition 1.1 and Theorems 1.2 to 1.4 
illustrate properties that hold in more general groups.
We record some of these in Proposition 3.4, 
for future reference. Two preliminary lemmas are required.

\bigskip

\proclaim{3.1.~Lemma}
Let $G$ be a topological group.
Let $K$ be a compact subset 
and let $V_1, \hdots, V_n$ be open subsets of $G$
such that $K \subset \bigcup_{j=1}^n V_j$.
Then there exists a neighbourhod $U$ of $e$ such that,
for all $x \in K$, there exists $j \in \{1,\hdots,n\}$
with $Ux \subset V_j$.
\endproclaim

\demo{Remark} This is strongly reminiscent of (and inspired by)
the so--called {\it Lebesgue number lemma}, according to which
if $\Cal V$ is an open covering of a compact metric space $X$,
there exists $\epsilon > 0$ such that every subset of $X$
of diameter less than $\epsilon$ is contained in 
an element~ $V$ of $\Cal V$. See e.g. \cite{Munkr--75}, \S~3.7.
\enddemo

\demo{Proof} For each $y \in K$ choose $j=j(y)$ such that $y \in V_j$. Then
$V_jy^{-1}$ is a neighbourhood of $e$ and 
there exist neighbourhoods $U_y, U'_y$ of $e$
such that $U_yU'_y \subset V_jy^{-1}$.
Set $W_y = U'_y y$; this is a neighbourhood of $y$ and
 $U_yW_y \subset V_j$. 
\par

From the open cover $(W_y)_{y \in K}$ we extract a finite
subcover indexed by $y_1,\hdots,y_N \in K$, say. 
Set $U = U_{y_1} \cap \cdots \cap U_{y_N}$. 
Then, for any $x \in K \subset \bigcup_{k=1}^N W_{y_k}$, 
we have $Ux\subset U_{y_k}W_{y_k} \subset V_j$ 
for some $k\in\{1,\dots,N\}$ and $j=j(y_k)$.
\hfill $\square$
\enddemo
\bigskip

Recall  that a closed subgroup $C$ of $G$
is {\it cocompact} if there exists a compact subset $K$ of $G$
such that $G = CK$. A {\it lattice} in a locally compact group $G$ 
is a discrete subgroup $\Lambda$ such that there exists
a $G$--invariant finite probability measure on $G / \Lambda$.

For subsets $A,B,\hdots$ and elements $g,h,\hdots$ of a group $G$,
we denote by $\langle A,B,\hdots,g,h,\hdots \rangle$ the subgroup of $G$
generated by $A \cup B \cup \hdots \cup \{g,h,\hdots\}$
(this subgroup need not be closed). 

\bigskip

\proclaim{3.2.~Lemma}
Let $G$ be a topological group 
that is compactly generated,
let $K$ be a compact generating set of $G$
such that $e \in K$ and $K^{-1} = K$,
and let $W$ be a relatively compact nonempty 
open subset of $G$.
Let $c_1,\hdots,c_{n} \in G$ be such that
$$
\overline{W} K \, \subset \, \bigcup_{j=1}^{n} c_j W .
\tag{3.1}
$$
Then there exists a symmetric neighbourhood $U$ 
of $e$ in $G$ such that,
for any $d_1,\hdots,d_{n} \in G$
with $d_j \in Uc_j$ for all $j =1,\hdots,n$, 
we have
$$
\langle d_1, \hdots, d_{n} \rangle W \, = \, G .
$$
In particular, the closed subgroup
$\overline{\langle d_1, \hdots, d_{n} \rangle}$ is cocompact.
\endproclaim

\demo{Proof} Let us first check that
Condition (3.1) is open in $c_1,\hdots,c_n$.
\par

By the previous lemma, there exists 
a neighbourhood $U$ of $e$ in $G$
such that, for any $x \in \overline{W}K$,
there exists $j \in \{1,\hdots,n\}$
with $Ux \subset c_jW$.
There is no loss of generality if we assume that,
moreover, $U^{-1} = U$.
For $j \in \{1,\hdots,n\}$,
consider any $d_j \in Uc_j$.
Let $u_j \in U$ be such that $d_j = u_jc_j$.
For $x \in \overline{W}K$ with $Ux \subset c_jW$,
we have $u_j^{-1}x \in c_jW$,
namely $x \in d_jW$.
Hence
$$
\overline{W} K \, \subset \, \bigcup_{j=1}^{n} d_j W .
\tag{3.2}
$$

\medskip

Set $F = \{d_1,\hdots,d_{n}\}$;
now (3.2) reads
$FW \supset \overline{W}K$.
By induction on $m$, it follows  that
$$
F^{m}W \, \supset \,
F^{m-1} \overline{W} K \, \supset \,
F^{m-1} W K \, \supset \,
\overline{W} K^{m}
$$
for any $m \ge 2$. Consequently
$\langle F \rangle W \supset \overline{W} \langle K \rangle$;
note that, at this point, we have used the hypothesis
$e \in K$ and $K^{-1} = K$.
Hence $\langle F \rangle W \supset \overline{W}G = G$,
and in particular $\overline{ \langle F \rangle}$
is cocompact.
\hfill $\square$
\enddemo

\bigskip

\proclaim{3.3.~Remark} Upon replacing $U$ by
$U' = \bigcap_{k=1}^n c_k^{-1}Uc_k$,
we can add the following conclusion to Lemma~3.2:
for any $d_1,\hdots,d_n \in G$ if $d_j \in c_j U'$
for $j=1,\hdots,n$, 
the closed subgroup $\overline{\langle d_1,\hdots,d_n \rangle}$
is cocompact.
\endproclaim

For a Hausdorff topological group,
the following two properties are clearly equivalent: (i)
every neighbourhood of the identity 
contains a non--trivial subgroup; (ii)
every neighbourhood of the identity 
contains a non--trivial closed subgroup.
(Indeed, for any neighborhood $U$ of the identity,
there exists a neighborhood $V$ of the identity
such that $\overline{V} \subset U$;
for any subgroup $S$  inside $V$, 
the closure $\overline{S}$ is a closed subgroup inside $U$.)
Recall that, by definition, a group
has {\it no small subgroup}, or is {\it NSS},
if these properties do {\it not} hold.

\bigskip

\proclaim{3.4.~Proposition} Let $G$ be a locally compact group.
\par

(i) The space $\Cal A (G)$ of closed abelian  subgroups of $G$
is closed in $\Cal C (G)$.

(ii) The space $\Cal K (G)$ of cocompact closed subgroups of $G$
is open in $\Cal C(G)$ if and only if $G$ is compactly generated.
\par

(iii) If $G$ is NSS,  the space $\Cal D (G)$ of discrete subgroups of $G$ 
is open in $\Cal C (G)$.
\par

\noindent
(It follows that,  if $G$ is compactly generated,
the space of cocompact lattices of $G$
is open in $\Cal C (G)$.)
\endproclaim

\demo{Remarks} 
Claim (ii) is essentially a result of Stieglitz and Oler \cite{Oler--73}.
Compare (ii) with Parts (i) and (ii) in 
\cite{Bourb--60}, chapitre~8, \S~5, no~3, th\'eo\-r\`eme~1:
if $\Cal U (G)$ denotes the space of 
unimodular closed subgroups of $G$,
then $\Cal U (G)$ is closed in $\Cal C (G)$.
If, moreover, $G$ is compactly generated,
then $\Cal K (G) \cap \Cal U (G)$ is open in $\Cal U (G)$.
\par
Compare (iii) with Parts (i) and (iii) in
\cite{Bourb--60}, chapitre~8, \S~5,  no~4, th\'eor\`eme~2:
if $G$ is NSS,
then $\Cal D (G)$ is locally closed in $\Cal C (G)$.
If, moreover, $G$ is  compactly generated,
then $\Cal D (G) \cap \Cal K (G)$ 
is locally closed in $\Cal C (G)$.
\enddemo

\demo{Proof} 
(i) Towards showing that the complement of $\Cal A (G)$ in $\Cal C (G)$
is open, we fix $C \in \Cal C (G) \smallsetminus \Cal A (G)$
and $x,y \in C$ such that $xy \ne yx$.
There exist neighbourhoods $U_x,U_y$ of $e$ in $G$ such that,
for all $x' \in xU_x$ and $y' \in yU_y$, we have $x'y' \ne y'x'$.
Set $U = U_x \cap U_y$.
\par

   For any $D \in \Cal V_{\{x,y\},U^{-1}}(C)$
we have $\{x,y\} \subset DU^{-1}$, see (1.1);
in other words, there exist $x',y' \in D$
such that $x' \in xU$ and $y' \in yU$.
In particular, $D$ is not abelian.

\medskip 

(ii) 
Suppose first that $\Cal K (G)$ is open in $\Cal C (G)$.
Since $G$ is clearly in $\Cal K (G)$,
there exist a compact subset $K$ and a nonempty open subset $U$ in $G$
such that
$$
\Cal V_{K,U}(G) \, = \, \left\{ D \in \Cal C (G) \mid
K \subset DU \right\} \, \subset \, \Cal K (G) .
$$
The closed subgroup $\overline{ \langle K \rangle}$
generated by $K$ is in $\Cal V_{K,U}(G)$,
{\it a fortiori} in $\Cal K (G)$.
[Note that $\langle K \rangle$
need not be closed.
This justifies the introduction of
$\langle K,V \rangle$ below.]
Let $V$ be any relatively compact open neighbourhood of $e$ in $G$.
On the one hand, the subgroup $\langle K,V \rangle$
generated by $K \cup V$ is open in $G$,
and therefore also closed;
on the other hand, we have 
$\overline{ \langle K \rangle} \subset \langle K,V \rangle$,
and therefore $\langle K,V \rangle \subset \Cal K (G)$.
Thus, there exists a compact subset $L$ of $G$ such that 
$\langle K,V \rangle L = G$.
It follows that the compact subset
$K \cup \overline{V} \cup L$ generates $G$;
in particular, $G$ is compactly generated.
\medskip

Suppose now that $G$ is compactly generated,
say by some compact subset $K$.
Without loss of generality, we assume that
$e \in K$ and $K^{-1} = K$.
Let $V$ be a relatively compact open neighbourhood of $e$.
\par

We fix $C \in \Cal K (G)$ and choose a compact subset $L$ of $G$
such that $CL = G$. Observe that $LV$ is open in $G$;
set $M = L\overline{V}$, which is also a compact subset of $G$.
Since we have
$$
MK \, \subset G \, = \, CLV \, = \, \bigcup_{c \in C} cLV 
$$
with $MK$ compact and the $cLV$'s  open,
there exists a finite family $c_1,\hdots,c_{n}$ in $C$
such that 
$MK \subset \bigcup_{j=1}^{n} c_j LV$.
\par

By Lemma~3.2 and Remark~3.3, there exists 
an open neighbourhood $U$ of $e$
such that, whenever $d_1,\hdots,d_n$ in $G$ are such that
$d_j \in c_j U$ for $j=1,\hdots,n$,
any {\it closed} subgroup  of $G$ containing $d_1,\hdots,d_n$
is cocompact.
Set $F = \{c_1,\hdots,c_n\}$.
Then $c_1,\hdots,c_n \in DU^{-1}$
for any closed subgroup
$D \in \Cal V_{F,U^{-1}}(C)$, see (1.1);
in other words, any $D \in \Cal V_{F,U^{-1}}(C)$
contains elements $d_1,\hdots,d_n$ 
such that $d_j \in c_j U$ for all $j$.
We have shown that the neighbourhood $\Cal V_{F,U^{-1}}(C)$
of $C$ is contained in $\Cal K (G)$.
This ends the proof of (ii).

\medskip

(iii) Since locally compact NSS groups are metrisable
(see for example \cite{Kapla--71}, Chapter II, Theorem 2),
there exists a left--invariant distance function
$\delta : G \times G \longrightarrow \bold R_+$
defining the topology of $G$.
For a closed subgroup $C$ of $G$, 
define the minimum distance
$$
\operatorname{min}_{\delta}(C)
\, = \,
\min \{ \delta(e,c) \mid c \in C, c \ne e \}
\, \in \,
\bold R_+^* \cup \{\infty\} .
$$
Following the standard convention, we set
$\operatorname{min}_{\delta}(\{e\}) = \infty$.
It is easy to check that the mapping
$\Cal C (G) \longrightarrow \bold R_+^* \cup \{\infty\}$
given by
$C \longmapsto \operatorname{min}_{\delta}(C)$
is continuous, and it follows that
$\Cal D (G)$ is open in $\Cal C (G)$.
\hfill $\square$
\enddemo

\bigskip

\proclaim{3.5.~Remarks}
(i) The space $\Cal L (G)$ of lattices of $G$
need not be open in $\Cal C (G)$.
\par
(ii) There are classes of groups $G$ in which all lattices are cocompact,
and $\Cal L (G)$ is open in $\Cal C (G)$ for such groups.
This applies for example to soluble Lie groups 
with countably many connected components 
(Theorem~3.1 in  \cite{Raghu--72}).
\par
(iii) A locally compact group which is not unimodular 
does not contain any lattice
(Remark~1.9 in  \cite{Raghu--72}).
In particular, for a solvable Lie group $G$,
the space $\Cal L (G)$ is empty \lq\lq in most cases\rq\rq .
\endproclaim

An example to which Remark 3.5.i applies is the group $PSL_2(\bold R)$   
of fractional linear transformations of the Poincar\'e half--plane
$\Cal H = \left\{ z \in \bold C \mid \operatorname{Im}(z) > 0 \right\}$.
For each real number $s \ge 2$, 
let $\Gamma_s$ denote the subgroup of $PSL_2(\bold R)$
generated by the transformations
$$
\left[ \matrix 1 & s \\ 0 & 1 \endmatrix \right]
\, : \, z \longmapsto z+s
\hskip.5cm \text{and} \hskip.5cm
\left[ \matrix 1 & 0 \\ s & 1 \endmatrix \right]
\, : \, z \longmapsto \frac{z}{sz+1} .
$$
It is standard that $\Gamma_s$ 
is a non--abelian free group of rank $2$ for any $s \ge 2$,
that $\Gamma_2$ is a lattice in $PSL_2(\bold R)$,
and that $\operatorname{Vol}(\Gamma_s \backslash \Cal H) = \infty$ for
$s>2$ (so then $\Gamma_s$ is not a lattice);
see for example Exercise II.33 in \cite{Harpe--00}.
\par

Thus, in $\Cal C (PSL_2(\bold R))$,
any neighbourhood of the lattice $\Gamma_2$
contains for $\epsilon$ small enough
a non--lattice $\Gamma_{2+\epsilon}$.

\bigskip

Let us assume that the locally compact group $G$ 
is  metrisable, and that its topology is defined 
by a left--invariant distance function $\delta$ relative to which
closed balls are compact.
Then  $\Cal C(G)$ is also metrisable,
and there is a basis of neighbourhoods
of any point $C$ in $\Cal C (G)$ which consists of the subsets
$$
\Cal V_{R,\epsilon}(C) \, = \,
\left\{ 
D \in \Cal C (G) \hskip.2cm \Bigg\vert \hskip.2cm
\aligned
\delta(x,D) < \epsilon \hskip.2cm &\text{for all} \hskip.2cm
x \in C \cap \overline{B}_R \hskip.2cm \text{and} \hskip.2cm
\\
\delta(y,C) < \epsilon \hskip.2cm &\text{for all} \hskip.2cm
y \in D \cap \overline{B}_R
\endaligned 
\right\} .
\tag{3.3}
$$
Here, $B_R$ denotes for $R > 0$ the open ball
$\{g \in G \mid \delta(g,e) < R\}$.
\par

There is another way to describe this topology.
In $\Cal C (G)$, a sequence $(C_n)_{n \ge 1}$ converges to $C$
if and only if:
\roster
\item"(3.4a)"
for any strictly increasing map $\varphi : \bold N \to \bold N$,
and for any sequence $(g_{\varphi (n)})_{n \ge 1}$
converging to some $g \in G$,
with $g_{\varphi(n)} \in C_{\varphi(n)}$ for all $n \ge 1$,
we have $g \in C$;
\item"(3.4b)"
for any $g \in C$, there exists a sequence $(g_n)_{n \ge 1}$
converging to $g$, with $g_n \in C_n$ for all $n \ge 1$.
\endroster
This has been observed by many authors; a proof 
can be found in  
\cite{BenPe--92, Section~E.1} or
in \cite{GuiR\'e--06, Section 2.1}.
\medskip

\bigskip
\head{\bf
4.~Closed subgroups of the group $\bold C$
}\endhead
\medskip

The present section is an exposition of \cite{HubPo--79}.
We describe the space $\Cal C (\bold C)$ in four steps:
first we describe the closed subspace 
$\Cal C_{nl}(\bold C)$
of closed subgroups that are {\it not lattices} and
the open subspace $\Cal L (\bold C)$ of lattices,
then we analyze the ways in which a sequence of lattices 
can converge to a non--lattice,
and finally we give a global description of 
$\Cal C(\bold C) = \Cal C_{nl}(\bold C) \cup \Cal L (\bold C)$.

\medskip
\head{\bf
4.I.~Closed subgroups in $\bold C$ which are not lattices
}\endhead
\medskip

There is an obvious identification of the
space $\Cal C_{\bold R}(\bold C)$ 
of closed subgroups of $\bold C$ isomorphic to $\bold R$ 
with the real projective line $\bold P^1$ and thus
 $\Cal C_{\bold R}(\bold C)$
is   homeomorphic to a circle.
Each closed subgroup $C \subset \bold C$ isomorphic to $\bold Z$
is contained in a unique group
$\overline{C} \in \Cal C_{\bold R}(\bold C)$
and is determined by $\overline C$ and the \lq\lq norm\rq\rq \ 
$|C| :=\min \{ \vert z \vert \in \bold R_+^* \, : \,
z \in  C, \, z \ne 0 \}$. Correspondingly,
the space $\Cal C_{\bold Z}(\bold C)$ 
is  homeomorphic to the direct product
of $\bold P ^1$ with an open interval. Moreover, $C\to \{0\}$
as $|C|\to\infty $, and for $\lambda\in[0,1]$
we have $\lambda C\to \overline C$ as $\lambda\to 0$.
Thus the space of closed subgroups of $\bold C$
isomorphic to one of $\{0\}$, $\bold Z$, $\bold R$
is homeomorphic to a cone
$$
\Cal C_{nl}^- (\bold C) \, = \,
\big( \bold P^1 \times [-1,0] \big)
\hskip.1cm / \hskip.1cm
\big( (x,-1) \sim (y,-1) \big) ,
$$
with the vertex $(*,-1)$ corresponding to $\{0\}$,
points $(x,t)$ with $-1 < t < 0$ corresponding to infinite cyclic groups,
and the base $\bold P^1 \times \{0\}$ 
corresponding to groups isomorphic to $\bold R$.
(The minus  sign indicates that $\Cal C_{nl}^- (\bold C)$
is the lower hemisphere of a $2$--sphere that plays an important
role in what follows.)
\par

A closed subgroup $C \subset \bold C$ isomorphic to $\bold R \oplus  \bold Z$
has an identity component $C^0 \in \Cal C_{\bold R}(\bold C)$
and is determined by $C^0$ and the \lq\lq norm\rq\rq \ 
$\min \{ \vert z \vert \in \bold R_+^* \, : \, z \in C, \, z \notin C^0 \}$;
it is convenient to parametrise $C$ by $C^0$ 
and {\it the inverse} of this norm. This parameterisation gives a homeomorphism
from $\Cal C_{\bold R \oplus \bold Z}(\bold C)$ 
to  $\bold P ^1 \times ]0,1[$. This extends to an identification
of the space of closed subgroups of $\bold C$
isomorphic to one of $\bold R$, $\bold R \oplus \bold Z$, $\bold C$
with the cone
$$
\Cal C_{nl}^+ (\bold C) \, = \,
\big( \bold P^1 \times [0,1] \big)
\hskip.1cm / \hskip.1cm
\big( (x,1) \sim (y,1) \big),
$$
with the vertex $(*,1)$ corresponding to $\bold C$,
points $(x,t)$ with $0 < t < 1$ 
to groups isomorphic to $\bold R \oplus \bold Z$, 
and the base of the cone $\bold P^1 \times \{0\}$ 
corresponding to groups isomorphic to $\bold R$. 

By combining these observations we obtain the following proposition.
This proposition is in the paper of John Hubbard and Ibrahim Pourezza
\cite{HubPo--79} but Hubbard informs us that they learned it from
 Adrien Douady. Moreover, it was Douady who suggested that they determine
the homeomorphism type of $\Cal C (\bold R^2)$. 
Apparently, this problem was of interest to Bourbaki at the time.

\bigskip

\proclaim{4.1.~Proposition}
The space of closed subgroups of $\bold C$ which are not lattices
is homeomorphic to a $2$--sphere
$$
\Cal C (\bold C) \smallsetminus \Cal L (\bold C) \, = \,
\Cal C_{\{0\}, \bold Z, \bold R, \bold R \oplus \bold Z, \bold C}
(\bold C) \, = \,
\Cal C_{nl}^- (\bold C) \bigcup_{ \Cal C_{\bold R} (\bold C) }
\Cal C_{nl}^+ (\bold C) \, \approx \,
\bold S^2.
$$
\endproclaim

\medskip
\head{\bf
4.II.~Lattices in $\bold C$
}\endhead
\medskip

   Classically, one defines for any lattice $L \subset \bold C$
two complex numbers
$$
g_2(L) \, = \, 60 \sum_{z \in L, z \ne 0} z^{-4}
\hskip.3cm \text{and} \hskip.3cm
g_3(L) \, = \, 140 \sum_{z \in L, z \ne 0} z^{-6} 
\tag{4.1}
$$
and one denotes by
$\Delta(L) \, = \, g_2(L)^3 - 27 g_3(L)^2$
the {\it discriminant} of $L$. 
The surface
$$
\Sigma \, = \, \{ (a,b) \in \bold C^2 \mid
a^3 -27b^2 = 0 \} 
$$
has an isolated singularity (cusp) at the origin
and  is smooth elsewhere.
Set
$$
T \, = \, \Sigma \cap \bold S^3 ,
$$
where
$\bold S^3 = \{(a,b) \in \bold C^2
\hskip.1cm \vert \hskip.1cm
\vert a \vert^2 + \vert b \vert^2 = 1 \}$
is the unit $3$--sphere. The smooth curve
$T$  is a trefoil knot.
It is a classical and basic result that $\Delta(L) \ne 0$
and that, moreover, the mapping
$$
\underline{g} \hskip.1cm : \hskip.1cm
\left\{
\aligned
\Cal L (\bold C) \, &\longrightarrow \,
\hskip.5cm \bold C^2 \smallsetminus \Sigma
\\
L \hskip.3cm &\longmapsto \, \left(g_2(L),g_3(L)\right)
\endaligned
\right.
\tag{4.2}
$$
is a homeomorphism.
There are two classical methods of proving this.
One method
uses the modular function, often denoted by $J$;
see for example \S~VIII.13 in \cite{SacZy--65} 
or \S~4 of Chapter II.4 in \cite{HurCo--64}. The 
other
method 
provides an explicit inverse to $\underline{g}$;
this associates to $(a,b) \in \bold C^2 \smallsetminus \Sigma$ the period lattice of the holomorphic
1-form $\frac{dX}{Y}$ on the genus-1
 plane  projective curve with equation
$Y^2 Z = 4X^3 - a X Z^2 -b Z^3$;
see for example \S~7.D in \cite{Mumfo--76}. 

\bigskip

Let $\Cal L^{\operatorname{umod}}(\bold C)$ be the subspace
of $\Cal L(\bold C)$ of unimodular lattices.
The natural action of the group $SL_2(\bold R)$
on $\Cal L^{\operatorname{umod}}(\bold C)$
is transitive
and the isotropy subgroup of the lattice $\bold Z [i]$
is $SL(2,\bold Z)$;
hence $\Cal L^{\operatorname{umod}}(\bold C)$
is homeomorphic to $SL_2(\bold R)/SL_2(\bold Z)$.
The universal covering of this space
is the universal covering 
$\widetilde{SL}_2(\bold R)$ of $SL_2(\bold R)$,
which is homeomorphic to an open $3$--disc.
It follows
\footnote{
This is in sharp contrast with the situation for $n \ge 3$.
In that case, the universal covering of the space
$SL_n(\bold R) / SL_n(\bold Z)$ of lattices in $\bold R^n$
is the two--sheeted covering $\widetilde{SL}_n(\bold R)$,
which is homotopic to its maximal compact subgroup 
$\operatorname{Spin}(n)$. 
In particular, $\widetilde{SL}_n(\bold R)$
is not contractible.
}
that the higher homotopy groups
$\pi_j(\Cal L (\bold C))$, $j \ge 2$, are trivial
and that $\pi_1(\Cal L (\bold C))$ is isomorphic to
the inverse image $\widetilde{SL}_2(\bold Z)$
of $SL_2(\bold Z)$ in the universal covering
of the group $SL_2(\bold R)$.

\bigskip

There are two continuous actions
of the multiplicative group $\bold C^*$
that are both natural and  important in the present context:
$$
\aligned
\bold C^* 
\hskip.2cm \text{acts on} \hskip.2cm
\Cal C(\bold C) 
\hskip.5cm \text{by} \hskip.5cm
&\hskip.5cm
(s,C) \hskip.1cm \longmapsto \hskip.1cm
\sqrt{s} \, C
\\
\bold C^* 
\hskip.2cm \text{acts on} \hskip.2cm
\bold C^2 
\hskip.5cm \ \text{by} \hskip.5cm
&\big( s,(a,b)\big) \hskip.1cm \longmapsto \hskip.1cm
\big( s^{-2}a, s^{-3}b \big) .
\endaligned
\tag{4.3}
$$
Several remarks are in order:
\roster
\item"---"
the subgroup $\sqrt{s} \, C$ is well defined
even though $\sqrt s$ is only defined up to a sign,
since $-C = C$;
\item"---"
the action of $\bold C^*$ on 
$\Cal C_{\bold Z^2, \bold R \oplus \bold Z, \bold Z} (\bold C)$
defined this way is faithful
(the action defined by $(s,C) \longmapsto sC$
would not be);
\item"---"
each of the subspaces
$\Cal L (\bold C)$, $\Cal C_{\bold R}(\bold C)$, 
$\Cal C_{\bold Z}(\bold C)$, $\{0\}$,
$\Cal C_{\bold R \oplus \bold Z}(\bold C)$, $\{\bold C\}$
of $\Cal C (\bold C)$ is invariant by $\bold C^*$;
\item"---"
the action of $\bold C^*$ on $\Cal C_{\bold Z}(\bold C)$
is free and transitive;
\item"---"
the hypersurface $\Sigma$ in $\bold C^2$ 
is $\bold C^*$--invariant;
\item"---"
the mapping $\underline{g}$ of (4.2)
is $\bold C^*$--equivariant
(this carries over to $\underline{g}'$, see (4.5) below).
\endroster
The actions of the subgroups $\bold R_+^*$ (positive reals) and 
$\bold S^1 = \left\{z \in \bold C^* \mid \vert z \vert = 1 \right\}$
obtained by restriction will also play a role below.  
\par

The action of $\bold R_+^*$ on $\Cal L (\bold C)$ is free,
and its orbits are transverse to 
$\Cal L^{\operatorname{umod}}(\bold C)$;
similarly, the action of $\bold R_+^*$ 
on $\bold C^2 \smallsetminus \Sigma$ is free,
and its orbits are transverse to $\bold S^3 \smallsetminus T$.
It follows that $\Cal L (\bold C)$
is homeomorphic to a direct product
\footnote{
There are at least three tempting choices
for a subspace $\Cal L^*(\bold C)$ which intersects
every $\bold R_+^*$--orbit exactly once,
so that $\Cal L (\bold C)$ is homeomorphic to the direct product
$\Cal L^*(\bold C) \times \bold R_+^*$;
and each $\Cal L^*(\bold C)$ is homeomorphic
to $\bold S^3 \smallsetminus  T$.
One choice is the space $\Cal L^{\operatorname{umod}}(\bold C)$
introduced here;
another choice is the space $\Cal L^{\operatorname{short}=1}(\bold C)$
of lattices $L$ whose shortest vector has norm $1$,
namely of lattices with
$\min \left\{ \vert z \vert \hskip.1cm : \hskip.1cm z \in L, \hskip.1cm z \ne 0\right\} = 1$;
and a third choice is the space of lattices $L$ with
$\vert g_2(L) \vert ^2 + \vert g_3(L) \vert ^2 = 1$.
Each of these three choices has its own virtues.
For example, is it apparent that there are sequences in
$\Cal L^{\operatorname{umod}}(\bold C)$
which converge to a subgroup isomorphic to $\bold R$,
and sequences in $\Cal L^{\operatorname{short}=1}(\bold C)$
which converge to subgroups isomorphic to $\bold Z$.
}
$\Cal L^{\operatorname{umod}}(\bold C) \times \bold R_+^*$,
and that we have a homeomorphism
$$
\Cal L^{\operatorname{umod}}(\bold C) \, \approx \,
\bold S^3 \smallsetminus  T .
\tag{4.4}
$$
In particular, the fundamental group $\widetilde{SL}_2(\bold Z)$
of the left--hand space is isomorphic to
the group of the trefoil knot,
also known as the Artin braid group on three strings.
(The idea of the  argument leading to (4.4) is due
to Daniel Quillen; see \cite{Milno--71, \S~10}.)
\par

The action of $\bold R_+^*$ on $\Cal C_{\bold Z}(\bold C)$
is also free (as already observed), 
and its orbits are  transverse to 
the subspace of subgroups of the form $\bold Z w$ 
with $w$ of modulus one in $\bold C$;
similarly, the action of $\bold R_+^*$ on $\Sigma \smallsetminus \{(0,0)\}$ is free,
and its orbits are transverse to $T$.
The three spaces $\Cal C_{\bold Z}(\bold C)$,
$\Sigma \smallsetminus \{(0,0)\}$, and  $T \times \bold R_+^*$
are naturally homeomorphic to each other.

\medskip
\noindent{\it
An alternative description of 
$SL_2(\bold R) / SL_2(\bold Z)$.
}

For the reader familiar with the geometry of $3$--manifolds, 
we describe an alternative manner of seeing 
that $SL_2(\bold R)/SL_2(\bold Z)$,
better viewed now as $PSL_2(\bold R)/PSL_2(\bold Z)$,
is homeomorphic to the complement of the trefoil knot. 
\par
 
The action of $PSL_2(\bold R)$ on the upper half--plane $\Cal H$ 
by fractional linear transformations induces 
a simply--transitive action on the unit tangent bundle
$U\Cal H$ (in the hyperbolic metric). 
Thus we may identify $PSL_2(\bold R)$ with $U \Cal H$ and 
$SL_2(\bold R) / SL_2(\bold Z)$ with $U\Cal H / PSL_2(\bold Z)$. 
The latter $3$--manifold is a Seifert fibre space 
whose base orbifold is the familiar $\Cal H / PSL_2(\bold Z)$, 
a once--punctured sphere with two cone points
where the local groups are cyclic of order $2$ and $3$, respectively. 
It is well--known that this Seifert fibred space 
is the complement of the trefoil knot.

Indeed the trefoil knot is the $(2,3)$--torus knot, 
and for any $p,q$ coprime,
the complement in $\bold S^3$ 
of the $(p,q)$--torus knot $T_{p,q}$
is a Seifert fibred space 
whose base orbifold is a disc with two cone points, 
one of order $p$ and one of order $q$. 
Perhaps the easiest way to see this is to observe that 
the non--singular fibres of
the map from $\bold S^3 \subset\bold C^2$
to $\bold C \cup \{\infty\} =\bold S^2$ (= the complex projective line) 
given by $(z_1,z_2) \longmapsto z_1^p/z_2^q$ 
are $(p,q)$--torus knots;
deleting one non-singular fibre gives the structure we seek.

\medskip

Returning now to the main argument, 
note that the mapping $\underline{g}$ of (4.2)
extends to a homeomorphism
$$
\underline{g'} \hskip.1cm : \hskip.1cm
\Cal L (\bold C) \cup \Cal C_{\bold Z}(\bold C)
\longrightarrow \bold C^2 \smallsetminus \{(0,0)\}
\tag{4.5}
$$
which is onto and $\bold R_+^*$--equivariant.
For any $w \in \bold C^*$, we have
$$
\aligned
g_2(\bold Z w ) \, &= \,
60 \sum_{n \in \bold Z, n \ne 0} (nw)^{-4} 
\phantom{1} \, = \,
\frac{4 \pi^4}{3 w^4}
\\
g_3(\bold Z w ) \, &= \,
140 \sum_{n \in \bold Z, n \ne 0} (nw)^{-6} 
\, = \,
\frac{2^3 \pi^6}{3^3 w^6}
\endaligned
\tag{4.6}
$$
(by an easy computation,
or by Formula (12.8) of Chapter VIII in \cite{SacZy--65}).
Summing up:

\bigskip

\proclaim{4.2.~Proposition}
With the notation above,
the space of lattices in $\bold C$ 
is homeomorphic to $\bold C^2 \smallsetminus \Sigma$,
equivalently
$\bold R^*_+ \times \left( \bold S^3 \smallsetminus T \right)$.
Its fundamental group is isomorphic to the trefoil knot group
and its universal covering is a $4$--disc.
\par
The homeomorphism 
$\Cal L (\bold C) \approx \bold C^2 \smallsetminus \Sigma$
extends to a homeomorphism 
$\Cal C_{\bold Z, \bold Z^2}(\bold C) 
\approx \bold C^2 \smallsetminus \{(0,0)\}$,
and further to a homeomorphism
$\Cal C_{\{0\}, \bold Z, \bold Z^2}(\bold C) \approx \bold C^2$.
\endproclaim

\medskip
\head{\bf
4.III.~Degeneration of lattices in $\bold C$
}\endhead
\medskip

Let $L$ be a lattice in $\bold C$.
Choose $w_1,w_2 \in L$ 
with $w_1$ of smallest absolute value 
among elements of $L \smallsetminus \{0\}$,
with $w_2$ of smallest absolute value
among elements of $L \smallsetminus \bold Z w_1$,
and with $\operatorname{Im}(w_2 / w_1) > 0$.
Then $(w_1,w_2)$ is a positively oriented
$\bold Z$--basis of $L$.
Define the {\it minimal length} $\ell_1(L) := \vert w_1 \vert$,
the {\it second minimal length} $\ell_2(L) := \vert w_2 \vert$,
and the {\it distortion coefficient} 
$\kappa(L) = \ell_2(L) / \ell_1(L)$ of $L$.
Clearly, these depend on $L$ only, 
and not on the choice of $w_1$ and $w_2$
(unless $\kappa(L) = 1$, the subgroup
$\bold Z w_1$ is also well--defined).
The definition of $\ell_1$ 
carries over to $\Cal C_{\bold Z}(\bold C)$,
and that of $\ell_2$ ({\it cum grano salis})
to $\Cal C_{\bold R \oplus \bold Z} (\bold C)$.
The functions
$$
\aligned
\ell_1 \, &: \,  \Cal L (\bold C) \cup 
            \Cal C_{\bold Z}(\bold C) 
\, \longrightarrow \, ]0,\infty[ ,
\\ 
\ell_2 \, &: \,  \Cal L (\bold C) \cup 
            \Cal C_{\bold R \oplus \bold Z} (\bold C)
\, \longrightarrow \, ]0,\infty[, 
\\
\kappa \, &: \,  \Cal L (\bold C) \longrightarrow [1,\infty[
\endaligned
$$
are continuous.
Observe that $w_2/w_1$ always belongs to the 
standard fundamental domain for the action of $PSL_2(\bold Z)$
on the Poincar\'e half--plane, hence
the angle spanned by $w_1$ and $w_2$
is always between $\pi/3$ and $2\pi/3$.

\bigskip

\proclaim{4.3.~Proposition}
Let $(L_n)_{n \ge 1}$ be a sequence in $\Cal L (\bold C)$.
\roster
\item"(i)"
$\lim_{n \to \infty}L_n = \{0\}$ 
if and only if $\ell_1(L_n) \to \infty$;
\item"(ii)"
$\lim_{n \to \infty}L_n = \bold C$ 
if and only if $\ell_2(L_n) \to 0$;
\item"(iii)"
if $\lim_{n \to \infty}L_n = C$ for some $C \in \Cal C(\bold C)$
with $C \notin \Cal L (\bold C)$, $C \ne \{0\}$,  $C \ne \bold C$,
then $\kappa(L_n) \to \infty$.
\endroster
\endproclaim

\demo{Proof}
We leave (i) and (ii) to the reader and prove (iii).
We shall assume that $\kappa (L_n)$ does not tend to $\infty$
and reach a contradiction.
\par
By assumption, there exists a subsequence $(L_{n_k})_{k \ge 1}$
such that $\kappa(L_{n_k})$ converges to some $\kappa_0 < \infty$.
Also, $\ell_1(L_{n_k})$ is bounded above because $C \ne \{0\}$,
and $\ell_2(L_{n_k})$ is bounded below because $C \ne \bold C$;
since their quotients $\kappa(L_{n_k})$ are bounded,
both $\ell_1(L_{n_k})$ and $\ell_2(L_{n_k})$ 
are bounded simultaneously above and below.
Upon extracting a further appropriate subsequence,
we can assume that $L_{n_k}$ contains vectors
$w_1(n_k), w_2(n_k)$ defined as above
and that we have limits in $\bold C$, say
$w_1 = \lim_{k \to \infty} w_1(n_k)$, 
$w_2 = \lim_{k \to \infty} w_2(n_k)$.
But then the limit group $C$ is 
the lattice spanned by $w_1$ and $w_2$,
and this is the desired contradiction.
\hfill $\square$
\enddemo
\bigskip

Our next proposition is a refinement of (iii) above;
the proof is left to the reader.
We denote by $[w]$ the class in $\bold C^* / \{\pm \operatorname{id} \}$
of a vector $w \in \bold C^*$,
and by $\bold R w$ its class 
in $\bold P^1$.

\bigskip

\proclaim{4.4.~Proposition}
Let $L_n \to C$ be as in the previous proposition,
and let $v \in \bold C^*$.
\roster
\item"(i)" 
$C = \bold Z v$ if and only if 
$\kappa(L_n) \to \infty$ and 
$[w_1(L_n)] \to [v] \in \bold C / \{\pm \operatorname{id} \}$.
\item"(ii)"
$C = \bold R v \oplus \bold Z (iv)$ if and only if
$\kappa(L_n) \to \infty$ and 
$[w_2(L_n)] \to [iv] \in \bold C / \{\pm \operatorname{id} \}$.
\item"(iii)"
$C = \bold R v$ if and only if
$\ell_1(L_n) \to 0$, $\ell_2(L_n) \to \infty$,
and $\bold R w_1(L_n) \to \bold R v \in \bold P^1$.
\endroster
\endproclaim

\demo{Comment on (iii)}
If $\ell_1(L_n) < \ell_2(L_n)$,
recall that $w_1(L_n)$ is well defined up to sign,
so that $\bold R w_1(L_n)$ is well defined in $\bold P^1$.
Since $L_n \to C$ by hypothesis, 
$\bold R w_1(L_n) \to \bold R v$ in $\bold P^1$.
\enddemo

\bigskip

\head{\bf
4.IV.~The homeomorphism $\Cal C (\bold C) \approx \bold S^4$
of \cite{HubPo--79}
}\endhead
\medskip

In preparation for Theorem~4.6, 
which is a detailed version of Theorem~1.2,
it is convenient to record the following identifications.
The $4$--sphere of Theorem 1.2 is identifed with
the one--point compactification $\bold C^2 \cup \{\infty\}$
of $\bold C^2$, 
its north and south poles to $\infty$ and $(0,0)$
respectively,
and its equator to the unit sphere
$\bold S^3 = \{(a,b) \in \bold C^2 \mid 
\Vert (a,b) \Vert = 1 \}$,
where $\Vert (a,b) \Vert = \sqrt{ \vert a \vert ^2 + \vert b \vert ^2}$.
We denote by $B$ the open unit ball
$\{(a,b) \in \bold C^2 \mid 
\Vert (a,b) \Vert < 1 \}$,
and we denote by 
$$
\gamma \, : \,  \bold C^2 \longrightarrow 
\Cal C_{\{0\},\bold Z,\bold Z^2}(\bold C)
$$ 
the inverse of the homeomorphism $\underline{g'}$
which appears in~(4.5) and in Proposition 4.2.
\par

As our aim is to define a homeomorphism $f$
from $\bold C^2 \cup \{\infty\}$ to $\Cal C (\bold C)$,
we will first define its restriction to $\overline{B}$.
We want the image  $f(\bold S^3 \smallsetminus T)$
to be the set of unimodular lattices; in contrast,
$\gamma(\bold S^3 \smallsetminus T)$
contains lattices with arbitrarily large coareas
(compare with Footnote~3).
The construction of $f$ involves several auxiliary mappings.
\medskip

The coarea is traditionally defined as a function
$\Cal L (\bold C) \longrightarrow ]0,\infty[$.
We extend it by defining
$\operatorname{coarea}(C) = \infty$
if $C = \{0\}$ or $C \cong \bold Z$
and $\operatorname{coarea}(C) = 0$
if $C \cong \bold R \oplus \bold Z$ or $C = \bold C$;
the resulting extension
$\Cal C (\bold C) \smallsetminus \Cal C_{\bold R}(\bold C)
\longrightarrow [0,\infty]$
is continuous.
Therefore, it makes sense to define
$$
\Cal C^{\operatorname{coarea} \ge 1}(\bold C)
\, = \,
\{ C \in \Cal C (\bold C) \smallsetminus \Cal C_{\bold R}(\bold C)
\mid
\operatorname{coaera}(C) \ge 1 \} ,
$$
which is a subset of $\Cal C_{\bold Z^2, \bold Z, \{0\}}(\bold C)$
containing $\Cal C_{\bold Z, \{0\}}(\bold C)$.
\par

The second auxiliary mapping is the retraction
$$
\pi \, : \, \bold C^2 \smallsetminus \{(0,0)\}
\hskip.2cm \longrightarrow \hskip.2cm \bold S^3 ,
\hskip.5cm
(a,b) \, \longmapsto \, \bold R_+^*(a,b) \cap \bold S^3
$$
which assigns to $(a,b)$ the intersection with $\bold S^3$
of its $\bold R_+^*$--orbit, see (4.3).
Observe that $\pi$ is $\bold S^1$--equivariant
(for the actions of $\bold S^1$ viewed as a subgroup of $\bold C^*$).
\par

We define a third continuous mapping
$$
\varphi \, : \, \bold C^2 \smallsetminus \{(0,0)\} 
\hskip.2cm \longrightarrow \hskip.2cm \ ]0,\infty] ,
\hskip.5cm
(a,b) \, \longmapsto \,  
\sqrt{ \operatorname{coarea}\big( \gamma(\pi(a,b)) \big) } .
$$
Observe that $\varphi$ is $\bold C^*$--invariant;
moreover, for $(a,b) \in \bold C^2 \smallsetminus \{0,0\}$, 
we have $(a,b) \in \Sigma$ 
if and only if $\pi(a,b) \in T$,
if and only if $\gamma(\pi(a,b)) \cong \bold Z$,
if and only if
$\operatorname{coarea}\big( \gamma(\pi(a,b)) \big) = \infty$.
\par

The fourth and last auxiliary mapping is a continuous function
$$
h \, : \, \overline{B} \longrightarrow \bold R_+ \cup \{\infty\}
$$
which the following properties:
\roster
\item"(i)"
$h(0,0) = 0$;
\item"(ii)"
for $(a,b) \in \overline{B} \smallsetminus \{(0,0)\}$,
there exists an order--preserving homeomorphism
$\eta_{\varphi (a,b)}$ from $[0,1]$ onto $[0,\varphi(a,b)]$
such that
$h(a,b) = \eta_{\varphi(a,b)} (\Vert (a,b) \Vert)$.
\endroster
The notation indicates that the homeomorphism
$\eta_{\varphi(a,b)}$ depends on $(a,b)$
through its image by $\varphi$ only;
in particular, the homeomorphisms associated to 
$(a,b)$, $(t^{-2}a,t^{-3}b)$, and 
$(e^{-2i\theta}a, e^{-3i\theta}b)$
are identical, for any $t$ and $\theta$ in $\bold R$.
It follows from (i) and (ii) that, for $(a,b) \in \overline{B}$:
\roster
\item"(iii)"
$h(a,b) = \infty$ if and only if $(a,b) \in T$;
\item"(iv)"
$h(a,b) = 0$ if and only if $(a,b) = (0,0)$;
\item"(v)"
$\frac{1}{h(a,b)} \gamma(\pi(a,b))
\in \Cal L^{\operatorname{umod}}(\bold C)$
if and only if 
$(a,b) \in \bold S^3 \smallsetminus T$
\endroster
The other properties of the mapping $h$ 
do not play any important role below.
One possibility (out of many others)
would be to require  $\eta_{\varphi(a,b)}$
to be a fractional linear homeomorphism 
tangent to the identity at the origin,
and therefore to set
$$
h(a,b) \, = \, \frac{ 
\Vert (a,b) \Vert 
}{
1 + \left( \frac{1}{\varphi(a,b)} - 1 \right) \Vert (a,b) \Vert
}
\tag{4.7}
$$
if $(a,b) \ne (0,0)$, and $h(0,0) = 0$.

Checking that the map defined by (4.7) is continuous
is equivalent to checking that the homeomorphisms 
$\eta_c : [0,1] \longrightarrow [0,c]$
depend continuously on $c \in ]0,\infty]$,
with respect to the uniform topology
on the space of mappings
from $[0,1]$ to the compact interval $[0,\infty]$.
\par

For a closed subgroup $C$ of $\bold C$, we denote by 
$$
\sigma(C) \, = \, \{ z \in \bold C \mid 
\overline{z} \in C \}
$$
the complex conjugate  of $C$.

\bigskip

\proclaim{4.5.~Proposition}
Let $f : \overline{B} \longrightarrow
\Cal C_{\bold Z^2, \bold R, \bold Z, \{0\}}(\bold C)$
be the mapping defined by
$$
\matrix
&f(a,b) \, &= \, &\frac{1}{h(a,b)} \gamma(\pi(a,b) )
\, &= \, &\frac{1}{\varphi(a,b)} \gamma(a,b)
&\hskip.5cm &\text{for} \hskip.2cm
&(a,b) &\in &\bold S^3 \smallsetminus T,
\hskip.6cm
\\
&f(a,b) \, &= \, &\bold R w 
&&&\hskip.3cm &\text{for} \hskip.2cm
&(a,b) &\in &T,
\hskip1.5cm
\\
&f(a,b) \, &= \, &\frac{1}{h(a,b)} \gamma(\pi(a,b))
&&&\hskip.3cm &\text{for} \hskip.2cm
&(a,b) &\in &B \smallsetminus \{(0,0)\},
\\
&f(0,0) &= &\{0\},
&&& &&&
\endmatrix
$$
where $w$ is defined by $\gamma(a,b) = \bold Z w$
in the second case.
\par

Then $f$ is a homeomorphism of $\overline{B}$ onto
$\Cal C^{\operatorname{coarea} \ge 1}(\bold C) 
\cup \Cal C_{\bold R}(\bold C)$.
Moreover, $f$ is compatible with complex conjugation:
$f(\overline{\vphantom{b}a},\overline{b}) = \sigma(f(a,b))$
for all $(a,b) \in \overline{B}$.
\endproclaim

\demo{Proof} To show that $f$ is a homeomorphism as claimed,
since $\overline{B}$ is compact,
it is sufficient to show that $f$ is continuous, injective,
and that its image is 
$\Cal C^{\operatorname{coarea} \ge 1}(\bold C) 
\cup \Cal C_{\bold R}(\bold C)$.
\medskip

{\it First step: $f$ is injective with image as indicated.}
The domain of $f$ is the disjoint union of four subsets
$$
\{(0,0)\}, \hskip.5cm
T, \hskip.5cm
\Sigma \cap (B  \smallsetminus \{(0,0)\}), \hskip.5cm
\overline{B} \smallsetminus (\Sigma \cap \overline{B}) .
$$
From the definitions of $\gamma$, $\varphi$, $h$, and $f$,
it is clear that these subsets are mapped by $f$ to
$$
\{0\} , \hskip.5cm
\Cal C_{\bold R}(\bold C), \hskip.5cm
\Cal C_{\bold Z}(\bold C), \hskip.5cm
\Cal L^{\operatorname{coarea} \ge 1}(\bold C)
$$
respectively.
As it is obvious that $f$ induces a bijection
from the one--point subset $\{(0,0)\}$ of $\overline{B}$
onto the one--point subset $\{0\}$ of $\Cal C(\bold C)$,
there are three more cases to check.
\medskip

(i) The group $\bold S^1$ acts freely transitively
on both $T$ and $\Cal C_{\bold R}(\bold C)$.
Since $f$ is $\bold S^1$--equivariant,
it follows that $f$ induces a bijection from $T$
onto $\Cal C_{\bold R}(\bold C)$.
\medskip

(ii) Let $(a,b) \in \Sigma \cap (B  \smallsetminus \{(0,0)\})$.
Then $\pi(a,b) \in T$,
there exists $w \in \bold C$, $w \ne 0$,
such that $f(a,b) = \bold Z w$,
and $w$  is well--defined up to a change of sign.
Moreover, $\varphi(a,b) = \infty$
and $h(a,b) = \frac{ \Vert (a,b) \Vert }{1 - \Vert (a,b) \Vert}$.
Thus $f$ induces a bijection 
from the part $\bold R_+^* (a,b) \cap B$ in $B$
of the $\bold R_+^*$--orbit of $(a,b)$
onto the set of infinite cyclic subgroups of $\bold C$
contained in $\bold R w$.
It follows that $f$ induces a bijection
from $ \Sigma \cap (B  \smallsetminus \{(0,0)\})$
onto $\Cal C_{\bold Z}(\bold C)$.
\medskip
 
(iii) Let $(a,b) \in 
\overline{B} \smallsetminus (\Sigma \cap \overline{B})$.
Then $\pi(a,b) \in \bold S^3 \smallsetminus T$
and $L_0 := \gamma(\pi(a,b)) \in \Cal L (\bold C)$.
Thus $f$ induces a bijection from $\bold R_+^*(a,b) \cap \overline{B}$
onto the set of lattices of the form $\frac{1}{t} L_0$,
with $0 < t \le \sqrt{\operatorname{coarea}(L_0)}$,
namely with 
$\operatorname{coarea}  \left( \frac{1}{t} L_0 \right) \ge 1$. 
It follows that $f$ induces a bijection from  
$\overline{B} \smallsetminus (\Sigma \cap \overline{B})$
onto $\Cal L^{\operatorname{coarea} \ge 1}(\bold C)$.
 
 \medskip
 
{\it Second step: $f$ is continuous.}
It is clear that $f$ is continuous at
any point of $\overline{B} \smallsetminus (T \cup \{(0,0)\})$, 
because the mappings $\gamma$, $\pi$, $\varphi$, 
and $h$ are continuous (and, moreover, $h(a,b) < \infty$). 
It remains to check the continuity of $f$
first at the origin and  at the points of $T$.
\medskip
 
(i) To show that $f$ is continuous at $(0,0)$
it is enough by Proposition 4.3.i to show 
that $\ell_1(f(a,b)) \to \infty$ when $(a,b) \to (0,0)$.
Since the function
$$
\bold S^3 \, \longrightarrow \bold R_+^*,
\hskip.5cm (a,b) \, \longmapsto \ell_1(\gamma(a,b))
$$
is continuous and positive on a compact domain,
its minimum is positive. Since 
\newline
$\lim_{(a,b) \to (0,0)} h(a,b) = 0$,
we have
$$
\ell_1(f(a,b)) \, = \,
\frac{1}{h(a,b)} \ell_1 ( \gamma(\pi(a,b))) \, \to \, \infty
\hskip.5cm \text{if} \hskip.5cm (a,b) \, \to \, (0,0).
 $$
Thus $f$ is indeed continuous at the origin.
\medskip
 
(ii) Let $(a_0,b_0) \in T$.
If $(a,b)$ tends to $(a_0,b_0)$ inside $T$,
then $f(a,b) \cong \bold R$ tends to $ f(a_0,b_0) \cong \bold R$
by the continuity of the $\bold S^1$--action
(which is transitive on both $T$ and $\Cal C_{\bold R}(\bold C)$)
and the $\bold S^1$--equivariance of $f$.
If $(a,b)$ tends to $(a_0,b_0)$ inside $\Sigma \cap B$,
it is a straightforward consequence of the definition of $f$
and of the continuity of $h$
that $f(a,b) \cong \bold Z$ tends to $f(a_0,b_0) \cong \bold R$.
Therefore, we can assume from now on 
that $(a,b)$ tends to $(a_0,b_0)$ inside 
$\overline{B} \smallsetminus (\Sigma \cap \overline{B})$, 
and we have to show that
$f(a,b) \cong \bold Z^2$ tends to $f(a_0,b_0) \cong \bold R$;
for this, we are going to check that the three conditions
of Proposition 4.4.iii are satisfied.
\par

We have
$\ell_1 (f(a,b)) = \frac{1}{h(a,b)} \ell_1(\gamma(\pi(a,b)))$.
Now the function $\ell_1 \circ \gamma$ is continuous on $\bold S^3$,
so  $\ell_1(\gamma(\pi(a,b)))$ is bounded,
and $h(a,b) \to \infty$.
Hence $\ell_1(f(a,b)) \to 0$ when $(a,b) \to (a_0,b_0)$.
\par

We have
$1 \le \operatorname{coarea}(f(a,b)) \le
\ell_1(f(a,b)) \ell_2(f(a,b))$.
The previous point implies that
$\ell_2(f(a,b)) \to \infty$ when $(a,b) \to (a_0,b_0)$.
\par

For $(a,b)$ inside
$\overline{B} \smallsetminus (\Sigma \cap \overline{B})$ 
and near enough to $(a_0,b_0)$,
we have $\ell_1(\gamma(a,b)) < \ell_2(\gamma(a,b)),$
so  $\bold R w_1(\gamma(a,b))$ is well-defined in $\bold P^1$;
moreover
$$
\bold R w_1 (\gamma(a,b)) \, = \,
\bold R w_1 (\gamma(\pi(a,b))) \, = \,
\bold R w_1 (f(a,b)) .
$$
As the mappings $\gamma$ and 
$(a,b) \longmapsto \bold R w_1 (\gamma(a,b))$
are continuous in appropriate domains,
we have $\bold R w_1(f(a,b)) \to \bold R v$.

\medskip

{\it Third step.} For each $C \in \Cal L (\bold C) \cup \Cal C_{\bold Z}(\bold C)$
we have $g_2(\sigma(C)) = \overline{g_2(C)}$, 
and similarly for $g_3$.
It follows that
$f(\overline{\vphantom{b}a},\overline{b}) = \sigma(f(a,b))$
for all $(a,b) \in \overline{B}$.
We leave the details to the reader. 
\hfill $\square$
\enddemo
\bigskip

We can now state the final result of this expository section.
We denote by $\delta$ the inversion of $\bold C^2 \cup \{\infty\}$
which exchanges $(0,0)$ and $\infty$,
and which is defined on other points by
$\delta(a,b) = \frac{(a,b)}{\Vert (a,b) \Vert^2}$.
For a closed subgroup $C$ of $\bold C$, we denote by 
$$
C^*  \, = \, \{ z \in \bold C \mid 
\operatorname{Im}(\overline{z}w) \in \bold Z
\hskip.2cm \text{for all} \hskip.2cm
w \in C \} 
$$
the  dual of $C$.
It is straightforward to check that this duality 
is a homeomorphism which
exchanges the following pairs of groups:
$$
\matrix
\{0\} && \leftrightsquigarrow && \bold C
\\
\bold Z w && \leftrightsquigarrow && \bold R w \oplus \bold Z(i/\overline{w})
\\
\bold R w && \leftrightsquigarrow && \bold R w
\\
\bold Z w_1 \oplus \bold Z w_2 && \leftrightsquigarrow && \
\frac{1}{\operatorname{Im}(\overline{w_1}w_2)}
\Big( \bold Z w_1 \oplus \bold Z w_2 \Big).
\endmatrix
$$
In particular, the fixed point set in $\Cal C (\bold C)$ 
of the duality is exactly the image 
$\Cal L^{\operatorname{umod}}(\bold C) \cup \Cal C_{\bold R}(\bold C)$
by $f$ 
of the fixed point set $\bold S^3$ of the inversion $\delta$.
The main result of this section (and of \cite{HubPo--79}) follows:

\bigskip

\proclaim{4.6.~Theorem}
Let $f : \bold C^2 \cup \{\infty\} \longrightarrow \Cal C (\bold C)$
be the extension of the mapping $f$ of Proposition 4.5
defined by $f(a,b) = f(\delta(a,b))^*$ when $(a,b) \notin \overline{B}$.
Then $f$ is a homeomorphism  
from $\bold C^2 \cup \{\infty\}$, which is homeomorphic to $\bold S^4$,
to the space $\Cal C (\bold C)$ of closed subgroups of $\bold C$, 
equipped with the Chabauty topology.
Moreover, $f$  has the following properties:
\roster
\item"(i)"
The images of the subspaces
$$
\bold C^2 \smallsetminus \Sigma, \hskip.2cm
\bold S^3 \smallsetminus T, \hskip.2cm
T, \hskip.2cm
(\Sigma \cap B) \smallsetminus \{(0,0)\}, \hskip.2cm
\{(0,0)\}, \hskip.2cm
\Sigma \cap (\bold C^2 \smallsetminus \overline{B}), \hskip.2cm
\{\infty\}
$$
of $\bold C^2 \cup \{\infty\}$ 
are respectively the subspaces
$$
\Cal L (\bold C), \hskip.2cm
\Cal L^{\operatorname{umod}}(\bold C), \hskip.2cm
\Cal C_{\bold R} (\bold C), \hskip.2cm
\Cal C_{\bold Z} (\bold C), \hskip.2cm
\{0\}, \hskip.2cm
\Cal C_{\bold R \oplus \bold Z}(\bold C), \hskip.2cm
\{ \bold C \}
$$
of $\Cal C (\bold C)$.
\item"(ii)"
Inversion in $\bold C^2$ corresponds to duality in $\Cal C (\bold C)$:
$$
f(\delta(a,b)) \, = \,  f(a,b)^*
\hskip.5cm \text{for all} \hskip.2cm (a,b) \in \bold C^2 \cup \{\infty\} .
$$
\item"(iii)"
For the action of  $\bold S^1$ defined by (4.3),
the mapping $f$ is equivariant:
$$
f(e^{-i2\theta}a,e^{-i3\theta}b) \, = \,  e^{i\theta/2}f(a,b)
\hskip.5cm \text{for all} \hskip.2cm 
(a,b) \in \bold C^2
\hskip.2cm \text{and} \hskip.2cm 
\theta \in \bold R .
$$
\item"(iv)"
The homeomorphism $f$ is compatible with complex conjugation:
$$
f(\overline{\vphantom{b}a},\overline{b}) \, = \,  \sigma(f(a,b))
\hskip.5cm \text{for all} \hskip.2cm (a,b) \in \bold C^2 \cup \{\infty\} .
$$
\item"(v)"
For the action of $\bold R_+^*$ defined by (4.3)
and for $(a,b) \in \bold C^2 \smallsetminus \Sigma$,
the lattices $f(a,b)$ and $\underline{g}^{-1}(a,b)$,
see (4.2),
are in the same $\bold R_+^*$--orbit.
\endroster
\endproclaim

\bigskip
\head{\bf
5.~Generalities on the real Heisenberg group
}\endhead
\medskip

The remainder of this article is dedicated to the study of the
3-dimensional Heisenberg group $H$ and its subgroups. 
\par
We begin by recalling a few basic properties of  $H$.
Its centre $Z(H)$ can be identified with the second factor $\bold R$
in the description $H = \bold C \times \bold R$ of (1.3);
moreover,  $Z(H)$ coincides with the commutator subgroup $[H,H]$.
There is a  canonical projection
$$
p \, : \, H \longrightarrow H/Z(H) = \bold C ,
\tag{5.1}
$$
as in (1.3), 
and $H$ is an extension of $\bold C$ by~$\bold R$
which does not split
(unlike the extension corresponding to Formula (1.4)).
It is useful to have the formula for commutators
$$
(z,t)(z',t')(z,t)^{-1}(z',t')^{-1}
\, = \,
\big( 0, \operatorname{Im}(\overline{z}z') \big) .
\tag{5.2}
$$
\par

We denote by $\operatorname{Aut}(H)$
the group of continuous 
\footnote{
Any automorphism of $H$ happens to be continuous,
but this fact does not play any role in our presentation.
}
automorphisms of $H$.
Any $\Phi \in \operatorname{Aut}(H)$ preserves the center 
and thus descends to a continuous automorphism of
$H / Z(H) \cong \bold R^2$. Hence we have a homomorphism
$$
\Pi \, : \, \operatorname{Aut}(H) \longrightarrow GL_2(\bold R) .
\tag{5.3}
$$
This is onto. 
The inner automorphisms act on $H$  by
$$
\operatorname{Int}_{(w,*)} (z , t)
\, = \,
(w,*) (z,t) (w,*)^{-1}
\, = \,
(z \hskip.1cm , \hskip.1cm t + \operatorname{Im}(\overline{w}z)) 
\tag{5.4}
$$
where $\operatorname{Int}_{(w,*)}$ is written for
$\operatorname{Int}_{(w,s)}$,
with $s$ arbitrary;
these automorphisms form a normal subgroup $\operatorname{Int}(H)$
of $\operatorname{Aut}(H)$ isomorphic to $H/Z(H) \cong \bold R^2$.
It follows from (5.4) that $\operatorname{Int}(H) \subset \ker (\Pi)$,
so that we have a sequence 
$$
\{1\} 
\hskip.1cm \longrightarrow \hskip.1cm
\operatorname{Int}(H)
\hskip.1cm \longrightarrow \hskip.1cm
\operatorname{Aut}(H)
\overset{\Pi}\to{\longrightarrow} \hskip.1cm
GL_2(\bold R)
\hskip.1cm \longrightarrow \hskip.1cm
\{1\} .
\tag{5.5}
$$

\bigskip

\proclaim{5.1.~Proposition}
The sequence (5.5) is exact and split,
so that the group of continuous automorphisms of $H$ 
is a semi--direct product
$$
\operatorname{Aut}(H)
\, = \,
\operatorname{Int}(H) \rtimes GL_2(\bold R)
\, \cong \,
\bold R^2   \rtimes GL_2(\bold R) .
\tag{5.6}
$$
In particular, in its natural topology, $\operatorname{Aut}(H)$
has two connected components.
\endproclaim

\demo{Note} We write $GL_2(\bold R)$ on the right--hand side of 
$\operatorname{Int}(H) \rtimes GL_2(\bold R)$,
but nevertheless we see $GL_2(\bold R)$ as acting from the left
on $\bold R^2 \cong \operatorname{Int}(H)$.
\enddemo

\demo{Proof}
The group $\operatorname{Aut}(H)$ contains a copy of $GL_2(\bold R)$,
acting by
$$
\left( \matrix a & b \\ c & d \endmatrix \right)
(x+iy \hskip.1cm , \hskip.1cm t)
\, = \,
\big((ax+by) + i(cx+dy) \hskip.1cm , \hskip.1cm
(ad - bc)t \big) .
\tag{5.7}
$$
Consequently, the homomorphism $\Pi$ in (5.5) has a section.
As we have already observed that
$\operatorname{Int}(H) \subset \ker (\Pi)$,
it remains to justify the opposite inclusion.
\par

Let $\Phi \in \ker (\Pi)$.
Since $\Phi$ acts as the identity on $H/Z(H)$, it is of the form
$$
(x+iy,t) \hskip.2cm \longmapsto \hskip.2cm (x+iy,t+\varphi(x+iy,t))
$$
for some mapping $\varphi : H \longrightarrow \bold R$.
The multiplication identity
$$
\Phi  (x+iy,t) \Phi (x'+iy',t') 
\hskip.2cm = \hskip.2cm 
\Phi  \big((x+iy,t)(x'+iy',t')\big)
$$
reduces to
$$
\varphi(x+iy,t) + \varphi(x'+iy',t') \hskip.2cm = \hskip.2cm
\varphi\big((x+iy,t)(x'+iy',t')\big) .
$$
Thus, $\varphi$ is a continuous group homomorphism with values in 
the abelian group $\bold R$.
It follows that $\varphi$ factors through
a homomorphism
$H/[H,H] \cong \bold R^2 \longrightarrow \bold R$,
namely a linear form on $\bold R^2$,
say $(x,y) \longmapsto -vx + uy$ for some $u,v \in \bold R$.
If $w = u+iv$, it follows from (5.4) that
$\Phi  = \operatorname{Int}_{w,*} \in \operatorname{Int}(H)$.
\hfill $\square$
\enddemo
\bigskip

We record here a few more observations concerning $\operatorname{Aut}(H)$.
\roster
\item"---"
The centre $Z(H)$ is fixed elementwise
by the derived group of $\operatorname{Aut}(H)$,
namely by $\bold R^2 \rtimes SL_2(\bold R)$.
\item"---"
The projection $p$ of (5.1) is $GL_2(\bold R)$--equivariant,
where the action on $\bold C$ is the standard action.
\item"---" The simple group $SL_2(\bold R)$ is a Levi factor
in the identity component of the group $\operatorname{Aut}(H)$.
There are many choices for such a factor,
and each one corresponds to a choice of a \lq\lq supplement\rq\rq \
$\bold C$ to the centre as one writes $H = \bold C \times \bold R$;
indeed, once the Levi factor $SL_2(\bold R)$ is chosen,
the slices $\bold C^* \times \{t_0\}$ are the $SL_2(\bold R)$--orbits
in $H$.
\endroster

\bigskip

\proclaim{5.2.~Proposition}
(i) Let $C$ be a closed subgroup of $H$ such that $C \cap Z(H) \ne \{e\}$.
Then the subgroup $p(C)$ of $\bold C$ is closed.
\par

In particular, if $C$ is any non--abelian closed subgroup of $H$,
then $p(C)$ is closed.
\smallskip

(ii) Let $C$ be a non--abelian closed subgroup of $H$.
Then, either $C$ contains $Z(H)$, or $C$ is a lattice in $H$.
In other words:
$$
\Cal C (H) \, = \,
\Cal L (H) \cup \Cal A (H) \cup \cc (H) .
$$
\smallskip

(iii) The assignment $C \longmapsto  p(C)$ defines a continuous
map 
$$
p_* \, : \, 
\Cal C (H) \smallsetminus \Cal A(H) \longrightarrow \Cal C(\bold C) .
$$

\smallskip

(iv) The spaces $\Cal L (H)$ and  $\Cal L_{!!} (H)$
are open in $\Cal C (H)$.

\smallskip

(v) The assignment $C \longmapsto C \cap Z(H)$ defines a continuous map 
$$
\Cal C (H) \smallsetminus \Cal A(H) \, \longrightarrow \, \Cal C(Z(H)) .
$$
\endproclaim

\demo{Remark} In (i), the condition that $C$ is non--abelian
cannot be removed. 
Consider for example the subgroup $A_{\theta}$ in $H = \bold C \times \bold R$   
generated by $(\theta,1)$ and $(1,1)$, 
where $\theta$ is an irrational real number.
This subgroup is closed, isomorphic to $\bold Z^2$,
and its projection in $\bold C$ is dense in a real line.
\par

Similarly, the map of (v) is not continuous on the whole of 
$\Cal C (H)$.

The mapping $p_*$ of (iii) does not extend continuously on $\Cal C (H)$;
see Proposition~8.8.
\enddemo

\demo{Proof} (i) Let $(z_n)_{n \ge 1}$ be a sequence in $p(C)$
which converges to a limit $z \in \bold C$.
There exists a sequence $(t_n)_{n \ge 1}$ in $\bold R$
such that $(z_n,t_n) \in C$ for all $n \ge 1$.
Upon multiplying each $(z_n,t_n)$ by $(0,t'_n)$
for appropriate $t'_n \in C \cap Z(H)$,
we can assume that the sequence $(t_n)_{n \ge 1}$ is bounded;
upon extracting a subsequence, we can assume that
the sequence $(t_n)_{n \ge 1}$ is convergent in $\bold R$.
It follows that the sequence $\big( (z_n,t_n) \big)_{n \ge 1}$
converges to some element $(z,t) \in C$
with $z = p(z,t) \in p(C)$.
\par

The last claim of (i) follows from the fact that,
whenever $C$ contains two non--commuting elements,
their commutator is in $\big(C \cap Z(H)\big) \smallsetminus \{e\}$.
\medskip

(ii) Let $C$ be in $\Cal C (H) \smallsetminus \Cal A(H)$, so
$p(C) \in \Cal C_{\bold Z^2, \bold R \oplus \bold Z, \bold C}
(\bold C)$.
Recall from Subsection 4.IV that the coarea of $p(C)$
is well--defined and lies in $[0,\infty[$.
We distinguish two cases.
\par

First case: $\operatorname{coarea}(p(C)) = 0$.
For every $\epsilon > 0$, there exist two non-commuting elements
$(z,t),(z',t') \in C$ such that
$0 < \vert \operatorname{Im}(\overline{z}z') \vert < \epsilon$.
The commutator of these two elements is 
$(0,\operatorname{Im}(\overline{z}z'))
\in C \cap Z(H)$. Since $C$ is closed, it follows that
 $Z(H) \subset C$.
\par

Second case: $\operatorname{coarea}(p(C)) > 0$.
Now  $p(C) \in \Cal L (\bold C)$, and
by Formula (5.2)  the commutator group $[C,C]$ 
is generated by $(0,\operatorname{coarea}(p(C)))$.
If $C \cap Z(H)$ is not the whole of $Z(H)\approx\bold R$
then, being closed, it
must be cyclic (hence discrete) and contain
 $[C,C]$ as a subgroup of finite index, $n$ say,
in which case $C \in \Cal L_n(H)$.

\medskip

(iii)  Let $(C_n)_{n \ge 1}$ be a sequence 
converging to some $C_0$
in $\Cal C (H) \smallsetminus \Cal A (H)$. 
By (i),  $p(C_n)$ is a closed subgroup of $\bold C$,
for all $n \ge 1$.
We have to show that $(p(C_n))_{n \ge 1}$
converges to $p(C_0)$.
\par

Let $\varphi : \bold N \longrightarrow \bold N$
be a strictly increasing map and let
$(y_{\varphi(n)})_{n \ge 1}$ be a sequence in $\bold C$
converging to some $y \in \bold C$,
with $y_{\varphi(n)} \in p(C_{\varphi(n)})$ for all $n \ge 1$.
For $n \ge 1$, choose $x_{\varphi(n)} \in C_{\varphi(n)}$
such that $p(x_{\varphi(n)}) = y_{\varphi(n)}$.
Upon multiplying $x_{\varphi(n)}$ by an appropriate element
in $C_{\varphi(n)} \cap Z(H)$, see Claim~(ii),
we can assume that the sequence $(x_{\varphi(n)})_{n \ge 1}$
is bounded in $H$, and therefore has a subsequence
which converges to some $x \in H$.
As $(C_n)_{n \ge 1}$ converges to $C_0$,
we have $x \in C$, and therefore $y = p(x) \in p(C)$.
Thus Condition (3.4a) from Section~3 applies.
\par
It is straightforward to check Condition (3.4b).

\medskip

(iv) We know from Remark~3.5.ii that  $\Cal L (H)$ is open in $\Cal C (H)$. 
(And we also know that $\Cal L (H)$ is non--empty;
more on this in Section~7.)
\par

Consider
$\Cal L_{!!} (H) = \{C \in \Cal C (H) \smallsetminus \Cal A (H)
\mid p(C) \in \Cal L (\bold C) \}$.
It follows from  (iii) and the openness
of $\Cal L (\bold C)$
 in $\Cal C (\bold C)$ that $\Cal L_{!!} (H)$
is open in $\Cal C (H) \smallsetminus \Cal A (H)$.
As the latter space is open in $\Cal C (H)$, by Proposition 3.4,
$\Cal L_{!!} (H)$ is also open in $\Cal C (H)$.

\medskip

(v) Let $C_n \longrightarrow C_0$
be as in the beginning of the proof of (iii).
We have to show that $\left( C_n \cap Z(H) \right)_{n \ge 1}$
converges to $C_0 \cap Z(H)$.
Since Condition (3.4.a) is straightforward to check,
we need only verify Condition (3.4.b).
Let $t \in C_0 \cap Z(H)$; we have to find
$t_k \in C_k \cap Z(H)$ so that
$\lim_{k \to \infty} t_k = t$.
\par

Choose $u,v \in C_0$ such that $[u,v] \ne e$.
For all $k \ge 1$, choose also $g_k, u_k, v_k \in C_k$
such that $\lim_{k \to \infty} g_k = t$,
$\lim_{k \to \infty} u_k = u$, and
$\lim_{k \to \infty} v_k = v$.
Observe that $\lim_{k \to \infty} p(g_k) = 0$;
also, $p(u)$ and $p(v)$ are $\bold R$--linearly independent
(since $[u,v] \ne e$), so $p(u_k), p(v_k)$
are $\bold R$--linearly independent for $k$ large enough.
Set $I = \{ k \ge 1 \mid g_k \notin C_k \cap Z(H) \}$.
In case $I$ is finite,
we can set $t_k = g_k$ if $k \notin I$
and choose $t_k$ arbitrarily if $t \in I$.
From now on, for convenience, we assume that $I$
is an infinite set of integers.
\par

As $p(g_k) \ne 0$ and $\lim_{k \in I, k \to \infty}p(g_k) = 0$,
it follows that the coarea of the lattice generated in $\bold C$
by $p(g_k), p(u_k), p(v_k)$ tends to $0$,
and therefore that $C_k \cap Z(H)$ is more and more dense in $Z(H)$
when $k \in I, k \to \infty$. Thus we can assure
that $\lim_{k \to \infty}t_k = t$ by defining
$t_k = g_k$ when $k \notin I$
and choosing $t_k \in C_k \cap Z(H)$
sufficiently close to $g_k$
when $k \in I$.
\hfill $\square$ 
\enddemo
\bigskip

Let $\Lambda$ be a non--abelian closed subgroup of $H$
which does not contain $Z(H)$.
It follows from Proposition 5.2 that $\Lambda$ is a lattice;
moreover, the commutator subgroup $[\Lambda,\Lambda]$ 
is of finite index in the centre 
$Z(\Lambda) = \Lambda \cap Z(H)$ of $\Lambda$.
For each integer $n \ge 1$, recall from Section~1 that
we denote by $\Cal L_n(H)$
the space of lattices $\Lambda$ such that
$[Z(\Lambda) : [\Lambda,\Lambda]] = n$.
We define
$$
p_n \, : \, \Cal L_n(H) \longrightarrow \Cal L (\bold C)
\tag{5.8}
$$
to be the projection induced by $p$.

\bigskip

\proclaim{5.3.~Proposition}
(i) For any $n \ge 1$, the subspace $\Cal L_n (H)$ 
is open in $\Cal C (H)$.

\smallskip

(ii) For  $(z,t),\,(z',t')\in H$, 
if $z,z'$ generate a lattice in $\bold C$,
then $(z,t),\,(z',t')$ generate a lattice in $H$.
\endproclaim

\demo{Proof} (i) The function $J : \Cal L_{!!} (H) \longrightarrow \bold R_+$
 defined by
 $$
 J(D) \hskip.2cm =  \hskip.2cm \left\{ \aligned
  \frac{1}{n} \hskip.2cm &\text{if} \hskip.2cm D \in \Cal L_n (H) 
 \\
 0 \hskip.2cm &\text{if} \hskip.2cm D \in \Cal L_{\infty} (H) 
\endaligned \right.
 $$
  is continuous because  
$$
J(D) \, = \, 
 \frac{ \min \{ \vert t \vert \mid t \in D \cap Z(H), t \ne e \} }
 {\operatorname{coarea}(p(D))}
 $$
 and each of the  functions
 $$
 \aligned
& \Cal L_{!!} (H) \longrightarrow \Cal C (Z(H)), \hskip.2cm
    D \longmapsto D \cap Z(H), 
\\
& \Cal C_{\bold Z, \bold R}(Z(H)) \longrightarrow \bold R_+, \hskip.2cm
   D \longmapsto \inf \{ \vert t \vert \mid t \in D, t \ne e \}
\\
& \Cal L (\bold C) \longrightarrow \bold R_+^*, \hskip.2cm
   L \longmapsto \operatorname{coarea} (L)
\endaligned
$$
is continuous.
\par

Thus $\Cal L_n (H)$ is open in $\Cal L (H)$
since it is the inverse image of the point $\frac{1}{n}$,
which is open in the image of $J$. 
As $\Cal L (H)$ is open in $\Cal C (H)$,
the space $\Cal L_n(H)$ is also open in $\Cal C (H)$.

\medskip

(ii) Denote by $\Gamma$
the subgroup of $H$ generated by two elements $a = (z, t)$ and $b = (z', t')$,
and set $\Lambda = p(\Gamma)$.
Since $H$ is nilpotent of class two, each element $\gamma\in\Gamma$
can be written as $\gamma=a^q b^r [a, b]^s$, where $p, q, r \in \bold Z$.

Assume that $\Lambda$ is a lattice in $\bold C$. 
Then $p(\gamma) =  qz + rz'=0$ if and only if $q=r=0$, 
and therefore $\Gamma \cap \ker(p)$ is the infinite cyclic group
generated by $[a, b] = \left( 0, \operatorname{Im}(\overline{z}z')\right)$.
Since $p(\Gamma) = \Lambda$ and 
$\Gamma \cap \ker(p) = \Gamma \cap Z(H)$ are lattices
in $\bold C$ and $Z(H)$ respectively, $\Gamma$ is a lattice in $H$.
\hfill $\square$ 
\enddemo

\bigskip

\head{\bf
6.~Closed subgroups of the  Heisenberg group which are not lattices
}\endhead
\medskip

Our analysis of $\Cal C(H)$ is spread over the next three sections.
In the present section we describe the subspace
formed by the closed subgroups that are not lattices.
In Subsections 6.I and 6.II, we deal with closed subgroups
which have Zariski closures of dimension $1$ and~$2$, respectively;
observe that, since $H$ is torsion--free,
the trivial group $\{e\}$ is the only subgroup
with Zariski closure of dimension $0$.
In  6.III, we consider
the subgroups which contain the centre $Z(H)$;
some, but not all, appear already in 6.I and 6.II.
In 6.IV, we describe the action 
of the group $\operatorname{Aut}(H)$
on the non--lattice part of $\Cal C (H)$.

\head{\bf
6.I.~Groups isomorphic to $\bold R$ and $\bold Z$
}\endhead
\medskip

{\bf (I.i) One--parameter subgroups in $H$.}
Any such subgroup is a real line through the origin
in $\bold C \times \bold R$,
$$
\left\{ (sz_0,st_0) \in H \mid s \in \bold R \right\}
\hskip.5cm \text{for some} \hskip.3cm
(z_0,t_0) \in H, \hskip.1cm (z_0,t_0) \ne (0,0).
$$
The centre $Z(H)$ corresponds to $z_0 = 0$.
In the Chabauty topology, these groups constitute a space
$$
\Cal C_{\bold R}(H) \approx \bold P^2
\tag{6.1}
$$
homeomorphic to a real projective plane.
\medskip

{\bf (I.ii) Infinite cyclic subgroups.}
Any $h \ne e$ in $H$ generates a closed subgroup isomorphic to $\bold Z$,
and each of these subgroups has exactly two generators.
Thus the subgroups isomorphic to $\bold Z$ constitute a space
$$
\Cal C_{\bold Z}(H) \, = \, (H \smallsetminus \{e\}) / \{\operatorname{id},J\}
\, \approx \, \bold P^2 \times ]0,\infty[
\tag{6.2}
$$
homeomorphic to the direct product of a real projective plane
and an open interval.
Here~$J$ denotes the  anti-automorphism of $H$
mapping each element $h = (z,t)$ to its inverse $h^{-1} = (-z,-t)$.
\medskip

{\bf (I.iii) The closure of the space of cyclic subgroups.}
The closure in $\Cal C (H)$
of   $\Cal C_{\bold Z}(H)$
is the space of subgroups that are
of  type (0),~(I.i) or (I.ii).
It is a closed cone on a projective plane
$$
\Cal C_{\{0\},\bold Z, \bold R}(H)
\, \approx \,
\big( \bold P^2 \times [0,\infty] \big)
\hskip.1cm / \hskip.1cm
\big((x,0) \sim (y,0)\big) ,
\tag{6.3}
$$
with the vertex $(*,0)$ of the cone corresponding
to  $\{e\}$, the
points in $\bold P^2 \times ]0,\infty[$ 
corresponding to
to infinite cyclic subgroups,
and the points in $\bold P^2 \times \{\infty\}$
to one-parameter subgroups.

\medskip
\head{\bf
6.II.~Other abelian groups: those isomorphic to
$\bold R^2$, $\bold R \oplus \bold Z$, and $\bold Z^2$
}\endhead
\medskip

{\bf (II.i) Maximal abelian subgroups.}
The multiplication formula (1.3) shows that the maximal abelian subgroups of
$H$ are of the form
$$
\left\{ (sz_0,t) \in H \mid (s,t) \in \bold R^2 \right\}
\hskip.5cm \text{for some} \hskip.3cm
z_0 \in \bold C^* ;
$$
they are all isomorphic to $\bold R^2$.
There is an obvious homeomorphism
from the space of such subgroups to the
space of real lines through the origin in $\bold C$,
$$
\Cal C_{\bold R^2}(H) \, \approx \, \bold P^1.
\tag{6.4}
$$
Note that every maximal abelian subgroup  of $H$ contains the
centre $Z(H) \cong \bold R$.
\par

It is convenient to fix
a left--invariant Riemannian metric on $H$
for which the submanifold
$\bold C \times \{0\}$ of $H$ 
(caveat~: it is {\it not} a subgroup of $H$!)
is orthogonal
to the centre $\{0\} \times \bold R = Z(H)$.
The induced Remannian metric on each
 $A \in \Cal C_{\bold R^2}(H)$
makes it
isometric to the Euclidean plane $\bold R^2$.
There are two isometric isomorphisms $A\to\bold C$
extending the isomorphism $Z(H)\to\bold R$
implicit in the notation $\{0\} \times \bold R = Z(H)$;
these correspond to a choice of
orientation on $A$.
We define
$\widehat{\Cal C}_{\bold R^2}(H)$ to be the set
of such groups $A$ {\it together with an orientation}
(equivalently, choice of isomorphism to $\bold C$).
There is an obvious bijection from
$\widehat{\Cal C}_{\bold R^2}(H)$ to
 the circle~$\bold S^1$
of unit vectors in the slice $\bold C\times\{0\}$
of $H = \bold C \times \bold R$.
Change of orientation map  defines a free
$(\bold Z /2 \bold Z)$--action on $\widehat{\Cal C}_{\bold R^2}(H)$,
and we have a two--sheeted covering
$$
\bold S^1 \approx \widehat{\Cal C}_{\bold R^2}(H)
\, \longrightarrow \,
\Cal C_{\bold R^2}(H) \approx \bold P^1
\tag{6.5}
$$
which, tracing through the
identifications, is the standard degree-two 
covering of $\bold P^1$ 
(see the Remark just before Theorem~1.3).

\medskip

\noindent{\it Notation:}
Any closed abelian subgroup $A$ of $H$ that is not contained in $Z(H)$
is contained in a unique closed subgroup of $H$ isomorphic to $\bold R^2$
and we shall  denote this $\overline{A}$.
When we place a \lq\lq hat\rq\rq on a letter, such as $\hat B$,
this will denote an oriented group, i.e.~an element of
$\widehat{\Cal C}_{\bold R^2}(H)$.

\medskip

{\bf (II.ii) Subgroups isomorphic to $\bold R \oplus \bold Z$.}
The space $\Cal C_{\bold R \oplus \bold Z}(H)$
is the total space of a fibration
$$
\Cal C_{\bold R \oplus \bold Z}(H)
\, \longrightarrow \,
\Cal C_{\bold R^2}(H)
\approx \bold P^1 ,
\hskip.2cm
A \longmapsto \overline{A} \hskip.3cm
\text{with fibre} \hskip.3cm
\Cal C_{\bold R \oplus \bold Z}(\bold C)
\approx \bold P^1 \times ]0,\infty[ .
\tag{6.6}
$$

To avoid confusion between the roles
played by the two copies of $\bold P^1$ in (6.6),
we will write $\bold P^1_Z$ for
the one corresponding to $\Cal C_{\bold R^2}(H)$ 
(maximal abelian subgroups of $H$,
containing $Z(H)$)
and write $\bold P^1_{\bold R}$ for the one
corresponding to the choice of the $\bold R$-factor
in $\Cal C_{\bold R \oplus \bold Z}(\bold C)$.
Observe that the circle $\bold S^1$ 
which appears a few lines below double covers $\bold P^1_Z$.
\par

Consider the two--fold cover 
$\widehat{\Cal C}_{\bold R \oplus \bold Z}(H):=
\{(A,\text{orientation on $\overline{A}$})\}$
of the total space of the previous fibration.
The corresponding fibration
$$
\widehat{\Cal C}_{\bold R \oplus \bold Z}(H) 
\longrightarrow \widehat{\Cal C}_{\bold R^2}(H)
$$
is  trivial. 
Indeed, an oriented plane $\hat B \subset H$ 
has a unique positive basis consisting of 
the unit central vector $(0,1)$
and a unit vector $(z_0,0)$ orthogonal to $Z(H)$;
the resulting canonical isomorphism from $\hat B$ onto $\bold C$
sends $(0,1) \in \hat B$ to $1 \in \bold C$
and $(z_0,0) \in \hat B$ to $i \in \bold C$.
Recall from Subsection 4.I that
$\Cal C_{\bold R \oplus \bold Z}(\bold C)$
is homeomorphic to $\bold P^1_{\bold R} \times ]0,\infty[$
(the open interval, $]0,1[$ in 4.I,
is better rescaled as $]0,\infty[$ here). 
It follows that
$$
\widehat{\Cal C}_{\bold R \oplus \bold Z}(H) \, = \, 
\Cal C_{\bold R \oplus \bold Z}(\bold C) \times 
\widehat{\Cal C}_{\bold R^2}(H) \, \approx \, 
\left( \bold P^1_{\bold R} \times ]0,\infty[ \right) \times 
\bold S^1 . 
\tag{6.7}
$$
\par

The total space $\Cal C_{\bold R \oplus \bold Z}(H)$ 
of the fibration (6.6) is the quotient 
of $\widehat{\Cal C}_{\bold R \oplus \bold Z}(H)$
by the orienta\-tion--reversing $\bold Z / 2 \bold Z$--action. 

Consider a closed subgroup
 $A\cong \bold R \oplus \bold Z$
 of $H$. Complex conjugation exchanges
 the maps $A\to\bold C$ corresponding to
the two orientations of $\overline{A}$.
Thus, if 
$(\theta,\rho,\varphi) \in 
\bold P^1_{\bold R} \times ]0,\infty[ \times \bold S^1 \approx
\widehat{\Cal C}_{\bold R \oplus \bold Z}(H)$
are the coordinates corresponding to one orientation,
the other orientation corresponds to coordinates
$(\sigma(\theta),\rho,\varphi+\pi)$,
where $\theta \longmapsto \sigma(\theta)$
denotes the symmetry of 
$\Cal C_{\bold R \oplus \bold Z}(\bold C) 
\approx \bold P^1_{\bold R}$ that replaces the chosen
$\bold R$-factor with its complex conjugate,
and the coordinate $\rho$ is the {\it inverse}
of the appropriate minimal norm (as 
described in Subsection 4.I).
Hence
$$
\Cal C_{\bold R \oplus \bold Z}(H) 
\, \approx \,
\big( \bold P^1_{\bold R} \times ]0,\infty[ \times \bold S^1 \big) 
\big/ 
\big( (\theta, \rho, \varphi) \, \sim \, 
(\sigma(\theta), \rho, \varphi+\pi)) 
\, \approx \,
\bold K \times ]0,\infty[ .
\tag{6.8}
$$
Consequently, (6.6) is the direct product with $]0,\infty[$
of the standard fibration of the Klein bottle $\bold K$
over the projective line $\bold P^1_Z$,
in other words of the non--trivial circle bundle over the circle.

\medskip

{\bf (II.iii) The closure of $\Cal C_{\bold R \oplus \bold Z}(H)$.} 
The frontier of 
$\Cal C_{\bold R \oplus \bold Z}(H)$ in $\Cal C (H)$
is the union of two pieces, namely
$\Cal C_{\bold R}(H)$ and $\Cal C_{\bold R^2}(H)$.
Convergence to points in the first piece occurs
when the $\rho$--coordinate tends to~$0$, while
convergence to points in 
the second piece arise when $\rho\to\infty$. 
Moreover,
all frontier points $A$ arise as limits of sequences
$(A_n)$ with $\overline A_n$ constant --- equal to
$\widetilde{A}$, say --- with only the $\rho$-coordinate
varying.
 
Given a closed subgroup $A$ of $H$ isomorphic to $\bold R$,
there are two cases to distinguish.
If $A \ne Z(H)$, then we define $\widetilde{A}  \cong \bold R^2$ to be
the subgroup generated by $A$ and $Z(H)$.
But if $A = Z(H)$ then every subgroup in 
$\Cal C_{\bold R^2}(H) \approx \bold P^1_Z$
is an equally strong candidate for $\widetilde{A}$. 
Hence the factor 
$\bold K$ in $\Cal C_{\bold R \oplus \bold Z}(H)$ 
can be obtained from
$\Cal C_{\bold R}(H) \approx \bold P^2$ by
blowing-up the point $Z(H)$ and making the
identification  $A \mapsto \widetilde{A}$ elsewhere. 

Thus $\Cal C_{\bold R}(H) \approx \bold P^2$
is attached to the end of 
$\Cal C_{\bold R \oplus \bold Z}(H) 
\approx \bold K \times ]0,\infty[$
near $0$ by means of the map
  $\bold K\to \bold P^2$ that
  blows-down the circle that appears
  in (6.8) as the image of $\{Z(H)\} \times \bold S^1
  \subset\bold P^1_{\bold R} \times \bold S^1$.
\par

The attachment of $\Cal C_{\bold R^2}(H)\approx
\bold P^1$
to the end of 
$\Cal C_{\bold R \oplus \bold Z}(H)$ where
$\rho\to\infty$ is straightforward:
$\Cal C_{\bold R \oplus \bold Z, \bold R^2}(H)$
is homeomorphic to the quotient of 
$\bold K \times ]0,\infty]$ 
by the relation that identifies
 $(k,\infty)$ to $(k',\infty)$
whenever $k,k' \in \bold K$ are in the same fibre of the natural fibration 
$\bold K \longrightarrow \bold (\bold S^1 /
(\varphi \sim \varphi+\pi))$.

\par
Thus $\Cal C_{\bold R, \bold R \oplus Z, \bold R^2}(H)$,
the closure of $\Cal C_{\bold R \oplus \bold Z}(H)$
in $\Cal C(H)$, is homeomorphic to 
the space obtained from $\bold K \times [0,\infty]$ 
by blowing down a circle in $\bold K \times \{0\}$  and by collapsing $\bold K \times \{\infty\}$ to
$\bold P^1$ in the manner described above.

\medskip

{\bf (II.iv) Subgroups isomorphic to $\bold Z^2$.}
As in (II.ii), the space $\Cal C_{\bold Z^2}(H)$
is the total space of a fibration
$$
\Cal C_{\bold Z^2}(H)
\, \longrightarrow \,
\Cal C_{\bold R^2}(H)
\approx \bold P^1 ,
\hskip.2cm
A \longmapsto \overline{A}
\hskip.3cm \text{with fibre} \hskip.3cm
\Cal L(\bold C)
\approx \bold C^2 \smallsetminus \Sigma ,
\tag{6.9}
$$
where $\Sigma$ denotes as in Section~4
the complex affine curve of equation $a^3 - 27 b^2 = 0$.
\par

For $A \in \Cal C_{\bold Z^2}(H)$, 
the two orientations of $\overline{A}$ correspond, as before, to conjugate embeddings $A\to\bold C$.
The values taken by
the classical invariants $g_2$ and $g_3$ on the
image of these
two embeddings are conjugate. Hence
$$
\widehat{\Cal C}_{\bold Z^2}(H)
\, = \,
\Cal L(\bold C) \times \widehat{\Cal C}_{\bold R^2}(H)
\, \approx \,
\left( \bold C^2 \smallsetminus \Sigma \right) \times \bold S^1 
\tag{6.10}
$$
and
$$
\Cal C_{\bold Z^2}(H) 
\, \approx \,
\left(\left( \bold C^2 \smallsetminus \Sigma \right) \times \bold S^1 \right) 
\big/ 
\big( 
(g_2,g_3,\varphi) \sim (\overline{g_2},\overline{g_3},\varphi+\pi)
\big) .
\tag{6.11} 
$$

\medskip

{\bf (II.v) Subgroups of rank $2$.}
We have a fibration
$$
\Cal C_{\bold Z^2, \bold R \oplus \bold Z, \bold R^2}(H)
\, \longrightarrow \,
\Cal C_{\bold R^2}(H) \approx \bold P ^1 ,
\hskip.2cm
A \longmapsto \overline{A}
\hskip.3cm \text{with fibre} \hskip.3cm
\Cal C_{\bold Z^2, \bold R \oplus \bold Z, \bold R^2}(\bold C) .
\tag{6.12}
$$
\par

Set
$
\widehat{\Cal C}_{\bold Z^2, \bold R \oplus \bold Z, \bold R^2}(H) 
=
\widehat{\Cal C}_{\bold Z^2}(H) \cup
\widehat{\Cal C}_{\bold R \oplus \bold Z}(H) \cup
\widehat{\Cal C}_{\bold R^2}(H)
$.
As in Subsections (II.ii) and (II.iv) above, the two--fold cover of (6.12)
$$
\widehat{\Cal C}_{\bold Z^2, \bold R \oplus \bold Z, \bold R^2}(H)
\, \longrightarrow \, 
\widehat{\Cal C}_{\bold R^2}(H)
$$
is a trivial fibration.
Recall from Chapter~4 that 
$\Cal C_{\bold Z^2, \bold R \oplus \bold Z, \bold R^2}(\bold C)$ 
is homeomorphic to the
complement of a $2$--disc in a $4$--sphere:
$$
\Cal C_{\bold Z^2, \bold R \oplus \bold Z, \bold R^2}(\bold C)
\simeq
(\bold C^2 \cup \{\infty\}) \smallsetminus \Sigma_{-},
$$
where 
$\Sigma_{-} = \{(a,b) \in \bold C^2 \mid
a^3 = 27 b^2 \hskip.2cm \text{and} \hskip.2cm |a|^2+|b|^2 \le 1 \}$  
is a $2$--disc $\bold D^2$
(recall that $\simeq$ denotes a homotopy equivalence).
Thus, 
$$
\aligned
\widehat{\Cal C}_{\bold Z^2, \bold R \oplus \bold Z, \bold R^2}(H) 
\, &= \, 
\Cal C_{\bold Z^2, \bold R \oplus \bold Z, \bold R^2}(\bold C) 
\times \widehat{\Cal C}_{\bold R^2}(H)
\\
\, &\approx \,
\left( (\bold C^2 \cup \{\infty\}) \smallsetminus \Sigma_{-} \right) \times \bold S^1
\, \approx \, 
(S^4 \smallsetminus D^2) \times \bold S^1,
\endaligned
\tag{6.13}
$$
where the embedding of  
$\bold D^2$ in $\bold S^4$
is not locally flat.
Reflection in the real axis in $\bold C$ 
(which corresponds to changing
the orientation of $\overline{A}$) acts on 
${\Cal C}(\bold C) \approx  \bold C^2 \cup \{\infty\}$ as
$(a,b) \longrightarrow (\overline{\vphantom{b}a},\overline{b})$. 
Thus the space 
${\Cal C}_{\bold Z^2, \bold R \oplus \bold Z, \bold R^2}(H)$ 
is homeomorphic to
$$
\Big(
\left( (\bold C^2 \cup \{\infty\}) \smallsetminus \Sigma_{-} \right) \times \bold S^1
\Big)
\,  \big/ \,
\Big(
(a,b,\varphi)\sim (\overline{\vphantom{b}a},\overline{b},\varphi+\pi)
\Big).
\tag{6.14}
$$
In particular, the space 
$\Cal C_{\bold Z^2, \bold R \oplus \bold Z, \bold R^2}(H)$ 
is  a topological manifold of dimension~$5$.
\par\noindent

\medskip

{\bf (II.vi) Closure of subgroups of rank $2$.}
The closure in $\Cal C (H)$  of  $\Cal C_{\bold Z^2}(H)$
is the whole space $\Cal A (H)$ of abelian closed subgroups.
More precisely: 

\bigskip

\proclaim{6.1.~Proposition} 
(i) 
The space $\Cal A (H)$ of abelian closed subgroups of $H$
is homeomorphic to
$$
\Cal A (H) \, \approx \,
\Big(S^4\times S^1\Big) \big/ 
\Big( (x,\varphi)\sim (x,\varphi'),  
\hskip.2cm \text{if $x \in I$} \Big),
$$
where $I\subset S^4$ is a tame closed interval.
\smallskip

(ii) Every $\operatorname{Aut}(H)$--orbit in 
$\Cal C_{\bold Z^2} (H)$ 
is dense in $\Cal A (H)$. 
\endproclaim

\demo{Proof} 
(i)
Let us introduce the two spaces
$$
\aligned
\widehat{\Cal A}_{\operatorname{en}}(H)
\, &\Doteq \, 
\{ (A, \hat B) \in {\Cal A}(H) \times \widehat{\Cal C}_{\bold R^2}(H) 
\mid
A \subset \hat B \} ,
\\ 
\Cal A_{\operatorname{en}}(H)
\, &\Doteq \, 
\{ (A,B) \in {\Cal A}(H) \times \Cal C_{\bold R^2}(H) 
\mid
A \subset B \} ,
\endaligned
$$ 
where \lq\lq en\rq\rq \ stands for \lq\lq enhanced\rq\rq .
There is an obvious two-sheeted covering
$\widehat{\Cal A}_{\operatorname{en}}(H) \longrightarrow
\Cal A_{\operatorname{en}}(H)$
and an equally obvious projection
$\Cal A_{\operatorname{en}}(H) \longrightarrow
\Cal A (H)$.
\par

For each
$\hat B_0 \in \widehat{\Cal C}_{\bold R^2}(H) \approx \bold S^1$,
we have a canonical identification of $\hat B_0$ with $\bold C$, as in II.ii, so
 the subspace of  $\widehat{\Cal A}_{\operatorname{en}}(H)$ consisting
of pairs $(A,\hat B_0)$ is canonically identified 
with $\Cal C(\bold C)$;
moreover, the change of orientation in 
$\widehat{\Cal C}_{\bold R^2}(H)$ corresponds to 
the involution $\sigma$ that
 complex conjugation induces on $\Cal C(\bold C)$.
 Thus
 $$
\Cal A_{\operatorname{en}}(H) 
\, \approx \,
\Big( \Cal C(\bold C)  \times \bold S^1\Big)
\Big/
\Big(
(C,\varphi) \sim (\sigma(C),\varphi+\pi)
\Big) .
$$
And by employing
the homeomorphism $f^{-1}$ of Theorem 4.6 we get:
$$
\Cal A_{\operatorname{en}}(H) 
\, \approx \,
\Big( (\bold C^2 \cup \{\infty\}) \times \bold S^1 \Big)
\Big/
\Big(
((a,b),\varphi) \sim (  (\overline{\vphantom{b}a},\overline{b}) , \varphi+\pi )
\Big) .
\tag{6.15}
$$
\par

For $\varphi \in \bold S^1$, let $\rho_{\varphi}$ be the homeomorphism
of $\bold C^2 \cup \{\infty\}$ defined by
$$
\rho_{\varphi}(a_1+ia_2 , b_1+ib_2) \, = \,
\left( a_1 + i(a_2 \cos \varphi - b_2 \sin \varphi) 
\hskip.1cm , \hskip.1cm  
b_1 + i(a_2 \sin \varphi + b_2 \cos \varphi) \right) .
$$
Observe that the complex conjugate of $\rho_{\varphi}(a,b)$
is $\rho_{\varphi + \pi}(a,b)$.
Define, then, a homeomorphism $R$ of 
$(\bold C^2 \cup \{\infty\}) \times \bold S^1$ by
$$
R((a,b),\varphi) \, = \, (\rho_{\varphi}(a,b),\varphi) .
$$
Writing $s : ((a,b),\varphi) \longmapsto 
((\overline{\vphantom{b}a},\overline{b}),\varphi + \pi)$
for the involution of
$(\bold C^2 \cup \{\infty\}) \times \bold S^1$
that appears in (6.15), we have
$$
R^{-1} s R ((a,b),\varphi) \, = \,
R^{-1} s (\rho_{\varphi}(a,b),\varphi) \, = \,
R^{-1} (\rho_{\varphi + \pi}(a,b) , \varphi + \pi) \, = \,
((a,b),\varphi+\pi) .
$$
Thus we obtain a homeomorphism
$$
\aligned
\Cal A_{\operatorname{en}}(H) 
\, &\approx \,
\Big(  (\bold C^2 \cup \{\infty\} ) \times \bold S^1 \Big)
\Big/ \{s\}
\\
\,&\approx \,
\Big(  (\bold C^2 \cup \{\infty\} ) \times \bold S^1 \Big)
\Big/ \{R^{-1}sR\}
\\
\, &= \,
\phantom{\Big(}
(\bold C^2 \cup \{\infty\} ) \times \bold P^1 .
\endaligned
\tag{6.16}
$$
The points in the right--hand term of (6.15)
corresponding to pairs $(A,B) \in \Cal A_{\operatorname{en}}(H)$
with $A \subset Z(H)$ are 
the points $((a,b),\varphi)$
with $a,b\in \bold R$, $a,b \ge 0$, $a^2+b^2 \le 1$,
and $a^3=27b^2$.
(This follows from Section~4, since $A$ is isomorphic
to one of $\{0\},\bold Z,\bold R$,
and the corresponding $C \in \Cal C (\bold C)$
is invariant under the action of $\sigma$.)
\par

Near each $(A,B) \in \Cal A_{\operatorname{en}}(H)$
with $A \not\subset Z(H)$, the projection
$\Cal A_{\operatorname{en}}(H) \longrightarrow \Cal A(H)$
is a local homeomorphism.
And if $A \subset Z(H)$ then all pairs 
$(A,B) \in \Cal A_{\operatorname{en}}(H)$
have the same image in $\Cal A (H)$.
It follows that
$$
\Cal A (H) \, \approx \,
\Big( (\bold C^2 \cup \{\infty\}) \times \bold P^1 \Big)
\Big/
\Big( ((a,b),\psi) \sim ((a,b),\psi ')
\hskip.2cm \text{for} \hskip.2cm 
(a,b) \in I , \hskip.2cm 
\psi,\psi' \in \bold P^1 \Big)
$$
where $I$ is the tame arc
$$
I \, = \,
\left\{ (a,b) \in \bold R^2 \, \mid \,
a,b \ge 0, \hskip.2cm a^2 + b^2 \le 1, \hskip.2cm
a^3=27b^2 \right\}
$$
of $(\bold C^2 \cup \{\infty\}) \times \bold P^1$.
This ends the proof of (i).

\medskip

(ii) For any closed subgroup $C_0$ of $H$,
let us denote by $\overline{\operatorname{Aut}(H)C_0}$
the closure in $\Cal C (H)$ of its $\operatorname{Aut}(H)$--orbit.
Let $\bold R \times \bold R$ denote the closed subgroup
$\{(x+iy,t) \in H \mid y = 0\}$ of $H$, and note
that its $\operatorname{Aut}(H)$-orbit is the
whole of  $\Cal C_{\bold R^2}(H)$.
(Indeed the subgroup $SO(2)$ has only one orbit
in $\Cal C_{\bold R^2}(H)$.) Thus,
in order to prove (ii),
it  suffices to show that
for any closed subgroup~$C_0$ of~$\bold R \times \bold R$ 
isomorphic to $\bold Z^2$,
the subspace 
$\{ C \in \overline{\operatorname{Aut}(H)C_0}
\mid
C \subset \bold R \times \bold R\}$
coincides with
$\{ A \in \Cal A (H) \mid A \subset \bold R \times \bold R \}$.
We will do this in three steps.

\smallskip

{\it Step one.}
Let $P$ be the group of automorphisms of $\bold R \times \bold R$ that extend to
automorphisms of $H$.
We regard
 $P$ as a subgroup of the full automorphism group
$GL(\bold R \times \bold R)$ of $\bold R \times \bold R$
(which {\it should not} be confused with the subgroup of
$\operatorname{Aut}(H)$ previously denoted $GL_2(\bold R)$).
As any automorphism of $H$ leaves the
centre invariant,
 $P$ consists entirely of matrices of the form
$\left( \matrix * & 0 \\ * & * \endmatrix \right)$.
Now, from the diagonal matrices
$\left( \matrix a & 0 \\ 0 & b \endmatrix \right)$
in the subgroup $GL_2(\bold R)$ of $\operatorname{Aut}(H)$,
we can obtain any diagonal matrix
$\left( \matrix a & 0 \\ 0 & ab \endmatrix \right)$ in $P$. On the other hand, using inner automorphisms
$\operatorname{Int}_{(-iv,*)} : (x,t) \longmapsto (x,t+vx)$,
we can obtain any unipotent matrix
$\left( \matrix 1 & 0 \\ v & 1 \endmatrix \right)$ in $P$.
It follows that $P$ is the full group of lower triangular matrices
in $GL(\bold R \times \bold R)$.

\smallskip

{\it Step two.} 
Let $G^+$ be 
the subgroup of $GL(\bold R \times \bold R)$
consisting of matrices with positive determinant,
and let $P^+=P \cap G^+$.
As $G^+$ can be seen as 
the set of oriented bases of $\bold R \times \bold R$,
there is a projection 
$G^+ \longrightarrow \Cal L(\bold R \times \bold R)$
mapping a matrix to the lattice generated by its columns;
this is a principal bundle whose structure group
is $SL_2(\bold Z)$,
acting on $G^+$ by right-multiplication.
The group $G^+$ itself acts on the right 
on the projective line $\bold P^1$ of $\bold R \times \bold R$,
by
$$
\left\{
\aligned
\bold P^1 \times G^+ 
\hskip1.2cm &\longrightarrow \hskip1.5cm
\bold P^1
\\
\left( 
\bold R (x,t) , \left( \matrix a & b \\ c & d \endmatrix \right)
\right)
\hskip.2cm &\longmapsto \hskip.2cm
\bold R (xa + tc , xb + td ) .
\endaligned
\right.
$$
This action is transitive
and the isotropy group of the second factor 
$Z(H)$ of $\bold R \times \bold R$ is $P^+$.
Thus  we may identify $P^+ \backslash G^+$ with $\bold  P^1$.
The double coset space $P^+ \backslash G^+ / SL_2(\bold Z)$
can be viewed equally
as the set of $P^+$--orbits in the set 
$G^+/SL_2(\bold Z) \approx \Cal L (\bold R \times \bold R)$
and as the set of $SL_2(\bold Z)$--orbits
in the projective line $P^+ \backslash G^+ \approx \bold P^1$.

\smallskip

{\it Step three.} It is well known 
that the natural action of $SL_2(\bold Z)$
on the projective line $\bold P^1$ is minimal,
i.e.~all its orbits are dense.
Here is an {\it ad hoc} argument.
Let $t$ be an arbitrary point in $\bold P^1 = \bold R \cup \{\infty\}$
and let $F = \overline{SL_2(\bold Z)t}$ denote its orbit closure.
Since $\left( \matrix 1 & n \\ 0 & 1 \endmatrix \right) t \to \infty$
when $n \to \infty$, we have $\infty \in F$.
For $\frac{a}{b} \in \bold Q$, with $a$ and $b$ in $\bold Z$ and coprime,
there exist $c,d \in \bold Z$ with $ad-bc = 1$, namely with
$\left( \matrix a & c \\ b & d \endmatrix \right) \in SL_2(\bold Z)$,
and $\frac{a}{b} = \left( \matrix a & c \\ b & d \endmatrix \right) \infty \in F$,
so that $\bold Q \subset F$ and therefore $F = \bold P^1$.
\par

It follows from Step two that 
the action of $P^+$ on $\Cal L (\bold R \times \bold R)$
is also minimal.
In particular, the orbit closure 
$\overline{P^+ C_0}$ contains 
the subspace $\Cal L (\bold R \times \bold R)$,
and therefore coincides with 
$\{ A \in \Cal A (H) \mid A \subset \bold R \times \bold R \}$. 
{\it A fortiori,} we have
$$
\{C \in \overline{\operatorname{Aut}(H)C_0}
\mid
C \subset \bold R \times \bold R \} 
\, = \,
\{ A \in \Cal A (H) \mid A \subset \bold R \times \bold R \} ,
$$
as was to be shown.
\hfill $\square$
\enddemo

\bigskip

\head{\bf
6.III.~The space $\cc(H)$ of subgroups that contain the centre
}\endhead
\medskip

Each closed subgroup  $C$ of $H$ that contains $Z(H)$ 
is uniquely determined by its image $p(C)$ in $\bold C$,
and the assignment $C \longmapsto p(C)$ gives a homeomorphism
$$
\cc(H) \, \approx \, \Cal C(\bold C)
\, \approx \, \bold S^4 ,
\tag{6.17}
$$
where the homeomorphism of $\Cal C(\bold C)$ with a $4$--sphere
is that of \cite{HubPo--79}. More precisely,
$\cc(H) \approx \bold S^4 \approx \Sigma \bold S^3$
is naturally the union of six subspaces:
\roster
\item"(i)"
the central subgroup $Z(H)$ corresponds to the south pole;
\item"(ii)"
the subgroups isomorphic to $\bold R^2$ 
correspond to a trefoil $T$ in $\bold S^3$; 
\item"(iii)"
the subgroups isomorphic to $\bold R \oplus \bold Z$ 
correspond to points on half--meridians between points of $T$
and the south pole; 
\item"(iv)"
the whole group $H$ corresponds to the north pole; 
\item"(v)"
the subgroups $C \in \Cal C (H)$ 
with $p(C) \cong \bold R \oplus \bold Z$
correspond to points on half--meridians between points of $T$
and the north pole;
\item"(vi)"
the subgroups in $\Cal L_{\infty}(H)$, 
namely the inverse images in $H$ of lattices in $\bold C$, 
correspond  to the  complement
of $\Sigma T$ in $\Sigma \bold S^3$,
which is open and dense in $\cc(H)$.
\endroster
We introduce here more notation:
\roster
\item"---" 
the union of the subspaces in (i), (ii)
and (iii) is a closed disc 
$\bold D_-$ inside $\cc (H)$,
\item"---"
the union of the subspaces in (ii), (iv) and (v)
is a closed disc $\bold D_+$ inside $\cc (H)$.
\endroster
For future reference (Section 8),
we record part of the previous discussion
as follows:

\bigskip

\proclaim{6.2.~Lemma} With the above notation, 
$$
\cc (H) \cap \Cal A(H) 
\hskip.2cm = \hskip.2cm
\bold D_- 
\hskip.2cm = \hskip.2cm \left\{ 
C \in \cc (H) \mid p(C) \in 
\Cal C_{\{0\},\bold Z,\bold R}(\bold C) 
\right\}.
$$
\endproclaim

\medskip
\head{\bf
6.IV.~$\operatorname{Aut}(H)$--orbits in
$\Cal A (H) \cup_{\bold D_-} \cc (H)$
}\endhead
\medskip

In the next proposition, the items are numbered in accordance with Theorem 1.4.

\bigskip

\proclaim{6.3.~Proposition}
For the natural action of the automorphism group
$\operatorname{Aut}(H)$ on $\Cal C(H)$, we have:
\roster
\item"(i)" 
the one--point orbit $\{e\}$;
\item"(ii)"
two $\operatorname{Aut}(H)$--orbits on
$\Cal C_{\bold R}(H) \approx \bold P^2$, 
see (6.1);
\item"(iii)"
two $\operatorname{Aut}(H)$--orbits on
$\Cal C_{\bold Z}(H) \approx \bold P^2 \times ]0,\infty[$,
see (6.2);
\item"(iv)"
one $\operatorname{Aut}(H)$--orbit on
$\Cal C_{\bold R^2}(H)  \approx  \bold P^1$,
see (6.4);
\item"(v)"
two $\operatorname{Aut}(H)$--orbits on
$\Cal C_{\bold R \oplus \bold Z}(H) 
\approx \bold K \times ]0,\infty[$,
see (6.6) and (6.8);
\item"(vi)"
uncountably many $\operatorname{Aut}(H)$--orbits on
$\Cal C_{\bold Z^2}(H)$,
see (6.9) and (6.11);
\item"(vii)"
one orbit in 
$p_*^{-1}\left( \Cal C_{\bold R \oplus \bold Z}(\bold C) \right)$;
\item"(viii)"
the one--point orbit $H$;
\item"(ix)"
and one orbit in $\Cal L_{\infty}(H)$.
\endroster
Moreover, for any point in $\Cal C (H)$,
its orbit under the identity component of $\operatorname{Aut}(H)$ 
coincides with its orbit under the whole of $\operatorname{Aut}(H)$.
\endproclaim

\demo{Comments}
The details of the proof are left to the reader but
we wish to highlight
the following points.
\par

In $\Cal C_{\bold R}(H)$, one orbit is a single point, 
the centre; 
the other orbit is its complement.
In $\Cal C_{\bold Z}(H)$,
one orbit consists of central subgroups, 
the other of non-central ones.
Of the two orbits
in $\Cal C_{\bold R \oplus \bold Z}(H)$,
one is $2$--dimensional and consists of the
subgroups whose identity component
is $Z(H)$; the other orbit is $3$--dimensional.
\par

Concerning the uncountably many orbits in
$\Cal C_{\bold Z^2}(H)$:
 the space $\Cal C_{\bold Z^2}(H)$ is 
 5-dimensional
whereas the quotient of $\operatorname{Aut}(H)$ by
the isotropy subgroup of a point in $\Cal C_{\bold Z^2}(H)$
is only 4-dimensional.
\par

The orbits of $\operatorname{Aut}(H)$ in $\cc(H)$ are in natural bijection
with the 
orbits of $GL_2(\bold R)$ in $\Cal C (\bold C)$.
In particular, there are six orbits in $\cc(H)$, namely
\roster
\item"---"
the centre $Z(H)$, which is one 
of the two $\operatorname{Aut}(H)$--orbits 
in $\Cal C_{\bold R}(H)$, 
\item"---"
the space $\Cal C_{\bold R^2}(H)$,
\item"---"
one of the two $\operatorname{Aut}(H)$--orbits 
in $\Cal C_{\bold R \oplus \bold Z}(H)$, 
\item"---"
the whole group $H$,
\item"---"
the space $p_*^{-1}\left( \Cal C_{\bold R \oplus \bold Z}(\bold C) \right)$,
\item"---"
and the space $\Cal L_{\infty}(H)$.
\endroster
The dimension of these orbits is 
$0, 1, 2, 0, 2, 4$, respectively.
\hfill $\square$
\enddemo
\bigskip

The preceding proposition accounts for all of
the
 $\operatorname{Aut}(H)$--orbits in $\Cal C (H)$
 except for those in $\Cal L_n(H)$, which
 will be described in Proposition 7.2.

\bigskip

\bigskip
\head{\bf
7.~The structure of the space of lattices in the real Heisenberg group
}\endhead
\medskip

We remind the reader that
the space of lattices $\Cal L(H)$ is a disjoint
union of the subspaces $\Cal L_n(H)$
defined by declaring that a lattice $\Lambda$ is in $\Cal L_n(H)$ 
if $[\Lambda,\Lambda]$ has index $n$ in $\Lambda \cap Z(H)$. In this section we shall describe the
structure of the spaces $\Cal L_n(H)$.
The map
$$
p_n \, : \, \Cal L_n(H) \longrightarrow \Cal L (\bold C)
\tag{7.1}
$$
from Section~5 will play a prominent role in 
our discussion, so we remind the reader that 
this is the projection induced by $p:H\to\bold C$.

\medskip

\proclaim{7.1.~Example}
For an integer $n \ge 1$, the subgroup 
$\Lambda_n$ of $H$
generated by 
$$
\Big(1,0\Big), \hskip.5cm \Big(i,0\Big) 
\hskip.5cm \text{and} \hskip.5cm  
\Big(0,\frac{1}{n}\Big)
$$
is a lattice in $H$, and $\Lambda_n \in \Cal L_n (H)$. 
Moreover, as a subset of $H = \bold C \times \bold R$,
$$
\Lambda_{n} \, = \, \bold Z [i] \times \frac{1}{n}\bold Z
\hskip.5cm \text{if $n$ is even}
\tag{7.2}
$$
and
$$
\Lambda_{n} \, = \,
\left\{ 
\left(x+iy,\frac{t}{2n}\right) \in 
\bold Z [i] \times \frac{1}{2n}\bold Z
\hskip.1cm \Big\vert \hskip.1cm
xy  \equiv t \pmod{2}
\right\}
\hskip.5cm \text{if $n$ is odd.}
\tag{7.3}
$$
In all cases ($n$ even and odd),
the lattice $\Lambda_n$ corresponds to
$$
\left( \matrix
1 & \bold Z & \frac{1}{n}\bold Z \\
0 & 1 & \bold Z \\
0 & 0 & 1
\endmatrix \right)
\, \subset \,
\left( \matrix
1 & \bold R & \bold R \\
0 & 1 & \bold R \\ 
0 & 0 & 1
\endmatrix \right)
$$
in the matrix picture.
\endproclaim

\demo{Proof}
Let us assume that $n$ is odd,
and let us show the equality in (7.3).
The verification of the other claims is left to the reader.
\par

Denote by $\widetilde \Lambda_n$ the right--hand term of (7.3).
By definition, and since both the commutator $[(1,0),(i,0)]$ 
and $(0,\frac{1}{n})$ are central, 
$\Lambda_n$ is the set of elements of the form
$$
\Big(1,0\Big)^x \Big(i,0\Big)^y \Big(0,\frac{1}{n}\Big)^s \, = \,
\Big(x+iy , \frac{xy}{2} + \frac{s}{n}\Big)
\hskip.5cm \text{with} \hskip.2cm x,y,s \in \bold Z .
\tag{7.4}
$$
Set $t \Doteq 2n\Big(\frac{xy}{2} + \frac{s}{n} \Big)$
and observe that $t = nxy + 2s \equiv xy \pmod{2}$,
so that the element (7.4) is in $\widetilde \Lambda_n$.
This shows that $\Lambda_n \subset  \widetilde \Lambda_n$.
\par

Conversely, if $\Big( x+iy, \frac{t}{2n} \Big) \in \widetilde \Lambda_n$,
set $s \Doteq \frac{1}{2}(t-nxy)$, and observe that $s \in \bold Z$,
so that
$$
\Big( x+iy,\frac{t}{2n}\Big) \, = \, 
\Big( 1,0\Big)^x \Big(i,0\Big)^y \Big( 0,\frac{1}{n} \Big)^s \in \Lambda_n .
$$
This shows that $\widetilde \Lambda_n \subset \Lambda_n$.
\hfill $\square$
\enddemo

\bigskip

\proclaim{7.2.~Proposition}
For every positive integer $n$, the space $\Cal L_n(H)$ 
is a homogeneous space for the group $\operatorname{Aut}(H)$.
\endproclaim

\demo{Proof}
Let $\Lambda \in \Cal L_n(H)$
and let $\Lambda_n$ be as in Example 7.1.
We have to find $g \in \operatorname{Aut}(H)$
such that $g(\Lambda) = \Lambda_n$.
\par

Since $GL_2(\bold R)$ acts transitively 
on the space of lattices in $\bold C$,
there exists $g \in GL_2(\bold R) \subset \operatorname{Aut}(H)$
such that $p_n(g\Lambda) = p_n(\Lambda_n) = \bold Z [i] \subset \bold C$.
Thus, without loss of generality,
we can assume that $p_n(\Lambda) = \bold Z [i]$,
namely that $\Lambda$ has three generators of the form
$$
\Big(1,s\Big), \hskip.3cm \Big(i,t\Big), \hskip.3cm \Big(0,\frac{1}{n}\Big) .
$$
The inner automorphism $\operatorname{Int}_{(-t+is,0)}$ maps
these three generators to
$$
\big(1,0\Big), \hskip.3cm \Big(i,0\Big), \hskip.3cm \Big(0,\frac{1}{n}\Big) ,
$$
and therefore $\operatorname{Int}_{(-t+is,0)} (\Lambda) = \Lambda_n$. 
\hfill $\square$
\enddemo

\bigskip

A pleasant property of $\Lambda_n$  is that
it  corresponds to a familiar lattice in the matrix picture. 
But we shall see below that  the lattice of the next example
serves as a more convenient basepoint
when one studies the isotropy of the 
$\operatorname{Aut}(H)$ action on $\Lambda_n$.
We'll say more about this  in Remark 7.6.

\bigskip

\proclaim{7.3.~Example}
For an integer $n \ge 1$, the subgroup 
$\Lambda'_n$ of $H = \bold C \times \bold R$
generated by 
$$
\Big(1,\frac{1}{2}\Big), \hskip.5cm \Big(i,\frac{1}{2}\Big) 
\hskip.5cm \text{and} \hskip.5cm  
\Big(0,\frac{1}{n}\Big)
$$
is a lattice in $H$, and $\Lambda'_n \in \Cal L_n (H)$. Moreover
$$
\Lambda'_n \, = \Lambda_n
\tag{7.5}
$$
if $n$ is even, and
$$
\aligned
\Lambda'_{n} 
\, &= \,
\left\{ 
\left(x+iy,\frac{t}{2n}\right) \in 
\bold Z [i] \times \frac{1}{2n}\bold Z
\hskip.1cm \Bigg\vert \hskip.1cm
\aligned
&t \hskip.2cm \text{even if} \hskip.2cm \frac{1}{2}(x+iy) \in \bold Z [i] ,
\\
&t \hskip.2cm \text{odd otherwise}
\endaligned
\right\}
\\
\, &\ne \Lambda_n 
\endaligned
\tag{7.6}
$$
if $n$ is odd.
\endproclaim

\demo{Proof}
The equality in (7.5) follows from the fact that,
if $n$ is even, $(0,\frac{1}{2})$ is a power of $(0,\frac{1}{n})$.
We assume from now on that $n$ is odd.
\par

Denote by $\widetilde \Lambda'_n$ the right--hand term of (7.6).
By definition, $\Lambda'_n$ is the set of elements of the form
$$
\Big( 1,\frac{1}{2} \Big)^x \Big( i,\frac{1}{2} \Big)^y \Big( 0,\frac{1}{n} \Big)^s
\, = \,
\Big( x+iy , \frac{x+y+xy}{2} + \frac{s}{n} \Big)
\hskip.5cm \text{with} \hskip.2cm  x,y,s \in \bold Z .
\tag{7.7}
$$
Set $t \Doteq 2n\Big( \frac{x+y+xy}{2} + \frac{s}{n} \Big)$
and observe that 
$$
t \, = \, n(x + y + xy) + 2s \, \equiv \, (x+1)(y+1)-1 \pmod{2} ,
$$
so that the element (7.7) is in $\widetilde \Lambda'_n$.
This shows that $\Lambda'_n \subset \widetilde \Lambda'_n$.
The opposite inclusion is proved similarly.
(Compare with the proof of Example 7.1.)
\par

To check that $\Lambda'_n \ne \Lambda_n$,
observe for example that
$\left(1, \frac{1}{2n} \right)$ and $\left(i, \frac{1}{2n} \right)$
are in $\Lambda'_n$ and not in~$\Lambda_n$,
while $(1,0)$ and $(i,0)$ lie in $\Lambda_n$ but not in $\Lambda'_n$.
\hfill $\square$
\enddemo

\bigskip

We will describe now the action of the automorphism group of $H$
on spaces of lattices.
\par

Consider $w = (u,v) \in \bold R^2$,
corresponding to an inner automorphism $\operatorname{Int}_{(w,*)}$ 
as in (5.4),
and let 
$g = \left( \matrix a & b \\ c & d \endmatrix \right) 
\in GL_2(\bold R)$
be viewed as an automorphism of $H$ as in (5.7).
We will denote by
$$
\Phi_{w,g} \, = \, \left( \operatorname{Int}_{(w,*)} \right) g
\, \in \, \operatorname{Aut}(H) 
\tag{7.8}
$$
the resulting composition.
From Proposition~5.1, we have
$$
\Phi_{w,g} \Phi_{w',g'} \, = \, \Phi_{w + g(w'), gg'} .
$$
For an integer $n \ge 1$, we set
$$
S'_n \, = \, \left\{
\Phi_{w,g} \in \operatorname{Aut}(H)
\hskip.2cm \Big\vert \hskip.2cm
w \in \big( \frac{1}{n} \bold Z \Big)^2
\hskip.2cm \text{and} \hskip.2cm
g \in GL_2(\bold Z)
\right\} .
\tag{7.9}
$$
We leave it to the reader to check that $S'_n$
is a subgroup of $\operatorname{Aut}(H)$
and that we have a natural split extension
$$
\{0\} 
\hskip.1cm \longrightarrow \hskip.1cm
\Big( \frac{1}{n} \bold Z \Big)^2
\hskip.1cm \longrightarrow \hskip.1cm
S'_n
\hskip.1cm \longrightarrow \hskip.1cm
GL_2(\bold Z)
\hskip.1cm \longrightarrow \hskip.1cm
\{1\} ,
\tag{7.10}
$$
where the homomorphisms are given by
$w \longmapsto \Phi_{w,\operatorname{id}}$
and $\Phi_{w,g} \longmapsto g$.

\bigskip

\proclaim{7.4.~Proposition}
For each integer $n \ge 1$, 
the subgroup $S'_n$ of $\operatorname{Aut}(H)$
is the stabilizer of the lattice $\Lambda'_n \in \Cal L_n(H)$.
\endproclaim

\demo{Proof}
Denote by $\operatorname{Stab}'_n$ 
the stabilizer of $\Lambda'_n$ in $\operatorname{Aut}(H)$;
we have to show that $\operatorname{Stab}'_n = S'_n$.
\par
Consider $\Phi_{w,g} \in \operatorname{Aut}(H)$,
as in (7.8).
Then $\Phi_{w,g} = \Phi_{w,\operatorname{id}} \Phi_{0,g}$.
The actions of $\Phi_{0,g}$
and $\Phi_{w,\operatorname{id}}$
on the first two generators of $\Lambda'_n$ are given by 
$$
\Phi_{0,g} 
\left(  1 , \frac{1}{2}  \right)
\, = \, 
\left( a +i c , \frac{1}{2} \right) ,
\hskip.5cm
\Phi_{0,g} 
\left(  i , \frac{1}{2}  \right)
\, = \, 
\left( b +i d , \frac{1}{2}  \right) ,
\tag{7.11}
$$
and, for $w = (u,v)$, 
$$
\aligned
\Phi_{(u,v),\operatorname{id}}
\left( 1,  \frac{1}{2} \right)
-
\left( 1,  \frac{1}{2} \right)
\, &= \, 
-  nv \left(  0 , \frac{1}{n} \right) ,
\\
\Phi_{(u,v),\operatorname{id}} 
\left( i,  \frac{1}{2} \right)
-
\left( i,  \frac{1}{2} \right)
\, &= \, 
\phantom{ - } nu \left(  0 , \frac{1}{n} \right) .
\endaligned
\tag{7.12}
$$
\par

Suppose first that $\Phi_{w,g} \in S'_n$.
Since $(a,c) \notin (2 \bold Z)^2$
and $(b,d) \notin (2 \bold Z)^2$,
the equalities of (7.11) 
and similar equalities for $\Phi_{0,g^{-1}}$ imply that 
$\Phi_{0,g} \in \operatorname{Stab}'_n$.
Since $n (-v,u) \in \bold Z^2$,
the equalities of (7.12) imply that
$\Phi_{w,\operatorname{id}} \in \operatorname{Stab}'_n$.
Hence $\Phi_{w,g} \in \operatorname{Stab}'_n$.
\par

Suppose now that $\Phi_{w,g} \in \operatorname{Stab}'_n$.
Then $g$ preserves $p_n(\Lambda'_n) = \bold Z [i]$,
so that $g \in GL_2(\bold Z)$,  and  $\Phi_{0,g}   \in S'_n$.
Moreover, using the previous step, we have
$\Phi_{0,g}    \in S'_n \subset \operatorname{Stab}'_n$, 
hence
$\Phi_{w,\operatorname{id}} = \Phi_{w,g} (\Phi_{0,g})^{-1}  
\in \operatorname{Stab}'_n$.
The equalities of (7.12) imply that 
$w \in \left( \frac{1}{n} \bold Z \right)^2$.
It follows that $\Phi_{g,w} \in S'_n$.
\hfill $\square$
\enddemo

\bigskip

\proclaim{7.5.~Remark}
Proposition 7.4 implies that the extensions (7.10)
corresponding to different values of $n$ are semi--direct products,
indeed they are all isomorphic to the obvious extension
$$
\{0\} \, \longrightarrow \,
\bold Z^2  \, \longrightarrow \,
\bold Z^2 \rtimes GL_2(\bold Z)  \, \longrightarrow \,
GL_2(\bold Z)  \, \longrightarrow \,
\{1\}
\tag{7.13}
$$
where the  action of $GL_2(\bold Z)$ on $\bold Z^2$
is the standard one.
\endproclaim

\demo{Proof} For each real number $\rho > 0$, the mapping
$$
\bold R^2 \rtimes GL_2(\bold R) 
\, \longrightarrow \,
\bold R^2 \rtimes GL_2(\bold R)
\hskip.2cm , \hskip.5cm
(w,g) \, \longmapsto \, (\rho w,g)
$$
is a group automorphism.
It follows that, for various values of $n$, the split extensions
$$
\{0\} 
\hskip.1cm \longrightarrow \hskip.1cm
\Big( \frac{1}{n} \bold Z \Big)^2
\hskip.1cm \longrightarrow \hskip.1cm
S'_n
\hskip.1cm \longrightarrow \hskip.1cm
GL_2(\bold Z)
\hskip.1cm \longrightarrow \hskip.1cm
\{1\} 
$$
are pairwise isomorphic.
\hfill $\square$
\enddemo

\bigskip

\proclaim{7.6.~Remark} For all $n \ge 1$, the sequence
$$
\{0\} 
\hskip.1cm \longrightarrow \hskip.1cm
\Big( \frac{1}{n} \bold Z \Big)^2
\hskip.1cm \longrightarrow \hskip.1cm
\operatorname{Stab}(\Lambda_n)
\hskip.1cm \longrightarrow \hskip.1cm
GL_2(\bold Z)
\hskip.1cm \longrightarrow \hskip.1cm
\{1\} 
$$
splits.
\endproclaim

\demo\nofrills{} 
But the splitting is more cumbersome in this case, 
which is why we prefer to work with~$\Lambda'_n$.
\enddemo

\bigskip

\proclaim{7.7.~Proposition}
(i)
For every integer $n \ge 1$, the projection
$$
p_n \, : \, \Cal L_n(H)
\, \longrightarrow \, 
\Cal L (\bold C)
$$
of (7.1) is a 
fibre bundle with fibre a $2$--torus.
\smallskip

(ii)
For two integers $n,n'$,
the two fibre bundles $p_n,p_{n'}$ are isomorphic
over the identity  $\Cal L (\bold C)\to \Cal L (\bold C)$.
\smallskip

In particular, for any pair $n,n'$ of positive integers,
the total spaces $\Cal L_n(H)$ and $\Cal L_{n'}(H)$
are homeomorphic to each other.
\endproclaim

\demo{Proof}
(i)
Given $L_0\in \Cal L(\bold C)$, 
we choose a positively oriented $\bold Z$-basis
$(z_0,z_0')$ for it.
Each lattice $L$ in a small neighbourhood 
$U \subset \Cal L(\bold C)$ of $L_0$
has a unique positively oriented $\bold Z$-basis 
$(z,z')$ close to $(z_0,z_0')$,
and the commutator in $H$ of elements 
of the form $(z,*)$ and $(z',*)$ is 
$(0,\operatorname{Im}(\overline z z'))
=
(0,\operatorname{area}(\bold C/L))$, 
which is close to $(0,\operatorname{area}(\bold C/L_0))$.
\par

For any choice of $r,r' \in [0,\frac{1}{n}[$, let
$\Lambda_{r,r'} \in \Cal L_n(H)$ be the lattice generated by
$$
(z,r\operatorname{Im}(\overline{z}z')) , \hskip.2cm
(z',r'\operatorname{Im}(\overline{z}z')) , \hskip.2cm
\text{and} \hskip.2cm
(0,\frac{1}{n}\operatorname{Im}(\overline{z}z')) .
$$
Then $p(\Lambda_{r,r'}) = L$ and 
$\Lambda_{r,r'} = \Lambda_{s,s'}$ if and only if 
$n(r-s,r'-s') \in \bold Z \times \bold Z$. 
Thus $(z,z',r,r')$ provide
coordinates for a trivialisation near $L_0$.

\medskip

(ii)
Consider the automorphism $\Psi_n$ of $\operatorname{Aut}(H)$
defined by
$$
\Psi_n \left( \Phi_{w,g} \right) \, = \,  \Phi_{\frac{w}{n}, g} 
\hskip.1cm .
$$
In Proposition 7.4 we proved that 
the stabilizer of $\Lambda_n'$ 
in $\operatorname{Aut}(H)$ is $S_n'$, 
and it is clear from the definition of $S_n'$ that
$$
\Psi(S'_1) \, = \, S'_n .
\tag{7.14}
$$
If 
$\Pi : \operatorname{Aut}(H) \longrightarrow GL_2(\bold R)$
is as in (5.3),  observe that 
$$
\Pi \Psi_n \, = \, \Pi .
\tag{7.15}
$$
Denote by $\pi'_n$ the restriction of $\Pi$ to $S'_n$.
The diagram
$$
\matrix
S'_n && \longrightarrow &&
\operatorname{Aut}(H) && \longrightarrow &&
\operatorname{Aut}(H) / S'_n = \Cal L_n (H)
\\ 
&&&&&&&&
\\
\downarrow \pi'_n  &&&& 
\downarrow  \Pi     &&&& 
\downarrow p_n
\\
&&&&&&&&
\\
GL_2(\bold Z) && \longrightarrow &&
GL_2(\bold R) && \longrightarrow &&
GL_2(\bold R) / GL_2(\bold Z) = \Cal L (\bold C)
\endmatrix
\tag{7.16}
$$
is commutative.
By (7.14), if we apply the homomorphism $\Psi_n^{-1}$
to $S'_n$ and $\operatorname{Aut}(H)$,
and the identity to the groups and space of the bottom row of (7.16),
we find the diagram corresponding to (7.16)
for the value $n=1$.
It follows that the bundles $p_n$ and $p_1$
are isomorphic over the identity of $\Cal L (\bold C)$.
\hfill $\square$
\enddemo

\bigskip

Let us now point out a natural factorisation
of the spaces $\Cal L_n(H)$.
The group of positive real numbers acts naturally
on $\Cal L (\bold C)$, see Section~4,
as well as on~$H$ by $(s,(z,t)) \longmapsto (sz,s^2t)$,
see Section~5.
The latter action induces free actions on
$\Cal L (H)$ and $\Cal L_{\infty}(H)$,
and we have homeomorphisms
$$
\aligned
\Cal L_n(H) \, &\approx \,
\Cal L_n^{\operatorname{umod}}(H) \times \bold R_+^*
\hskip.5cm \text{for all} \hskip.1cm n \ge 1 ,
\\
\Cal L_{\infty}(H) \, &\approx \,
\Cal L_{\infty}^{\operatorname{umod}}(H) \times \bold R_+^* ,
\endaligned
\tag{7.17}
$$
where the superscript \lq\lq umod\rq\rq \ indicates that
the lattice $p_n(\Lambda)$ in $\bold C$ is unimodular
for $\Lambda \in \Cal L_n^{\operatorname{umod}}(H)$,
and similarly for $\Lambda \in \Cal L_{\infty}^{\operatorname{umod}}(H)$.
The mapping $p_n$ introduced above,
as well as the homeomorphism 
$p_{\infty} : \Cal L_{\infty}(H) \longrightarrow \Cal L (\bold C)$,
have restrictions, again denoted by the same letters
$$
\aligned
p_n \, : \, \Cal L_n^{\operatorname{umod}}(H)
   \, &\longrightarrow \,
   \Cal L^{\operatorname{umod}} (\bold C)
\hskip.3cm \text{for all} \hskip.1cm n \ge 1 ,
\\
p_{\infty} \, : \, \Cal L_{\infty}^{\operatorname{umod}}(H)
   \, &\overset{\approx}\to{\longrightarrow} \,
   \Cal L^{\operatorname{umod}} (\bold C) .
\endaligned
$$
These mappings are $SL_2(\bold R)$--equivariant;
since $SL_2(\bold R)$ acts transitively on the space
$\Cal L^{\operatorname{umod}} (\bold C)$,
it follows that these mappings are also onto.

\medskip

Let $\operatorname{SAut}(H)$ 
denote the subgroup of $\operatorname{Aut}(H)$
of those automorphisms $\Phi$ which induce on $H / Z(H)$
an orientation preserving and area preserving automorphism.
Then $\operatorname{SAut}(H)$ is isomorphic 
to the standard semi--direct product $\bold R^2 \rtimes SL_2(\bold R)$.
For each integer $n \ge 1$,
the proofs of Propositions~7.2 \& 7.4
show that $\Cal L^{\operatorname{umod}}_{n}(H)$
is a homogeneous space $\operatorname{SAut}(H) / J_{n}$
and that we have a split exact sequence
$$
\{0\} 
\, \longrightarrow \,
\Big( \frac{1}{n} \bold Z \Big)^2 
\, \longrightarrow \,
J_{n} = \Big( \frac{1}{n} \bold Z \Big)^2 \rtimes SL_2(\bold Z)  
\, \longrightarrow \,
SL_2(\bold Z)  \, \longrightarrow \,
\{1\} ,
\tag{7.18}
$$
where $J_n = S'_n \cap \operatorname{SAut}(H)$.
Let us mention that 
the group $J_1 = \bold Z^2 \rtimes SL_2(\bold Z)$
appears in analytical number theory as the
{\it Jacobi group}; see Section 1.1 in \cite{EicZa--85}.
\par

As in Section~4, we denote by $\widetilde{SL}_2(\bold R)$ 
the universal covering of $SL_2(\bold R)$,
by $\widetilde{SL}_2(\bold Z)$ the inverse image 
of $SL_2(\bold Z)$ in this group,
and by $\bold R^2 \rtimes \widetilde{SL}_2(\bold R)$
the standard semi--direct product, 
which can be identified with
the universal covering group of $\operatorname{SAut}(H)$.
The space $\Cal L^{\operatorname{umod}}_{n}(H)$
is also diffeomorphic to the homogeneous space
$$
\left( \bold R^2 \rtimes \widetilde{SL}_2(\bold R) \right)
\, / \,
\left( \bold Z^2 \rtimes \widetilde{SL}_2(\bold Z) \right) .
$$
As  $\bold R^2 \rtimes \widetilde{SL}_2(\bold R)$  
is a contractible space, it follows that 
$\Cal L^{\operatorname{umod}}_{n}(H)$
is an Eilenberg--MacLane space with fundamental group
isomorphic to $\bold Z^2 \rtimes \widetilde{SL}_2(\bold Z)$.
In view of the factorisations (7.17)
and of the last claim of Proposition 7.7,
we have proved the following proposition.

\bigskip

\proclaim{7.8.~Proposition}
For any positive integer $n$,
the space $\Cal L_n(H)$ is a $6$--dimensional connected manifold
which is an Eilenberg--MacLane space with fundamental group
$\bold Z^2 \rtimes \widetilde{SL}_2(\bold Z)$.
\endproclaim

\bigskip

The proof of Claims (i) to (iii) in Theorem~1.3 is now complete~:
see the previous proposition for~(i), about $\Cal L (H)$,
Proposition~6.1 for~(ii), about $\Cal A (H)$, 
and Subsection~6.III for~(iii), about $\cc (H)$.

\bigskip
\head{\bf
8.~On the global structure of the space $\Cal C (H)$
}\endhead
\medskip

In this section we shall describe 
the closure of the space of lattices in $H$ 
and complete the proof of the theorems stated in the introduction. 
We begin with an analysis of the closure of the strata $\Cal L_n(H)$. 
For each $n$, the frontier $\partial \Cal L_n(H)$ is
as described in Theorem 1.4; 
it contains all of the non--lattice subgroups of $H$  except those 
in $\Cal L_{\infty}(H)$.
In Proposition~8.5, we show that $\Cal L_{\infty}(H)$
lies in the closure of the union of the $\Cal L_n(H)$, 
thus completing the proof that the lattices
form an open dense subset of $\Cal C(H)$. 
This also leads easily to a  proof 
that $\Cal C(H)$ is not locally connected at points of 
$\Cal L_{\infty}(H)$.
\par

Recall from the proof of Proposition~6.1 that
we denote by $\bold R \times \bold R$ the closed subgroup
$\{(x+iy,t) \in H \mid y = 0 \}$ of $H$.

\medskip
\head{\bf
8.I.~The frontier of $\Cal L_n (H)$, the frontier of $\bigcup_{n \ge 1} \Cal L_n (H)$,
and related matters
}\endhead
\medskip

Let first $n$ be a fixed positive integer.

We consider a sequence $\left( A_k \right)_{k \ge 1}$
of closed subgroups of $\bold R \times \bold R$,
a closed subgroup $A$ of $\bold R \times \bold R$,
and a sequence $\left( t_k \right)_{k \ge 1}$ of nonzero real numbers
with the following properties:
\roster
\item"---"
$A_k \cong \bold Z^2$ for all $k \ge 1$ and $A \cong \bold Z^2$;
\item"---"
$A_k \cap Z(H) = \langle (0,t_k) \rangle$ for all $k \ge 1$;
\item"---"
$A_k \to A$ and $\vert t_k \vert \to \infty$ for $k \to \infty$.
\endroster
(Recall that $\langle \cdots \rangle$ indicates 
the subgroup generated by $\cdots$.)
An example of such data is provided by
$A_k = \langle (1,0),(-\frac{1}{k},1) \rangle$,
$A = \bold Z^2$,  and $t_k = k$.

\bigskip

\proclaim{8.1.~Lemma}
Let $\left( A_k \right)_{k \ge 1}$,  $A$, and  $\left( t_k \right)_{k \ge 1}$
be as above.
\smallskip

(i) There exists a sequence $\left( \Lambda_k \right)_{k \ge 1}$
in $\Cal L_n (H)$ such that $\Lambda_k \cap (\bold R \times \bold R) = A_k$
for all $k \ge 1$.
\smallskip

(ii) For any  $\left( \Lambda_k \right)_{k \ge 1}$ as in (i),
we have $\lim_{k \to \infty} \Lambda_k = A$ in $\Cal C (H)$.
\endproclaim

\demo{Proof}
(i) Let $k \ge 1$. Since $A_k \cap Z(H) \cong \bold Z$,
the projection $p(A_k) \subset \bold R \subset \bold C$
is a closed subgroup isomorphic to $\bold Z$.
Choose $u_k \in \bold R^*$ such that $p(A_k) = \langle u_k \rangle$;
observe that $A_k$ is generated by $(u_k,s_k)$ and $(0,t_k)$
for some $s_k \in \bold R$.
Define $\Lambda_k$ to be the subgroup of $H$ 
generated by $A_k$ and $(i n \frac{t_k}{u_k},0)$.
Then $\Lambda_k \in \Cal L_n(H)$ and
$\Lambda_k \cap (\bold R \times \bold R) = A_k$
for all $k \ge 1$.
\par

(ii) Let $\left( \Lambda_k \right)_{k \ge 1}$ be as in (i).
Three observations are in order.
\par

Firstly, the group $p(A)$ is not reduced to $\{0\}$,
because $A \cong \bold Z^2$ is closed in $H$
(note that $p(A)$ need not be closed in $\bold C$);
it follows that $\sup_{k \ge 1} \vert u_k \vert < \infty$.
Secondly,
$\operatorname{coarea}(p(\Lambda_k)) = n\vert t_k \vert$
is arbitrarily large when $k \to \infty$. 
Finally
$$
\operatorname{Im}(z) \, \in \, \bold Z
\frac{\operatorname{coarea}(p(\Lambda_k))}{\vert u_k \vert}
\hskip.5cm \text{for all} \hskip.2cm
z \in p(\Lambda_k) , 
$$
so that
$$
\vert \operatorname{Im}(z) \vert \, \ge \,
\frac{\operatorname{coarea}(p(\Lambda_k))}{\vert u_k \vert}
\hskip.5cm \text{for all} \hskip.2cm
z \in p(\Lambda_k), \hskip.2cm z \notin \bold R .
$$
It follows from the first two observations 
that the right--hand term of the last inequality
tends to $\infty$ when $k$ tends to $\infty$.
\par

Hence 
$$
\lim_{k \to \infty} \Lambda_k \, = \,
\lim_{k \to \infty} \Lambda_k \cap p^{-1}(\bold R) \, = \,
\lim_{k \to \infty} A_k \, = \,
A,
$$
as claimed.
\hfill $\square$
\enddemo
\bigskip

Recall from Section~1 that 
$\Cal L_{!!} (H) = \left( \bigcup_{n=1}^{\infty} \Cal L_n (H) \right)
 \cup \Cal L_{\infty} (H)$,
from Lemma~6.2 that we have defined closed discs 
 $\bold D_-, \bold D_+ \subset \cc (H)$,
and that the interior of $\bold D_+$ is
$$
\overset{\circ} \to {\bold D}_+ 
\, = \,
 \left\{ C \in \cc (H) \mid 
 p(C) \in \Cal C_{\bold R \oplus \bold Z, \bold C}(\bold C) \right\}  
 \, = \,
 \cc (H) \smallsetminus \left( \Cal L_{\infty}(H) \cup \bold D_- \right).
 $$
 
 \bigskip

\proclaim{8.2.~Lemma}
 Let $(C_k)_{k \ge 1}$ be a sequence in $\Cal L_{!!} (H)$ 
 and $C \in \overset{\circ} \to {\bold D}_+$.
 \par
 Then $\lim_{k \to \infty} C_k = C$ in $\Cal C (H)$
 if and only if $\lim_{k \to \infty} p(C_k) = p(C)$ in 
 $\Cal C (\bold    C)$.
 \endproclaim
 
 \demo{Proof} 
 If $\lim_{k \to \infty} C_k = C$, 
 then $\lim_{k \to \infty} p(C_k) = p(C)$ by the continuity of $p_*$,
 see Proposition~5.2.
 \par
 
If $\lim_{k \to \infty} p(C_k) = p(C)$,
then $\lim_{k \to \infty} \operatorname{coarea}(p(C_k)) =
\operatorname{coarea}(p(C)) = 0$,
hence $\lim_{k \to \infty} (C_k \cap Z(H)) = Z(H)$,
and it follows that $\lim_{k \to \infty} C_k = C$.
\hfill $\square$
\enddemo

\bigskip

\proclaim{8.3.~Proposition} (i) For each $n \ge 1$,
the frontier of $\Cal L_n (H)$ in $\Cal C (H)$ is
$$
\Cal A (H) \, \cup_{\bold D_-} \, 
\left( \cc (H) \smallsetminus \Cal L_{\infty}(H) \right) .
$$
\par
(ii) The closure of $\Cal L_{\infty} (H)$ in $\Cal C (H)$ is
$\cc (H)$.
\endproclaim

\demo{Proof}
We leave the proof of (ii) to the reader, and we split that of (i) in three steps.
\par

{\it First step: the frontier of $\Cal L_n (H)$ contains $\Cal A (H)$.}
Lemma~8.1 shows that this frontier contains
at least one point of $\Cal C_{\bold Z^2}(H)$.
Since $\Cal L_n(H)$ is invariant 
under the action of $\operatorname{Aut}(H)$, 
its frontier is also invariant, and contains therefore 
at least one $\operatorname{Aut}(H)$--orbit
in $\Cal C_{\bold Z^2}(H)$.
It follows from Proposition~6.1.ii that
the frontier of $\Cal L_n(H)$ contains $\Cal A (H)$.
\par

{\it  Second step: the frontier of $\Cal L_n (H)$ contains 
$\overset{\circ} \to {\bold D}_+$.}
This is clear from  Lemma~8.2,
since $\Cal L (\bold C)$ is dense in $\Cal C (\bold C)$.
\par

{\it Third step.}
It remains to show that $\Cal L_{\infty} (H)$
is disjoint of $\overline{ \Cal L_n (H)}$.
This will be a consequence of Lemma 8.4.iii below.
\hfill $\square$
\enddemo
\bigskip

For an open subset $\Cal W$ in $\Cal L (\bold C)$
and a positive integer $N$, set
$$
\Cal U_{\Cal W,N} \, = \,
\Big\{ D \in  \Big( \bigcup_{n=N+1}^{\infty} \Cal L_n (H) \Big)
\cup \Cal L_{\infty} (H)
\hskip.2cm \Big\vert \hskip.2cm 
p(D) \in \Cal W \Big\} .
$$
  
\bigskip

\proclaim{8.4.~Lemma} (i) For $\Cal W$ and $N$ as above,
the subset $\Cal U_{\Cal W,N}$  is open in $\Cal C (H)$.
\smallskip
 
(ii) For $C \in \Cal L_{\infty} (H)$,
the set $\{ \Cal U_{\Cal W,N} \}$,
with $\Cal W \ni p(C)$ and $N \ge 1$,
is a basis of neighbourhoods of $C$ in $\Cal C (H)$.
\smallskip

(iii) Let $\left( C_k \right)_{k \ge 1}$ 
be a sequence in $\Cal L (H)$
and $C \in \Cal L_{\infty} (H)$; for $k \ge 1$, let $n_k$ be 
the integer such that $C_k \in \Cal L_{n_k}(H)$.
\par
Then $\lim_{k \to \infty} C_k = C$ if and only if
$\lim_{k \to \infty} p(C_k) = p(C)$ and $\lim_{k \to \infty} n_k = \infty$.
\endproclaim
 
\demo{Proof}
(i) Recall from Proposition~5.3 and its proof that the functions
$$
J \, : \, \Cal L_{!!} (H) \longrightarrow \bold R_+ , \hskip.5cm
J(C) \, = \, \left\{
\aligned 
\frac{1}{n} \hskip.2cm &\text{if} \hskip.2cm D \in \Cal L_n (H) 
\\
0 \hskip.2cm &\text{if} \hskip.2cm D \in \Cal L_{\infty} (H) 
\endaligned \right.
$$
and
$$
p_* \, : \, \Cal L_{!!} (H) \longrightarrow \Cal L (\bold C) , \hskip.5cm
D \longmapsto p(D)
$$
are continuous. Hence the subset
$$
\Cal U_{\Cal W,N} \, = \,
J^{-1} \left( \Big[0,\frac{1}{N}\Big[ \right) \, \cap p_*^{-1}(\Cal W) .
$$
 is open in $\Cal C (H)$.
\medskip

(ii) Consider the  function $\ell : H \longrightarrow \bold R_+$
defined by 
$\ell(z,t) = \vert z \vert + \vert t \vert$
and the function $\delta : H \times H \longrightarrow \bold R_+$ defined by
$$
\delta \big( (z,t) , (z',t') \big) \, = \,
\ell \big( (z,t)^{-1} (z',t') \big) \, = \,
\vert z'-z \vert + \vert t' - t - \frac{1}{2} \operatorname{Im}(\overline{z}z') \vert .
$$
The function $\delta$ {\it is not a distance function,}
because it does not satisfy the triangle inequality.
But we go on as if it was (compare with (3.3)); thus, we set
$$
B_R \, = \, \{ (z,t) \in H \hskip.2cm \vert \hskip.2cm \ell(z,t) < R \} 
\hskip.5cm \text{for all} \hskip.2cm R > 0 
$$
and we define for all $C \in \Cal C (H)$ and $R,\epsilon > 0$
$$
\Cal V_{R,\epsilon}(C) \, = \,
\left\{ 
D \in \Cal C (H) \hskip.2cm \Bigg\vert \hskip.2cm
\aligned
&\delta((z,t),D) < \epsilon \hskip.2cm \text{for all} \hskip.2cm
(z,t) \in C \cap \overline{B}_R \hskip.2cm \text{and} \hskip.2cm
\\
&\delta((z',t'),C) < \epsilon \hskip.2cm \text{for all} \hskip.2cm
(z',t') \in D \cap \overline{B}_R
\endaligned 
\right\} .
\tag{8.1}
$$
We leave it to the reader to check that
$\left( \Cal V_{R,\epsilon}(C) \right)_{R > 0, \epsilon > 0}$
is a basis of neighbourhoods of $C$ in $\Cal C (H)$,
all of them being relatively compact.
\par

Consider  $C \in \Cal L_{\infty}(H)$,
a pair $R,\epsilon > 0$,
and the resulting open neighbourhood $\Cal V_{R,\epsilon}(C)$ of $C$.
Set $\eta = \frac{\epsilon}{R + 2}$ and
$$
\aligned
\Cal W \, &= \, 
\left\{ L \in \Cal V_{R,\eta}(p(C)) \cap \Cal L (\bold C)
\hskip.2cm \vert \hskip.2cm
\operatorname{coarea}(L) < 2 \operatorname{coarea}(p(C)) \right\}
\\
N \, &= \,  
\left\lceil \frac{2 \operatorname{coarea}(p(C)) } {\eta } \right\rceil 
\endaligned
$$
(where $\lceil \cdots \rceil$ indicates the upper integral part).
To complete the proof, it suffices to show that
$\Cal U_{\Cal W,N} \subset \Cal V_{R,\epsilon}(C)$.
In other words, we choose $D \in \Cal U_{\Cal W,N}$
and we have to show the two conditions on $C$ and $D$
which appear in the right--hand side of (8.1).
\par

For $z,z' \in \bold C$ with $\vert z \vert < R$ and $\vert z' - z \vert < \eta$,
note that 
$$
\vert \operatorname{Im}(\overline{z}z') \vert \, = \, 
\vert \operatorname{Im}(\overline{z}(z'-z)) \vert \, < \, 
R \eta .
\tag{8.2}
$$
For $(z,t),(z',t') \in B_R \subset H$
with $\vert z'-z \vert < \eta$ and $\vert t'-t \vert < \eta$,
note that
$$
\delta\big( (z,t),(z',t') \big) \, \le \,
\vert z'-z \vert + \vert t'-t \vert 
+ \frac{1}{2} \vert \operatorname{Im}(\overline{z}z') \vert  \, < \,
2\eta + \frac{1}{2} R \eta \, < \, \epsilon .
\tag{8.3}
$$
\par

Consider now $(z',t') \in D$ with $\ell(z',t') \le R$.
Since $p(D) \in \Cal V_{R,\eta}(p(C))$
and since $z' \in p(D)$ satisfies $\vert z' \vert \le R$,
there exists $z \in p(C)$ with $\vert z - z' \vert < \eta$.
We have $(z,t') \in C$
and $\delta( (z',t'),(z,t') ) < \epsilon$ by (8.2).
This shows that
$\delta( (z',t'),C) < \epsilon$.
\par

Consider finally $(z,t) \in C$ with $\ell(z,t) \le R$.
Since $p(D) \in \Cal V_{R,\eta}(p(C))$
and since $z \in p(C)$ satisfies $\vert z \vert \le R$,
there exists $z' \in p(D)$ with $\vert z' - z \vert < \eta$.
Since $D \in \Cal U_{\Cal W,N}$,
there exists $u \in \bold R$ such that
$(0,u) \in D \cap Z(H)$ and
$$
\vert u \vert \, \le \,
\frac{ \operatorname{coarea}(p(D)) }{N} \, \le \,
\frac{ 2 \operatorname{coarea}(p(C))}{N} ,
$$
namely with $\vert u \vert < \eta$ (by the choice of $N$).
Hence there exists $t' \in \bold R$
with $(z',t') \in D$ and $\vert t' - t \vert < \eta$.
We have
$$
\delta( (z,t) , (z',t') ) \, = \,
\ell( z'-z , t'-t-\frac{1}{2} \operatorname{Im}(\overline{z}z') ) \, < \,
\vert z'-z \vert + \vert t' - t \vert + \frac{1}{2} R \eta \, < \,
\epsilon .
$$
This shows that 
$\delta( (z,t) ,  D ) < \epsilon$,
and thus completes the proof of (ii).

\medskip
(iii) This claim is a straightforward consequence of Claim~(ii).
\hfill $\square$
\enddemo

\bigskip

\proclaim{8.5.~Proposition}
(i) The  frontier of
$\Cal L (H) = \bigcup_{n=1}^{\infty} \Cal L_n (H)$
contains $\Cal L_{\infty} (H)$.
\par

(ii) The space $\Cal C (H)$ is not locally connected.
Indeed, any point $C \in \Cal L_{\infty}(H)$
does not have any connected neighbourhood
contained inside $\Cal L_{!!} (H)$.
\par

(iii) The subspaces $\Cal L (H)$ and $\Cal L_{!!} (H)$ 
are open and dense in $\Cal C (H)$.
\par

(iv) The space $\Cal C (H)$ is arc connected.
\endproclaim

\demo{Proof}
Claim (i) follows from Lemma 8.4.ii. 
\medskip

(ii) Let   $C \in \Cal L_{\infty} (H)$ 
and let $U$ be a neighbourhood of $C$ contained in $\Cal L_{!!} (H)$.
By Lemma~8.4, $U$ contains a neighbourhood $\Cal U_{\Cal W,N}$.
The function $J$,
from the proof of Proposition~5.3, 
is continuous on $U$,
takes its values in $\{1,\hdots,\frac{1}{n},\hdots,0\}$,
and is not constant on $U$
since it is not contant on $\Cal U_{\Cal W,N}$.
Hence $U$ is not connected,
and this establishes Claim~(ii).
 \medskip
 
 (iii)
 The spaces $\Cal L (H)$ and $\Cal L_{!!} (H)$ are open,
 see Remark~3.5.ii and Proposition~5.2.v respectively.
 The closure of $\Cal L (H)$ is the whole of $\Cal C (H)$
 by Proposition~8.3 and by Claim~(i) above.
 
\medskip

(iv) First note that the space $\Cal D (H)$
of discrete subgroups of $H$ is arc--connected and contains $\{e\}$.
Indeed,  for  $s \in \bold R$, $s > 0$, denote by $\varphi_s$
the automorphism $(z,t) \longmapsto (sz,s^2t)$ of $H$.
Then, for each discrete subgroup $\Gamma$ of $H$,
we have $\lim_{s \to \infty} \varphi_s(\Gamma) = \{e\}$ in~$\Cal D (H)$.
If $\Lambda$ is a lattice then
we have also $\lim_{s \to 0} \varphi_s(\Lambda) = H$, 
so that $H$ lies in the path component of $\{e\}$.
\par

The path component of $\{e\}$ also contains $\Cal A(H)$,
see Proposition 6.1, while  the path component of $H$ contains
$\cc(H)$ since the latter is homeomorphic to $\bold S^4$. 
As $\Cal C(H)$ is the union of $\Cal D (H), \, \Cal A(H)$,
and $\cc(H)$, this completes the proof of (iv).
\hfill $\square$
\enddemo

\bigskip

\head{\bf
8.II.~The space $\Cal L_{!!} (H)$ as a bundle over $\Cal L (\bold C)$,
and the homotopy type of $\Cal C (H)$
}\endhead
\medskip

Let us denote by $\bold L$ the subset 
$\{1, \frac{1}{2}, \hdots, \frac{1}{n}, \hdots, 0\}$ 
of the real line,
endowed with the induced topology from $\bold R$,
which makes it a countable compact space.
For any topological space $X$, denote by $\bold L X$ 
the quotient of $X \times \bold L$ 
by the equivalence relation 
$(x,0) \sim (x',0)$, for all $x,x' \in X$
(we will resist writing more than once that $\bold L X$ is the 
\lq\lq discrete cone\rq\rq \ over~$X$).
Observe that if $X$ is compact, so is $\bold L X$.
For a fibre bundle  $p : E \longrightarrow B$, 
with fibre $F$, let 
$\bold L_{\operatorname{bu}}(p) : 
\bold L_{\operatorname{bu}} E 
\longrightarrow B$
be the fibre bundle, 
with fibre $\bold L F$,
of which the total space 
is the quotient of  $E \times \bold L$
by the equivalence relation
$(e,0) \sim (e',0)$ for $e,e' \in E$ with  $p(e) = p(e')$;
we will denote by $p(e)$ the class of $(e,0)$.

\bigskip

\proclaim{8.6.~Proposition} The projection 
$p : \Cal L_{!!}(H) \longrightarrow \Cal L (\bold C)$
is a  fibre bundle which is isomorphic to
$\bold L_{\operatorname{bu}}(p_1) : 
\bold L_{\operatorname{bu}} \Cal L_1 (H) 
\longrightarrow \Cal L (\bold C)$.
\endproclaim

\demo{Proof}
For $n \ge 1$, denote by 
$\Theta_n : \Cal L_n (H) \longrightarrow \Cal L_1 (H)$
the bundle isomorphism of Proposition~7.7.
Define $\Theta : \Cal L_{!!} (H) \longrightarrow \bold
L_{\operatorname{bu}}\Cal L_1 (H)$ by 
$$
\Theta(C) \, = \, \left\{ \aligned
\left( \Theta_n (C) , \frac{1}{n} \right) 
\hskip.2cm &\text{if} \hskip.2cm
C \in \Cal L_n (H)
\\
p(C) \hskip1cm &\text{if} \hskip.2cm C \in \Cal L_{\infty} (H) . 
\endaligned \right.
$$
Then $\Theta$ is clearly a bijection 
above the identity on $\Cal L (\bold C)$.
As its domain and its range are fibre bundles with compact fibres 
(which are of the form $\bold L T$, with $T$ a $2$--torus) 
over a locally compact basis (the space $\Cal L (\bold C)$),
it only remains to check that $\Theta$ is continuous at each point.
\par

The continuity of $\Theta$ at a point $C \in \Cal L_n (H)$ 
follows from the continuity of $\Theta_n$ and from the fact that
$\Cal L_n (H)$ is open in $\Cal L_{!!} (H)$.
The continuity of $\Theta$ at a point $C \in \Cal L_{\infty} (H)$
follows from Lemma~8.4.ii.
\hfill $\square$
\enddemo
\bigskip

Let $X$ be a topological space, $B$ an open dense subspace,
and $p : E \longrightarrow B$ a fibre bundle.
We define the {\it fibre--collapse continuation of $p$ over $X$}
as the space $E_* = E \sqcup (X \smallsetminus B)$
with the projection $p_* : E_* \longrightarrow B$
defined by $p_*(e) = p(e)$ for $e \in E$
and $p_*(x) = x$ for $x \in X \smallsetminus B$;
the topology on $E_*$ is defined by decreeing that $E$ is open, 
and $x \in X \smallsetminus B$ has a basis of open neighbourhoods
consisting of the sets $p_*^{-1}(U)$, 
with $U$ an open neighbourhood of $x$ in $X$.
With this definition:

\bigskip

\proclaim{8.7.~Proposition}
The bundle projection 
$p : \Cal L_{!!}(H) \longrightarrow \Cal L (C)$
extends to the continuous mapping
$$
p_* \, : \,
\Cal C (H) \smallsetminus \Cal A (H) = \Cal L_{!!}(H) \cup 
\overset{\circ} \to {\bold D}_+
\hskip.5cm  \longrightarrow \hskip.5cm
\Cal L (\bold C) \cup 
\Cal C_{\bold R \oplus \bold Z, \bold C}(\bold C)
$$
of Proposition 5.2.iii, which is isomorphic to
the fibre--collapse continuation
of $p$ over $\Cal L (\bold C) \cup 
\Cal C_{\bold R \oplus \bold Z, \bold C}(\bold C)$.
In particular, the restriction of $p_*$ to 
$\overset{\circ} \to {\bold D}_+$ is a homeomorphism onto
\par\noindent
$\Cal C_{\bold R \oplus \bold Z, \bold C}(\bold C)$.
\endproclaim

\bigskip

One might hope that $p_*$ 
has an extension to $\Cal C (H)$ 
showing that $\cc(H)$ is a retract of $\Cal C (H)$.
The next proposition shows that this optimism is too naive.
But all is not lost, as Proposition 8.9 shows that
the desired retraction does exist at the level of homotopy.

\bigskip

\proclaim{8.8.~Proposition}
The mapping $p_*$ of Proposition 8.7 does not have a continuous extension $\Cal C (H) \longrightarrow \Cal C (\bold C)$.
\endproclaim

\demo{Proof}
If such an extension
$\overline{p}_* : \Cal C (H) \longrightarrow \Cal C (\bold C)$
were to exist, it would be unique, since $\Cal A (H)$ is nowhere dense.
Moreover, $\overline{p}_*$ would be equivariant 
with respect to the homomorphism
$\Pi : \operatorname{Aut}(H) \longrightarrow GL_2(\bold R)$
of (5.3), since this is the case for $p_*$.
\par

Let $A_k = \langle (1,0), (-\frac{1}{k},1) \rangle$, $k \ge 1$,
and $A = \bold Z^2$ be as just before Lemma~8.1, 
and let $(\Lambda_k)_{k \ge 1}$ be
the corresponding sequence of lattices; see Lemma~8.1.ii. 
Then $\lim_{k \to \infty} p(\Lambda_k) = \bold R \subset \bold C$,
and therefore $\overline{p}_* (\bold Z^2) = \bold R$.
\par

Let $(\varphi_s)_{s > 0}$ be the automorphisms of $H$ from the proof of Proposition~8.5. 
Then $\lim_{s \to \infty} \varphi_s (\bold Z^2) = \{e\}$ and
$r\left(\varphi_s\right) (\bold R) = \bold R$ for all $s > 0$.
Hence $\overline{p}_*(\{e\}) = \bold R$.
But this is impossible, since the action of 
$\operatorname{Aut}(H)$ on $\{e\}$ is trivial
and the action of $GL_2(\bold R)$ on $\bold R$ is not.
\hfill $\square$
\enddemo
\bigskip

We  denote by 
$$
q_* \, : \,
\Cal C (H) \smallsetminus \Cal A (H) = \Cal L_{!!}(H) \cup 
\overset{\circ} \to {\bold D}_+
\hskip.5cm  \longrightarrow \hskip.5cm
\Cal L_{\infty} (H) \cup \overset{\circ} \to {\bold D}_+
$$
the composition of $p_*$ with the natural homeomorphism
from $\Cal L (\bold C) \cup 
\Cal C_{\bold R \oplus \bold Z, \bold C}(\bold C)$
onto 
$\Cal L_{\infty} (H) \cup \overset{\circ} \to {\bold D}_+$.
\par

Let $\Cal C (H) / \Cal A (H)$ denote the compact space
obtained from $\Cal C (H)$ by identifying all points
of $\Cal A (H)$ with each other. 
Let similarly $\cc (H) / \bold D_-$ denote the compact space
obtained from $\cc (H)$ by identifying all points
of $\bold D_-$ with each other;
observe that it is the union of
$\Cal L_{\infty} (H) \cup \overset{\circ} \to {\bold D}_+
= \cc (H) \smallsetminus \bold D_-$
and the point $[\bold D_-]$ (the class of $\bold D_-$).
\par

We denote by 
$$
q_*^{\operatorname{ext}} \, : \, 
\Cal C (H) / \Cal A (H) \longrightarrow
\cc (H) / \bold D_-
$$
the extension of $q_*$ mapping the point $[\Cal A (H)]$
to the point $[\bold D_-]$.

\bigskip

\proclaim{8.9.~Proposition}
(i) With the notation above, $q_*^{\operatorname{ext}}$
is a retraction.
\par

(ii) There exists a continuous map
$\psi : \cc (H) / \bold D_- \longrightarrow \cc (H)$
such that the composition
$\cc (H) \longrightarrow \cc (H) / \bold D_-
\overset{\psi}\to{\longrightarrow} \cc (H)$
is homotopic to the identity.
\par

(iii) The composition
$$
\Cal C (H) \hskip.2cm \longrightarrow \hskip.2cm 
\Cal C (H) / \Cal A (H) \hskip.2cm  
\overset{q_*^{\operatorname{ext}}}\to{\longrightarrow} \hskip.2cm 
\cc (H) / \bold D_- 
\hskip.2cm \overset{\psi}\to{\longrightarrow} \hskip.2cm  \cc (H)
\tag{8.4}
$$
is a weak retraction.

(iv) The fourth homotopy group $\pi_4(\Cal C (H))$
has an infinite cyclic quotient.
\endproclaim

\demo{Proof}
(i) For $C \in \Cal C (H) \smallsetminus \Cal A (H)$,
we have $q_*^{\operatorname{ext}}(C) = \langle C,Z(H) \rangle$;
this implies that the restriction of $q_*^{\operatorname{ext}}$
to $\cc (H) / \bold D_-$ is the identity. 
\par

Note that
$q_*^{\operatorname{ext}}(C) = p^{-1}\big( p(C) \big)$
for any $C \in \Cal C (H) \smallsetminus \Cal A (H)$.
We have to check that 
$q_*^{\operatorname{ext}}$
is continuous at each point.
The continuity of $q_*^{\operatorname{ext}}$ at points
distinct from $[\Cal A (H)]$ results from the continuity of
$C \longmapsto p(C)$ at these points (Proposition~5.2)
and the continuity of $p^{-1}$ as a mapping from
closed subsets of $\bold C$ to closed subsets of $H$.
It remains to check that $q_*^{\operatorname{ext}}$
is continuous at $[\Cal A (H)]$.
\par

Let $\left( C_k \right)_{k \ge 1}$ be a sequence
in $\Cal C (H) \smallsetminus \Cal A (H)$
such that 
$$
\lim_{k \to \infty} C_k \, = \,  [\Cal A (H)] \, \in \,
\Cal C (H) / \Cal A (H) .
$$
As the projection $p_* : \Cal C (H) \smallsetminus \Cal A (H)
\longrightarrow \Cal C (\bold C) \smallsetminus 
\Cal C_{\{0\},\bold Z,\bold R}(\bold C)$
is proper, see Proposition~8.7,
it follows that 
$\lim_{k \to \infty} q_*^{\operatorname{ext}} (C_k)
= [\bold D_-]$.

\medskip

(ii) In subsection 6.III we saw that there is a homeomorphism
from $\Cal C_{\ge Z}(H)$ to the 4-sphere 
that identifies $\bold D_-$ with 
a subset of a closed hemisphere $\bold S_-$. 
Thus the identity map of $\Cal C_{\ge Z}(H)$
induces a continuous map
$\Cal C_{\ge Z}(H) / \bold D_- \to \Cal C_{\ge Z}(H) / \bold S_-$, 
and we can define $\psi$ to be the composition of this map and
a homotopy inverse to  
$\Cal C_{\ge Z}(H) \to \Cal C_{\ge Z}(H) / \bold S_-$.
\medskip

(iii) The mapping of (8.4) is continuous
as a composition of continuous maps,
and its restriction to $\cc (H)$ is homotopic to the identity
by (i) and (ii).

\medskip

(iv) This is a straightforward consequence of (iii),
and of the homeomorphisms of $\cc (H)$ with $\bold S^4$.
\hfill $\square$
\enddemo


\bigskip



\bigskip
\Refs
\widestnumber\no{ReVaW--02}

\ref \no BenPe--92 \by R. Benedetti and C. Petronio
\book Lectures on hyperbolic geometry
\publ Springer \yr 1992
\endref

\ref \no BorJi--06 \by A. Borel and L. Ji
\book Compactifications of symmetric and locally symmetric spaces
\publ Birkh\"auser \yr 2006
\endref

\ref \no Bourb--63 \by N. Bourbaki
\book Int\'egration, chapitres 7 et 8
\publ Hermann \yr 1963
\endref


\ref \no CaEpG--87 \by R.D. Canary, D.B.A. Epstein, and P.L. Green
\paper Notes on notes of Thurston
\jour in \lq\lq Analytical and geometrical aspects of hyperbolic
spaces\rq\rq , D.B.A. Epstein Editor,
London Math. Soc. Lecture Notes Series {\bf 111} 
\yr Cambridge Univ. Press, 1987 \pages 3--92
\endref

\ref \no Chaba--50 \by C. Chabauty
\paper Limite d'ensembles et g\'eom\'etrie des nombres
\jour Bull. Soc. Math. France \vol 78 \yr 1950 \pages 143--151
\endref

\ref \no ChaGu--05 \by C. Champetier and V. Guirardel
\paper Limit groups as limits of free groups:
compactifying the set of free groups
\jour Israel J. Math. \vol 146 \yr 2005 \pages 1--76
\endref

\ref \no Champ--00 \by C. Champetier
\paper L'espace des groupes de type fini
\jour Topology \vol 39 \yr 2000 \pages 657--680
\endref


\ref \no CoGuP--07 \by Y. de Cornulier, L. Guyot, and W. Pitsch
\paper On the isolated points in the space of groups
\jour J. of Algebra \vol 307 \yr 2007 \pages 254--277
\endref

\ref \no CoGuP \by Y. de Cornulier, L. Guyot, and W. Pitsch
\paper The space of subgroups of an abelian group
\jour In preparation
\endref

\ref \no EicZa--85 \by M. Eichler and D. Zagier
\book The theory of Jacobi forms
\publ Birkh\"auser \yr 1985
\endref

\ref \no Fell--62 \by J.M.G. Fell
\paper A Hausdorff topology for the closed subsets 
of a locally compact non Hausdorff space
\jour Proc. Amer. Math. Soc. \vol 13 \yr 1962 \pages 472--476
\endref

\ref \no Ghys--07 \by E. Ghys
\paper Knots and dynamics
\jour Proceedings ICM, Madrid, 2006, Volume I
Plenary lectures and ceremonies, European Math. Soc.
\yr 2007 \pages 247--277
\endref

\ref \no Grigo--84 \by R. Grigorchuk
\paper Degrees of growth of finitely generated groups and the theory of
invariant means
\jour Math. USSR Izv. \vol 25 \yr 1985 \pages 259--300
\endref

\ref \no Gromo--81 \by M. Gromov
\paper Groups of polynomial growth and expanding maps
\jour Publ. Inst. Hautes \'Etudes Scient. \vol 53 \yr 1981 \pages 53--73
\endref

\ref \no GuiR\'e--06 \by Y. Guivarc'h and B. R\'emy
\paper Group--theoretic compactification of Bruhat--Tits buildings
\jour Ann. Scient. \'Ec. Norm. Sup. \vol 39 \yr 2006 \pages 871--920
\endref


\ref \no GuJiT--96 \by Y. Guivarc'h, L.~Ji, and J.C. Taylor
\book Compactifications of symmetric spaces
\publ Birkh\"auser \yr 1998
\endref

\ref \no Harpe--00 \by P. de la Harpe
\book Topics in geometric group theory
\publ The University of Chicago Press \yr 2000
\endref

\ref \no Harve--77 \by W.J. Harvey
\paper Spaces of discrete groups
\jour in \lq\lq Discrete groups and automorphic functions\rq\rq ,
W.J. Harvey Editor, Academic Press \yr 1977 \pages 295--438
\endref

\ref \no HurCo--64 \by A. Hurwitz and R. Courant
\book Vorlesungen \"uber allgemeine Funktionentheorie und elliptische
Funktionen
\publ 4th Edition, Springer \yr 1964
\endref

\ref \no HubPo--79 \by J. Hubbard and I. Pourezza
\paper The space of closed subgroups of $\bold R^2$
\jour Topology  \vol 18  \yr 1979 \pages 143--146
\endref

\ref \no Kapla--71 \by I. Kaplansky
\book Lie algebras and locally compact groups
\publ The University of Chicago Press \yr 1971
\endref

\ref \no MacSw--60 \by A.M. Macbeath and S. Swierczkowski
\paper Limits of lattices in a compactly generated group
\jour Canadian J. Math. \vol 12 \yr 1960 \pages 427--437
\endref

\ref \no Milno--71 \by J. Milnor
\book Introduction to algebraic K--theory
\publ Princeton Univ. Press \yr 1971
\endref

\ref \no Munkr--75 \by J.R. Munkres
\book Topology: a first course
\publ Prentice Hall \yr 1975
\endref
 

\ref \no Mumfo--76 \by D. Mumford
\book Algebraic geometry I, complex projective varieties
\publ Springer \yr 1976
\endref

\ref \no Oler--73 \by N. Oler
\paper Spaces of closed subgroups of a connected Lie group
\jour Glasgow Math. J. \vol 14 \yr 1973 \pages 77--79
\endref

\ref \no Raghu--72 \by M.S. Raghunathan
\book Discrete subgroups of Lie groups
\publ Springer \yr 1972
\endref

\ref \no SacZy--65 \by S. Saks and A. Zygmund
\book Analytic functions
\publ Second Edition, Polish Scientific Publishers \yr 1965
\endref

\ref \no Thurs--80 \by W.P. Thurston
\book The geometry and topology of three--manifolds
\publ Lectures Notes, Princeton University, 1980
\yr www.msri.org/publications/books/gt3m
\endref

\endRefs
\enddocument